%% file: roitner2.tex
\documentclass[12pt,reqno,a4wide]{amsart}

\usepackage{tikz} 
\usepackage{calrsfs}
\usepackage{mathrsfs}
\usetikzlibrary{decorations.pathmorphing}
\usetikzlibrary{arrows.meta}
\usepgflibrary {shadings}
\usepackage{array}
\usepackage{xcolor}
\usepackage{comment}
\newcommand{\seqnum}[1]{\href{https://oeis.org/#1}{\rm\underline{#1}}}

\usepackage{hyperref}

\usetikzlibrary{decorations.pathreplacing}
\usetikzlibrary{fit}

\usepackage{pgfplots}

\allowdisplaybreaks

\newcommand{\mD}{\mathcal{D}}

\catcode`,\active

\catcode`\,12

\begin{document}
\bibliographystyle{plain}

%
%

	\title[Motzkin paths of bounded height]
	{ Motzkin paths of bounded height with two   forbidden contiguous subwords of length two}

	\author[H. Prodinger ]{Helmut Prodinger }
	\address{Department of Mathematics, University of Stellenbosch 7602, Stellenbosch, South Africa
	and
NITheCS (National Institute for
Theoretical and Computational Sciences), South Africa.}
	\email{hproding@sun.ac.za}


	\begin{abstract}
		Motzkin excursions and meanders are revisited.  This is considered in the context of forbidden patterns.
		Previous work by Asinowski, Banderier, Gittenberger, and Roitner is continued. Motzkin paths of bounded height are
		considered, leading to matrix equations and also to continued fractions.
		The enumeration is done by properly setting up bivariate generating functions which can be expanded using
		the kernel method.

	\end{abstract}
	

\maketitle

\tableofcontents

\section{Introduction}

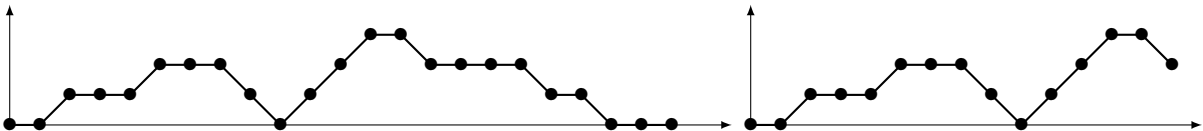
\begin{figure}[h]
	\begin{center}
		\begin{tikzpicture}[scale=0.4 ]
			
			\draw[thin,-latex] (0,0) -- (24,0);
			\draw[thin,-latex] (0,0) -- (0,4);		
			\draw [thick](0,0)-- (1,0) -- (2,1) --(3,1) --(4,1) --(5,2) --(6,2) --(7,2) --
			(8,1) --(9,0) --(10,1) --(11,2) --(12,3) --(13,3) --(14,2) -- (15,2) --(16,2) --(17,2) --(18,1) --
			(19,1) --(20,0) --(21,0) --(22,0) ;
			\node at (0,0) {$\bullet$};
			\node at (1,0) {$\bullet$};
			\node at (2,1) {$\bullet$};
			\node at (3,1) {$\bullet$};
			\node at (4,1) {$\bullet$};
			\node at (5,2) {$\bullet$};
			\node at (6,2) {$\bullet$};
			\node at (7,2) {$\bullet$};
			\node at (8,1) {$\bullet$};
			\node at (9,0) {$\bullet$};
			\node at (10,1) {$\bullet$};
			\node at (11,2) {$\bullet$};
			\node at (12,3) {$\bullet$};
			\node at (13,3) {$\bullet$};
			\node at (14,2) {$\bullet$};
			\node at (15,2) {$\bullet$};
			\node at (16,2) {$\bullet$};
			\node at (17,2) {$\bullet$};
			\node at (18,1) {$\bullet$};
			\node at (19,1) {$\bullet$};
			\node at (20,0) {$\bullet$};
			\node at (21,0) {$\bullet$};
			\node at (22,0) {$\bullet$};
		\end{tikzpicture}%
		\begin{tikzpicture}[scale=0.4 ]
			\draw[thin,-latex] (0,0) -- (15,0);
			\draw[thin,-latex] (0,0) -- (0,4);		
			\draw [thick](0,0)-- (1,0) -- (2,1) --(3,1) --(4,1) --(5,2) --(6,2) --(7,2) --
			(8,1) --(9,0) --(10,1) --(11,2) --(12,3) --(13,3) --(14,2);
			\node at (0,0) {$\bullet$};
			\node at (1,0) {$\bullet$};
			\node at (2,1) {$\bullet$};
			\node at (3,1) {$\bullet$};
			\node at (4,1) {$\bullet$};
			\node at (5,2) {$\bullet$};
			\node at (6,2) {$\bullet$};
			\node at (7,2) {$\bullet$};
			\node at (8,1) {$\bullet$};
			\node at (9,0) {$\bullet$};
			\node at (10,1) {$\bullet$};
			\node at (11,2) {$\bullet$};
			\node at (12,3) {$\bullet$};
			\node at (13,3) {$\bullet$};
			\node at (14,2) {$\bullet$};
		\end{tikzpicture}
		
	\end{center}
	\caption{A Motzkin excursion and a Motzkin meander (ending at level $2$).}
	\label{motz1}
\end{figure}

Motzkin excursions (Motzkin paths) consist of steps $U=(1,1)$, $D=(1,-1)$, $H=(1,0)$, start at the origin, never go below the
$x$-axis and return to the $x$-axis. Motzkin meanders (partial Motzkin paths, prefixes of Motzkin paths) are similar, but must not
necessarily return to the $x$-axis. The latter $H$ stands for `horizontal', sometimes one finds the symbol
$F$ (`flat'). The example excursion (Figure~\ref{motz1}) does not contain the (contiguous) pattern $UD$.
The notations `excursion' and 'meander' were Philippe Flajolet's favorites, see e. g. \cite{BF}; in his honor, they will be used here.

This paper deals with \emph{forbidden subpatterns}. In our analysis, we typically have \emph{two} contiguous patterns that are forbidden.
Our effort is a continuation of earlier work by Asinowski, Banderier, Gittenberger, Roitner (and perhaps others), \cite{ABBG, AsyWally,Vally}.
What we did:
\begin{itemize}
	\item[$\bullet$] We rederived appropriate generating function.
	\item[$\bullet$] We worked out  (almost?) all interesting cases explicitly.
	\item[$\bullet$] We considered  versions of \emph{bounded height}; the example excursion (Figure~\ref{motz1}) has height 3, which is 
	largest ordinate value that occurred.
	\item[$\bullet$] We obtained generating functions for (restricted) Motzkin meanders ending on a prescribed level.
		\item[$\bullet$] We provided continued fraction expansions for the generating functions of Motzkin excursions.
		Note that, since the generating functions are always algebraic of order 2, there \emph{must} be a periodic expansion.
		However, the expansions we obtained are motivated by the bounded height model and thus have a combinatorial meaning.
\item Unlike other papers where one often reads `other cases may be done in a similar fashion,' here, these `other cases' are really investigated.
\end{itemize}

What we did \emph{not} do, although it would be possible:
\begin{itemize}
	\item[$\bullet$] Derive generating functions of bounded restricted Motzkin meanders with prescribed end level.
	\item[$\bullet$] Consider questions of an asymptotic type, e.~g., the average height of restricted Motzkin excursions.
	\item[$\bullet$] Allowing the paths to go below the $x$-axis.
\end{itemize}

For the readers' convenience, we prepared tables of the cases that were considered. 
\begin{center}
	\begin{tabular}{ | l || c | c | c | c| l |}
		\hline
		forbidden & $UU, HH$&$UD, HH$&$DU, HH$&$DD, HH$      \\ \hline
		excursion &  A329666&A329676&A329666&A329666 \\		\hline
		meander &  A329667& A329675&A329668  & A329669  \\		\hline
	\end{tabular}
\end{center}

\begin{center}
	\begin{tabular}{ | l || c | c | c | c| l |}
		\hline
		forbidden & $UH, HU$&$DH, HD$&$UH, HD$&$DH, HU$      \\ \hline
		excursion &  A329701&A329701&A329702&A329701 \\		\hline
		meander & --- & ---& --- & ---  \\		\hline
	\end{tabular}
\end{center}

\begin{center}
	\begin{tabular}{ | l || c | c | c | c|c| l |}
		\hline
		forbidden & $UU, DD$&$UU, DU$ & $UU, HD$&$UU, UD$\\ \hline
		excursion &  A004149&A023431 &A217282&A023431\\		\hline
		meander & A308435 & ---&---      &---\\		\hline
	\end{tabular}
\end{center}

\begin{center}
	\begin{tabular}{ | l || c | c | c | c| l |}
		\hline
		forbidden & $UD, DD$&$DD, DU$ &$UD, DU$ &$UD,DH$\\ \hline
		excursion &  A023431&A023431&A004149 &A023431\\		\hline
		meander & --- & --- &  A308435&---   \\		\hline
	\end{tabular}
\end{center}

\begin{center}
	\begin{tabular}{ | l || c | c | c | c| l |}
		\hline
		forbidden & $DU, HD$&$DU,UH$  &$HU,DD$ & $UD,HU$\\ \hline
		excursion &  A004149&A217282 &A217282& ---   \\		\hline
		meander & A308435 & --- &  --- & ---   \\		\hline
	\end{tabular}
\end{center}

\begin{center}
	\begin{tabular}{ | l || c | c | c |c|  l |}
		\hline
		forbidden &$UU, HU$ & $DU, DH$&$DU, HU$&$UD, UH$&$UD,HD$ \\ \hline
		excursion &degenerate  &degenerate&degenerate&degenerate&degenerate  \\		\hline
	\end{tabular}
\end{center}

\begin{center}
	\begin{tabular}{ | l || c | c | c |c|  l |}
		\hline
		forbidden &$HU, HD$ & $HU, HH$ \\ \hline
		excursion &degenerate  &degenerate  \\		\hline
	\end{tabular}
\end{center}

Note that, even for sequences that were in OEIS, the description might not exactly match our
application.

Should some reasonable patterns have been overlooked (or duplicated), corrections will occur
in later versions of this draft.\footnote{Readers are encouraged to report such instances.}

The motivation to provide such a quasi-complete list originates from Drew Sills' amazing rendition
of Lucy Slater's account on Rogers-Ramanujan (type) identities, \cite{Drew}. 

\section{Outline of the approach}

\begin{itemize}
	\item All of our examples consist of three states; since we take care of the current altitude of the path as well,
	we rather talk about layers. One layer is responsible for $U$-steps, leading into it, another one for $H$-steps, 
	and another one for $D$-steps.
	
	\item We consistently use the letters $f,g,h$ (and $F,G,H$) for the three layers; for each case, a nice drawing is provided.
	There are recursions for $f_i=f_i(z)$, $g_i=g_i(z)$, $h_i=h_i(z)$, according the possible steps leading into a layer.
	
	\item Double (bivariate) generating functions $F(u)=F(u,z)=\sum_{i\ge0}f_iu^i$ etc. are introduced. The recursions
	translate into three functional equations for $F(u),G(u),H(u)$. They can be solved easily, but contain the quantities
	$F(0),G(0),H(0)$. Substitution of $u=0$ is not sufficient, but in each instance, there is a universal denominator 
	which factorizes as $(u-r_1(z))(u-r_2(z))$. One of these factors is always `bad,' in the sense that $1/(u-r_1(z))$ does not have
	a power series expansion around $(0,0)$. However, for combinatorial reasons, this cannot be, so the bad factor must have a bad
	counterpart in the numerator, and cancellation must occur. Our notation is consistent, so that $u-r_1(z)$ is the bad factor.
	
	\item Once the bad factors disappear, setting $u=0$ is possible, and thus $F(0),G(0),H(0)$ can be computed. 
	This procedure is typical for the kernel method, and many attractive examples have been collected in \cite{Prodinger-kernel, garden}.
	The total is
	always denoted by $S(u)=F(u)+G(u)+H(u)$, and $[u^j]S(u)$ is the generating function of restricted Motzkin meanders ending 
	at level $j$. So we get $S(0)$ (excursions), $S(1)$ meanders, and, as a bonus $[u^j]S(u)$.

\item Since there is the denominator, which is quadratic in $u$, we obtain recursions of order 2 for $f_j, g_j, h_j, s_j$. In order to introduce the
concept of \emph{bounded height}, we replace $f_j, g_j, h_j, s_j$ by 0 for $j>K$. The resulting recursion can always be written as a linear system,
involving a $(K+1)\times (K+1)$ matrix. We solve for $f_0, g_0, h_0, s_0$ and use the notations $\phi_{K+1}, \gamma_{K+1}, \eta_{K+1}, \sigma_{K+1}$ for the results.
This is always done by Cramer's rule, and expresses the solutions as a quotient of two determinants.

\item After a warm-up example, we discuss the issue that the layer that describes  a symbol $D$ leading into it (say it is the third layer), contains a state 
$h_K$ but it cannot be reached, since $f_{K+1},g_{K+1},h_{K+1},s_{K+1},$ are \emph{taboo}. So, strictly speaking, the bounded version of $s_j$ is slightly incorrect, 
as it allows contributions that shouldn't be allowed. So, $s^{[K+1]}=f^{[K+1]}+g^{[K+1]}+h^{[K+1]}$ (ad hoc notation) is computed, but correct would be
$f^{[K+1]}+g^{[K+1]}+h^{[K]}$. However, each of the 3 constituants are computed individually.

\item Since $(r_1/r_2)^K$ goes to 0 for large $K$, the bounded versions converge to the unbounded versions. Convergence is in terms of power series, that
more and more coefficients are correct.

\item The $(K+1)\times (K+1)$ matrix has an irregular first row and column. Otherwise, it is a  band matrix, and the determinant satisfies a recursion of order 2 (and thus allows for a Binet form). This recursion may be rewritten in terms of continued fractions, not uncommon when dealing with objects of bounded height.
\end{itemize}




\input{DHHD}
\input{UDDU}

\input{DDDU}

\input{UHHU}

\input{UUHH}

\input{UDHH}

\input{DUHH}

\input{DDHH}

\input{UHHD}

\input{DHHU}

\input{UUDD}
\input{UUDU}

\input{UUHD}

\input{UUUD}

\input{UDDD}

\input{UDHU}

\input{DUHD}

\input{DUUH}

\input{UDHH}
\input{HUDD}


\noindent

(Concerned with sequences  
\seqnum{A329666},
\seqnum{A329676},
\seqnum{A329667},
\seqnum{A329675},
\seqnum{A329668},
\seqnum{A329669},
\seqnum{A329701},
\seqnum{A329702},
\seqnum{A004149},
\seqnum{A023431},
\seqnum{A217282},
\seqnum{A023431},
\seqnum{A308435}.%
)

\bibliographystyle{plain}

\end{document}

%% file: DHHD.tex
\section{DH and HD are forbidden, A329701}
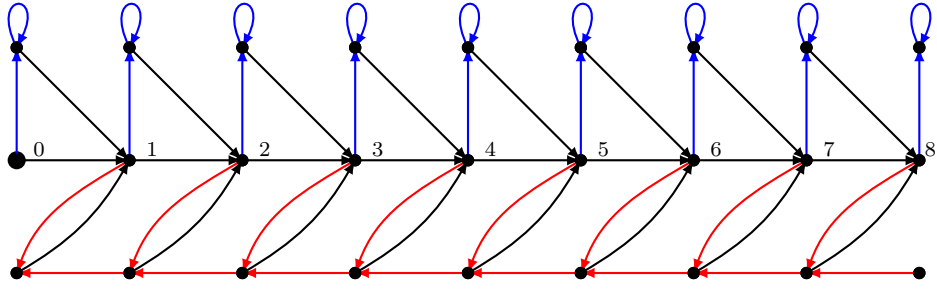
\begin{figure}[h]

	\begin{center}
		\begin{tikzpicture}[scale=1.5,main node/.style={circle,draw,font=\Large\bfseries}]

			\foreach \x in {0,1,2,3,4,5,6,7,8}
			{
				\draw (\x,0) circle (0.05cm);
				\fill (\x,0) circle (0.05cm);
				\draw (\x,-1) circle (0.05cm);
				\fill (\x,-1) circle (0.05cm);
				\draw (\x,1) circle (0.05cm);
				\fill (\x,1) circle (0.05cm);
			}

			\fill (0,0) circle (0.08cm);

			\foreach \x in {0,...,8}
			{
				\draw[thick,-latex,blue ] (\x,1)  ..  controls (\x-0.25,1+0.5) and  (\x+0.25,1+0.5) .. (\x,1+0) ;	
				\draw[thick,-latex,blue ] (\x,0) to (\x,1+0) ;	
			}
			
			\foreach \x in {0,...,7}
			{
				\draw[thick,latex-, ] (\x+1,0)  to (\x,1+0) ;	
			}
			
			\foreach \x in {0,...,7}
			{
				
				\draw[thick,-latex ] (\x,-1)[out=30,in=240] to (\x+1,0) ;	
				
			}

			\foreach \x in {0,1,2,3,4,5,6,7}
			{
				\draw[thick, -latex] (\x,0) to  (\x+1,0);	
				\draw[thick,red, latex-] (\x,-1) to  (\x+1,-1);	
				\draw[thick,red, -latex] (\x+1,0) [out=210,in=70]to  (\x,-1);	
				\node at  (\x+0.2,0.1){\tiny$\x$};
			}			
			
			\foreach \x in {0,1,2,3,4,5,6,7,8}
			{
			}

			\node at  (8+0.1,0.1){\tiny$8$};
			
			\foreach \x in {0,1,2,3,4,5,6,7,8}
			{
				\draw (\x,0) circle (0.05cm);
				\fill (\x,0) circle (0.05cm);
				\draw (\x,-1) circle (0.05cm);
				\fill (\x,-1) circle (0.05cm);
				\draw (\x,1) circle (0.05cm);
				\fill (\x,1) circle (0.05cm);
			}
		\end{tikzpicture}
	\end{center}
	\caption{Graph (automaton) to recognize  Motzkin paths with forbidden subwords DH and HD. }
	\label{purpelDHHD}
\end{figure}
The generating functions relative to Figure~\ref{purpelDHHD} satisfies the system
\begin{align*}
	f_i&=zf_i+zg_i,\ f_0=\frac{z}{1-z},\\
	g_{i+1}&=zf_i+zg_i+zh_i,\ g_0=1,\\
	h_i&=zg_{i+1}+zh_{i+1},
\end{align*}
which translates into bivariate generating functions as follows:
\begin{align*}
	F(u)&=zF(u)+zG(u),\\*
	G(u)&=1+zuF(u)+zuG(u)+zuH(u),\\*
	H(u)&=\frac zu\bigl(G(u)-1\bigr)+\frac zu\bigl(H(u)-H(0)\bigr).
\end{align*}
Solving the system leads to 
\begin{align*}
F(u)&=\frac{z(z^2u+z^2uH(0)-u+z)}{-z^3u-u+zu+zu^2+z-z^2},\\
G(u)&=\frac{(1-z)(z^2u+z^2uH(0)-u+z)}{-z^3u-u+zu+zu^2+z-z^2},\\
H(u)&=\frac{-z(-H(0)+zH(0)+zu+zuH(0))}{-z^3u-u+zu+zu^2+z-z^2}.
\end{align*}
We factorize the denominator as $z(u-r_1)(u-r_2)$ with
\begin{gather*}
	r_1=\frac{1-z+z^3-W}{2z},\ r_2=\frac{1-z+z^3+W}{2z},\\ W=\sqrt{(1+z-2z^2-z^3)(1-3z+2z^2-z^3)}.
\end{gather*}
Dividing out the factor $u-r_1$ according to the kernel method leads to the expressions
\begin{equation*}
	F(u)=\frac{z}{1-z-ur_1},\ G(u)=\frac{1-z}{1-z-ur_1},\ H(u)=\frac{r_1-z}{z(1-z-ur_1)},\ S(u)=\frac{r_1}{z(1-z-ur_1)},
\end{equation*}
whence
\begin{equation*}
	S(0)=\frac{r_1}{z(1-z)},\quad [u^j]S(u)=s_j=\frac1z\Big(\frac{r_1}{1-z}\Big)^{j+1};
\end{equation*}
the latter is the generating function of restricted Motzkin paths ending on level $j$.
The generating function $S(0)$ of restricted Motzkin excursions is in the OEIS (A329701), but
\begin{equation*}
	S(1)=1+2z+5z^2+11z^3+26z^4+60z^5+142z^6+334z^7+794z^8+1885z^9+\cdots,
\end{equation*}
the generating function of meanders,  is not.

Using the denominator $-z^3u-u+zu+zu^2+z-z^2$, we derive recursions
\begin{align*}
	zs_{j-2}+(z-z^3-1)s_{j-1}+z(1-z)s_j&=0, \ j\ge 2,\\
	(z-z^3-1)s_{0}+z(1-z)s_1&=-1,
\end{align*}
which, upon restricting ourselves to states not larger than  $K$, can be written as a linear system
{\small
\begin{equation*}
	\begin{pmatrix}
-1+z-z^3		 & z(1-z) & 0&\cdots  &0  \\
		z & -1+z-z^3 &z(1-z)& \cdots & 0 \\
		&		z & -1+z-z^3 &z(1-z) & 0 \\
		\vdots  & \ddots  & \ddots & \ddots  & \vdots\\
		0 & 0 & \cdots&z &-1+z-z^3
	\end{pmatrix}
	\begin{pmatrix}
		s_0\\
		s_1\\
		s_2\\
		\vdots\\
		s_K
	\end{pmatrix}=
	\begin{pmatrix}
		-1\\
		0\\
		0\\
		\vdots\\
		0
	\end{pmatrix}
\end{equation*}}
Denoting the determinant of the matrix without first row and column (a $K\times K$-matrix) by $\mD_K$, we derive a recursion
\begin{equation*}
	\mathcal{D}_i=(z-z^3-1  )\mathcal{D}_{i-1}-z^2(1-z)\mathcal{D}_{i-2},\ \mathcal{D}_0=1,\ \mathcal{D}_1=z-z^3-1 
\end{equation*}
and, using $t_i=\mD_{i-1}/\mD_{i}$ and $\tau_i=-z^2(1-z)t_i$, ($\tau_0=0$),
\begin{equation*}
	t_i=\frac{1}{z-z^3-1-z^2(1-z)t_{i-1}},\quad \tau_i=\frac{z^2(1-z)}{1-z+z^3-\tau_{i-1}}.
\end{equation*}
There is a Binet formula for the determinant
\begin{equation*}
	\mathcal{D}_j=\frac{(-z)^{j+1}}{W}\big(r_1^{j+1}-r_2^{j+1}\big),
\end{equation*}
and the solution of the component $s_0$ may be achieved by Cramer's rule as a quotient of determinants:
\begin{align*}
	s_0=\sigma_{K+1}&=\frac{-\mathcal{D}_{K}}{(z-z^3-1)\mathcal{D}_{K}-z^2(1-z)\mathcal{D}_{K-1}	}\\
	&=\frac{-1}{(z-z^3-1)-z^2(1-z)t_{K}	}=\frac{1}{1-z+z^3-\tau_{K}	}.
\end{align*}
So we have the continued fraction type expansion
\begin{equation*}
	\sigma_{K+1}=\frac{1}{1-z+z^3-\tau_{K}	},\quad   \tau_K=\frac{z^2(1-z)}{1-z+z^3-\tau_{K-1}},\ \tau_0=0.
\end{equation*}

Now we compute the individual quantities, according to restricted Motzkin paths ending in the first, resp.\ second, resp.\ third layer.
$f_0=\frac{z}{1-z}$ and $g_0=1$, regardless of $K$, so we concentrate on the third layer.

{\small
\begin{equation*}
	\begin{pmatrix}
		-1+z+z^2-z^3 & z(1-z) & 0&\cdots  &0  \\
		z & -1+z-z^3 &z(1-z)& \cdots & 0 \\
		&		z & -1+z-z^3 &z(1-z) & 0 \\
		\vdots  & \ddots  & \ddots & \ddots  & \vdots\\
		0 & 0 & \cdots&z &-1+z-z^3
	\end{pmatrix}
	\begin{pmatrix}
		h_0\\
		h_1\\
		h_2\\
		\vdots\\
		h_K
	\end{pmatrix}=
	\begin{pmatrix}
		-z^2\\
		0\\
		0\\
		\vdots\\
		0
	\end{pmatrix}
\end{equation*}
}
Again by Cramer's rule, 
\begin{equation*}
h_0=\eta_{K+1}=\frac{-z^2\mD_{K}}{(-1+z+z^2-z^3)\mD_{K}-z^2(1-z)\mD_{K-1}},\quad \eta_0=0.
\end{equation*}
This formula produces $\eta_1=\frac{z^2}{(1-z)(1-z^2)}$, which can be seen directly from the graph as a check. Further, we can see that
$\sigma_1=\frac{1}{(1-z)(1-z^2)}$, which also checks.
Finally we derive continued fraction type formulas for the individual quantities $\phi_K$, $\gamma_K$.
{\small
	\begin{equation*}
		\begin{pmatrix}
			-z		 & (1-z)^2 & 0&\cdots  &0  \\
			z & -1+z-z^3 &z(1-z)& \cdots & 0 \\
			&		z & -1+z-z^3 &z(1-z) & 0 \\
			\vdots  & \ddots  & \ddots & \ddots  & \vdots\\
			0 & 0 & \cdots&z &-1+z-z^3
		\end{pmatrix}
		\begin{pmatrix}
			f_0\\
			f_1\\
			f_2\\
			\vdots\\
			f_K
		\end{pmatrix}=
		\begin{pmatrix}
			-z\\
			0\\
			0\\
			\vdots\\
			0
		\end{pmatrix}
\end{equation*}}
{\small
	\begin{equation*}
		\begin{pmatrix}
			(z+1)(z-1)^2		 & z(z-1) & 0&\cdots  &0  \\
			z & -1+z-z^3 &z(1-z)& \cdots & 0 \\
			&		z & -1+z-z^3 &z(1-z) & 0 \\
			\vdots  & \ddots  & \ddots & \ddots  & \vdots\\
			0 & 0 & \cdots&z &-1+z-z^3
		\end{pmatrix}
		\begin{pmatrix}
			h_0\\
			h_1\\
			h_2\\
			\vdots\\
			h_K
		\end{pmatrix}=
		\begin{pmatrix}
			z^2\\
			0\\
			0\\
			\vdots\\
			0
		\end{pmatrix}
\end{equation*}}

\begin{align*}
f_0=\phi_{K+1}&=\frac{-z\mD_{K}}{-z\mD_{K}-z(1-z)^2\mD_{K-1}}
=\frac{1}{1+(1-z)^2t_{K}}=\frac{z^2}{z^2 -(1-z)\tau_K  },\ \phi_0=\frac z{1-z},\\
g_0=\gamma_{K}&=1,\\
h_0=\eta_{K+1}&=\frac{z^2\mD_{K}}{(z+1)(z-1)^2\mD_{K}-z^2(z-1)\mD_{K-1}}=
\frac{z^2}{1-z-z^2+z^3-\tau_{K}},\ \eta_0=0.
\end{align*}

%% file: UDDU.tex
\section{UD and DU are forbidden, A004149}
\begin{figure}[h]

	\begin{center}
		\begin{tikzpicture}[scale=1.5,main node/.style={circle,draw,font=\Large\bfseries}]

			\fill (0,-1) circle (0.08cm);

			\foreach \x in {0,1,2,3,4,5,6,7}
			{
				\draw[thick, -latex] (\x,0-1) to  (\x+1,0);   
				\draw[ thick,red, -latex] (\x+1,-1) to  (\x,-2);    
				\draw[thick,latex-,blue ] (\x,-1)  ..  controls (\x-0.3,-1+0.3) and  (\x-0.3,-1-0.3) .. (\x,-1) ;	
				\draw[ thick,red, - latex] (\x+1,-2)   to  (\x,-2);    
				\draw[ thick, latex-] (\x+1,0) to  (\x,0);    
				\node at  (\x+0.2,0.1){\tiny$\x$};
			}  
	
			\foreach \x in {8}
		{
		\draw[thick,latex-,blue ] (\x,-1)  ..  controls (\x-0.3,-1+0.3) and  (\x-0.3,-1-0.3) .. (\x,-1) ;
	}
		                      
			\foreach \x in {0,1,2,3,4,5,6,7,8}
			{
				
				\draw[ thick,blue, -latex] (\x,0)to  (\x,-1);    
				\draw[ thick,blue, -latex] (\x,-2)to  (\x,-1);    
				
			}                        
			
			\node at  (8+0.2,0.1){\tiny$8$};
			
			\foreach \x in {0,1,2,3,4,5,6,7,8}
			{
				\draw (\x,0) circle (0.05cm);
				\fill (\x,0) circle (0.05cm);
				\draw (\x,-1) circle (0.05cm);
				\fill (\x,-1) circle (0.05cm);
				\draw (\x,-2) circle (0.05cm);
				\fill (\x,-2) circle (0.05cm);
			}
		\end{tikzpicture}
	\end{center}
	\caption{Graph (automaton) to recognize  Motzkin paths with forbidden subwords UD and DU.}
	\label{purpelUDDU}
\end{figure}
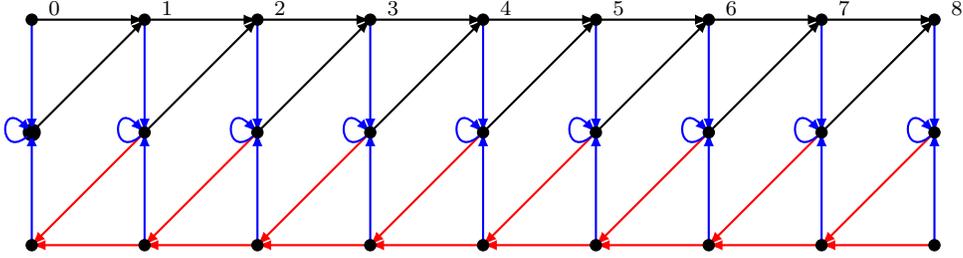
According to Figure \ref{purpelUDDU}, we have
\begin{align*}
	f_{i+1}&=zf_i+zg_i,\ f_0=0,\\*
	g_i&=[i=0]+zf_i+zg_i+zh_i,\\*
	h_i&=zg_{i+1}+zh_{i+1}
\end{align*}
and as a consequence
\begin{align*}
	F(u)&=zuF(u)+zuG(u),\\*
	G(u)&=1+zF(u)+zG(u)+H(u),\\*
	H(u)&=\frac zu(G(u)-G(0))+\frac zu(H(u)-H(0)).
\end{align*}
We have $zu^2+(-1+z-z^2-z^3)u+z=z(u-r_1)(u-r_2)$ with
\begin{equation*}
r_1=\frac{1-z+z^2+z^3-W}{2z},\ r_2=\frac{1-z+z^2+z^3+W}{2z},\ W=\sqrt{(1-z^4)(1-2z-z^2)}.
\end{equation*}
Dividing out the bad factor $u-r_1$ from the components of the solved system, we end up with
\begin{equation*}
	F(u)=-1+\frac{1}{1-ur_1},\ 
G(u)=\frac{r_1(1-zu)}{z(1-ur_1)},\
H(u)=\frac{r_1(1-z)-z}{z^2(1-ur_1)},\ 
S(u)=\frac{r_1-z}{z^2(1-ur_1)}.
\end{equation*}
From the denominator $zu^2+(-1+z-z^2-z^3)u+z$ we derive the recursion
\begin{align*}
	zs_{j-2}+(-1+z-z^2-z^3)s_{j-1}+zs_j&=0,\ j\ge2\\
	(-1+z-z^3)s_{0}+zs_1&=-1+z^2
\end{align*}
or, ignoring $s_i$ with $i>K$, in matrix form
{\footnotesize
\begin{equation*}
	\begin{pmatrix}
		-1+z-z^3& z& 0&\cdots  &0  \\
		z & -1+z-z^2-z^3 &z& \cdots & 0 \\
		&		z & -1+z-z^2-z^3 &z & 0 \\
		\vdots  & \ddots  & \ddots & \ddots  & \vdots\\
		0 & 0 & \cdots&z &-1+z-z^2-z^3
	\end{pmatrix}
	\begin{pmatrix}
		s_0\\
		s_1\\
		s_2\\
		\vdots\\
		s_K
	\end{pmatrix}=
	\begin{pmatrix}
		-1+z^2\\
		0\\
		0\\
		\vdots\\
		0
	\end{pmatrix}
\end{equation*}
}
As always, we consider the $K\times K$ determinant after erasing the first row resp.\ column and derive a recursion resp. solution for it:
\begin{equation*}
	\mD_{K}=(-1+z-z^2-z^3)\mD_{K-1}-z^2\mD_{K-2},\quad 	\mD_{K}=\frac{(-z)^{K+1}}{W}(r_1^{K+1}-r_2^{K+1}).
\end{equation*}
Further, introducing $t_i=\mD_{i-1}/\mD_{i}$ and $\tau_i=-z^2t_i$ we derive the continued fraction type recursion
\begin{equation*}
	t_{K}=\frac{1}{-1+z-z^2-z^3-z^2t_{K-1}}, \quad 	\tau_{K}=\frac{z^2}{1-z+z^2+z^3-\tau_{K-1}},\ \tau_0=0.
\end{equation*}
The $s_0$-component of the system can be solved by Cramer's rule:
\begin{align*}
	s_0=\sigma_{K+1}&=\frac{(-1+z^2)\mD_{K}}{(-1+z-z^3)\mD_{K}-z^2\mD_{K-1}}
	=\frac{(-1+z^2)}{(-1+z-z^3)-z^2t_K}
	=\frac{1-z^2}{1-z+z^3-\tau_{K}}.
\end{align*}

Working with $s_j^{[K]}$ (states larger than $K$ ignored) \emph{alone} might lead to small anormalities. According to the Figures \ref{purpelklein}
(which shows the instance $K=3$),
$f_j^{[K]}+g_j^{[K]}+h_j^{[K]}$ is considered. But it should really be
$f_j^{[K]}+g_j^{[K]}+h_j^{[K-1]}$ as the $K$-state in the third layer cannot be reached.

To achieve this, we have to compute $f_j^{[K]}$, $g_j^{[K]}$, $h_j^{[K]}$ separately.

\begin{figure}[h]

	\begin{center}
		\begin{tikzpicture}[scale=0.8,main node/.style={circle,draw,font=\Large\bfseries}]
			\definecolor{ashgrey}{rgb}{0.7, 0.75, 0.71}	
			\definecolor{arschcol}{rgb}{0.74, 0.81, 0.69}

			\path[fill=arschcol,line width=0pt,rounded corners=0.2cm] 
			(-0.5,0.4) to (3.4,0.4) to (3.4,-2.2) to
			(-0.5,-2.2)  -- cycle;

			\fill (0,-1) circle (0.08cm);

			\foreach \x in {0,1,2,3,4,5,6,7}
			{
				\draw[thick, -latex] (\x,0-1) to  (\x+1,0);   
				\draw[ thick,red, -latex] (\x+1,-1) to  (\x,-2);    
				\draw[thick,latex-,blue ] (\x,-1)  ..  controls (\x-0.3,-1+0.3) and  (\x-0.3,-1-0.3) .. (\x,-1) ;	
				\draw[ thick,red, - latex] (\x+1,-2)   to  (\x,-2);    
				\draw[ thick, latex-] (\x+1,0) to  (\x,0);    
				\node at  (\x+0.2,0.2){\tiny$\x$};
			}                        
			\foreach \x in {0,1,2,3,4,5,6,7,8}
			{
				
				\draw[ thick,blue, -latex] (\x,0)to  (\x,-1);    
				\draw[ thick,blue, -latex] (\x,-2)to  (\x,-1);    
				
			}                        
			
			\node at  (8+0.2,0.1){\tiny$8$};
			
			\foreach \x in {0,1,2,3,4,5,6,7,8}
			{
				\draw (\x,0) circle (0.05cm);
				\fill (\x,0) circle (0.05cm);
				\draw (\x,-1) circle (0.05cm);
				\fill (\x,-1) circle (0.05cm);
				\draw (\x,-2) circle (0.05cm);
				\fill (\x,-2) circle (0.05cm);
			}
			
			
		\end{tikzpicture}%
		\hspace{0.5cm}%
		\begin{tikzpicture}[scale=0.8,main node/.style={circle,draw,font=\Large\bfseries}]
			\definecolor{ashgrey}{rgb}{0.7, 0.75, 0.71}	
			\definecolor{arschcol}{rgb}{0.74, 0.81, 0.69}

			\path[fill=arschcol,line width=0pt,rounded corners=0.2cm] 
			(-0.5,0.4) to (3.4,0.4) to (3.4,-1.2) to
			(2.4,-2.2) to
			(-0.5,-2.2)  -- cycle;

			\fill (0,-1) circle (0.08cm);

			\foreach \x in {0,1,2,3,4,5,6,7}
			{
				\draw[thick, -latex] (\x,0-1) to  (\x+1,0);   
				\draw[ thick,red, -latex] (\x+1,-1) to  (\x,-2);    
				\draw[thick,latex-,blue ] (\x,-1)  ..  controls (\x-0.3,-1+0.3) and  (\x-0.3,-1-0.3) .. (\x,-1) ;	
				\draw[ thick,red, - latex] (\x+1,-2)   to  (\x,-2);    
				\draw[ thick, latex-] (\x+1,0) to  (\x,0);    
				\node at  (\x+0.2,0.2){\tiny$\x$};
			}                        
			\foreach \x in {0,1,2,3,4,5,6,7,8}
			{
				
				\draw[ thick,blue, -latex] (\x,0)to  (\x,-1);    
				\draw[ thick,blue, -latex] (\x,-2)to  (\x,-1);    
				
			}                        
			
			\node at  (8+0.2,0.1){\tiny$8$};
			
			\foreach \x in {0,1,2,3,4,5,6,7,8}
			{
				\draw (\x,0) circle (0.05cm);
				\fill (\x,0) circle (0.05cm);
				\draw (\x,-1) circle (0.05cm);
				\fill (\x,-1) circle (0.05cm);
				\draw (\x,-2) circle (0.05cm);
				\fill (\x,-2) circle (0.05cm);
			}
			
			
		\end{tikzpicture}
	\end{center}
	\label{purpelklein}
\end{figure}
Now $f_0=\phi_{K}=0$, regardless of $K$. Further

{\footnotesize
	\begin{equation*}
	\begin{pmatrix}
		-1+z-z^3& z& 0&\cdots  &0  \\
		z & -1+z-z^2-z^3 &z& \cdots & 0 \\
		&		z & -1+z-z^2-z^3 &z & 0 \\
		\vdots  & \ddots  & \ddots & \ddots  & \vdots\\
		0 & 0 & \cdots&z &-1+z-z^2-z^3
	\end{pmatrix}
	\begin{pmatrix}
		g_0\\
		g_1\\
		g_2\\
		\vdots\\
		g_K
	\end{pmatrix}=
	\begin{pmatrix}
		-1\\
		z\\
		0\\
		\vdots\\
		0
	\end{pmatrix}
\end{equation*}
}
and thus
\begin{equation*}
	g_0=\gamma_{K+1}=\frac{-\mD_{K}-z^2\mD_{K-1}}{(-1+z-z^3)\mD_{K}-z^2\mD_{K-1}},\ \gamma_0=\frac1{1-z}.
\end{equation*}
Similarly,
{\footnotesize
\begin{equation*}
	\begin{pmatrix}
		-1+2z-z^2+z^4		& z(1-z)& 0&\cdots  &0  \\
		z & -1+z-z^2-z^3 &z& \cdots & 0 \\
		&		z & -1+z-z^2-z^3 &z & 0 \\
		\vdots  & \ddots  & \ddots & \ddots  & \vdots\\
		0 & 0 & \cdots&z &-1+z-z^2-z^3
	\end{pmatrix}
	\begin{pmatrix}
		h_0\\
		h_1\\
		h_2\\
		\vdots\\
		h_K
	\end{pmatrix}=
	\begin{pmatrix}
		-z^3\\
		0\\
		0\\
		\vdots\\
		0
	\end{pmatrix}
\end{equation*}
}
and thus
\begin{equation*}
	h_0=\eta_{K+1}=\frac{-z^3\mD_{K}}{(-1+2z-z^2+z^4)\mD_{K}-z^2(1-z)\mD_{K-1}},\ \eta_0=0.
\end{equation*}
The generating function of Motzkin excursions
\begin{equation*}
S(0)=1+z+z^2+2z^3+4z^4+8z^5+16z^6+33z^7+69z^8+146z^9+312z^{10}+\cdots
\end{equation*}
is A004149 in OEIS, and of Motzkin meanders
\begin{equation*}
	S(1)=1+2z+4z^2+9z^3+20z^4+45z^5+102z^6+233z^7+535z^8+1234z^9+2857z^{10}\cdots
\end{equation*}
is A308435 in OEIS.

%% file: DDDU.tex
\section{DD and DU are forbidden, A023431}
\begin{figure}[h]

	\begin{center}
		\begin{tikzpicture}[scale=1.5,main node/.style={circle,draw,font=\Large\bfseries}]

			\fill (0,-1) circle (0.08cm);

			\foreach \x in {0,1,2,3,4,5,6,7}
			{
				\draw[thick, -latex] (\x,0-1) [in=-145, out=55]to  (\x+1,0);   
				\draw[ thick,red, -latex] (\x+1,-1) to  (\x,-2);    
				\draw[thick,latex-,blue ] (\x,-1)  ..  controls (\x-0.3,-1+0.3) and  (\x-0.3,-1-0.3) .. (\x,-1) ;	
				\draw[ thick,red, - latex] (\x+1,0)   to  (\x,-2);    
				\draw[ thick, latex-] (\x+1,0) to  (\x,0);    
				\node at  (\x+0.2,0.1){\tiny$\x$};
			}        
		
			\foreach \x in {8}
		{
			\draw[thick,latex-,blue ] (\x,-1)  ..  controls (\x-0.3,-1+0.3) and  (\x-0.3,-1-0.3) .. (\x,-1) ;	
		}

			\foreach \x in {0,1,2,3,4,5,6,7,8}
			{
				
				\draw[ thick,blue, -latex] (\x,0)to  (\x,-1);    
				\draw[ thick,blue, -latex] (\x,-2)to  (\x,-1);    
				
			}                        
			
			\node at  (8+0.2,0.1){\tiny$8$};
			
			\foreach \x in {0,1,2,3,4,5,6,7,8}
			{
				\draw (\x,0) circle (0.05cm);
				\fill (\x,0) circle (0.05cm);
				\draw (\x,-1) circle (0.05cm);
				\fill (\x,-1) circle (0.05cm);
				\draw (\x,-2) circle (0.05cm);
				\fill (\x,-2) circle (0.05cm);
			}
		\end{tikzpicture}
	\end{center}
	\caption{Graph (automaton) to recognize  Motzkin paths with forbidden subwords DD and DU.}
	\label{purpelDDDU}
\end{figure}
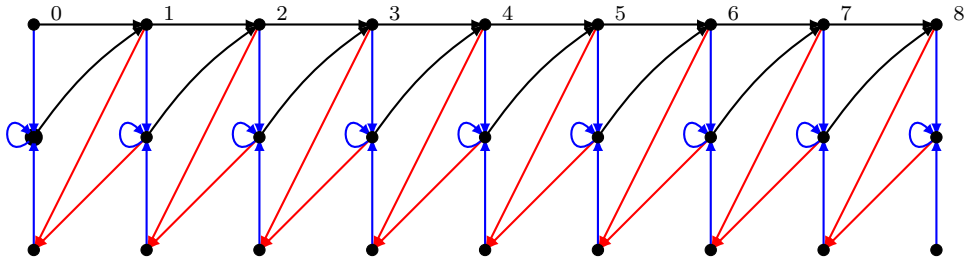
According to Figure \ref{purpelDDDU}, we find
\begin{align*}
	f_{i+1}&=zf_i+zg_i, \ f_0=0,\\
	g_i&=[i=0]+zf_i+zg_i+zh_i,\ i\ge0,\\
	h_{i}&=zf_{i+1}+zg_{i+1}
\end{align*}
and thus
\begin{align*}
	F(u)&=zuF(u)+zuG(u),\\*
	G(u)&=1+zF(u)+zG(u)+zH(u),\\*
	H(u)&=\frac zuF(u)+\frac zu(G(u)-G(0)).
\end{align*}
According to the denominator $z^2+(-1+z)u+zu^2=z(u-r_1)(u-r_2)$ we find the simplified versions
\begin{equation*}
	F(u)=-1+\frac{z}{z-ur_1},\ G(u)=1+\frac{r_1-z^2}{z(z-ur_1)},\ H(u)=\frac{r_1^2}{z(z-ur_1)},\ S(u)=\frac{r_1(1+r_1)}{z(z-ur_1)}
\end{equation*}
with
\begin{equation*}
	r_1=\frac{1-z-W}{2z},\ r_2=\frac{1-z+W}{2z},\quad W:=\sqrt{1-2z+z^2-4z^3}.
\end{equation*}
The denominator provides the recursion
\begin{align*}
	zs_{j-2}+(z-1)s_{j-1}+z^2s_j&=0,\ j\ge2\\
	-(1-z-z^3)s_0+z^2s_1&=-1-z^2
\end{align*}
and in matrix form
\begin{equation*}
	\begin{pmatrix}
		-1+z+z^3& z^2& 0&\cdots  &0  \\
		z & z-1 &z^2& \cdots & 0 \\
		&		z & z-1 &z^2 & 0 \\
		\vdots  & \ddots  & \ddots & \ddots  & \vdots\\
		0 & 0 & \cdots&z &z-1
	\end{pmatrix}
	\begin{pmatrix}
		s_0\\
		s_1\\
		s_2\\
		\vdots\\
		s_K
	\end{pmatrix}=
	\begin{pmatrix}
		-1-z^2\\
		0\\
		0\\
		\vdots\\
		0
	\end{pmatrix}
\end{equation*}
With the usual determinant we find
\begin{equation*}
	\mathcal{D}_{K}=(z-1)\mathcal{D}_{K-1}-z^3\mathcal{D}_{K-2}=\frac{(-z)^{K+1}}{W}(r_1^{K+1}-r_2^{K+1})
\end{equation*}
\begin{equation*}
	t_{K}=\frac1{(z-1)-z^3t_{K-1}},\ z^3t_{K}=\frac{z^3}{-1+z-z^3t_{K-1}},\\ \tau_K=\frac{z^3}{1-z-\tau_{K-1}},\quad \tau_0=0
\end{equation*}
\begin{align*}
	s_0=\sigma_{K+1}&=\frac{-(1+z^2)\mathcal{D}_{K}}{(-1+z+z^3)\mathcal{D}_{K}-z^3\mathcal{D}_{K-1}}
	=\frac{1+z^2}{1-z-z^3+z^3t_K}\\
	&=\frac{1+z^2}{1-z-z^3-\tau_K},\quad \sigma_0=\frac1{1-z}. \quad \tau_K=\frac{-z^3\mathcal{D}_{K-1}}{\mathcal{D}_{K}}
\end{align*}
Here it is again wise to compute the three components (with bound $K$) separately, since state $K$ in the third layer cannot be reached.
$f_0=\phi_K$ is always zero, but
\begin{equation*}
	\begin{pmatrix}
		-1+z+z^3& z^2& 0&\cdots  &0  \\
		z & z-1 &z^2& \cdots & 0 \\
		&		z & z-1 &z^2 & 0 \\
		\vdots  & \ddots  & \ddots & \ddots  & \vdots\\
		0 & 0 & \cdots&z &z-1
	\end{pmatrix}
	\begin{pmatrix}
		g_0\\
		g_1\\
		g_2\\
		\vdots\\
		g_K
	\end{pmatrix}=
	\begin{pmatrix}
		-1\\
		z\\
		0\\
		\vdots\\
		0
	\end{pmatrix}
\end{equation*}
Consequently, 
\begin{equation*}
\gamma_{K+1}=\frac{-\mD_{K}-z^3\mD_{K-1}}{(-1+z+z^3)\mD_{K}-z^3 \mD_{K-1}},\quad \gamma_0=\frac1{1-z}.
\end{equation*}
Similarly,
\begin{equation*}
	\begin{pmatrix}
		-1+2z-z^2+z^3		& z^2(1-z)& 0&\cdots  &0  \\
		z & z-1 &z^2& \cdots & 0 \\
		&		z & z-1 &z^2 & 0 \\
		\vdots  & \ddots  & \ddots & \ddots  & \vdots\\
		0 & 0 & \cdots&z &z-1
	\end{pmatrix}
	\begin{pmatrix}
		h_0\\
		h_1\\
		h_2\\
		\vdots\\
		h_K
	\end{pmatrix}=
	\begin{pmatrix}
-z^2\\
		z\\
		0\\
		\vdots\\
		0
	\end{pmatrix}
\end{equation*}
Consequently, 
\begin{equation*}
	\eta_{K+1}=\frac{-z^2\mD_{K}}{(-1+2z-z^2+z^3)\mD_{K}-z^3(1-z) \mD_{K-1}},\quad \eta_0=0.
\end{equation*}
The generating function $S(0)$ of excursions is A023431 in OEIS, and
\begin{align*}
S(1)&=\frac{1+z}{1-2z-z^2}+\frac{(z-1)r_1}{z(1-2z-z^2)}\\&=
1+2z+5z^2+12z^3+28z^4+66z^5+157z^6+374z^7+892z^8+2132z^9+\cdots,
\end{align*}
the generating function of meanders, is not in OEIS.

Since
\begin{equation*}
S(u)=\frac{r_1(1+r_1)}{z^2(1-\tfrac{ur_1}{z})}\quad\Longrightarrow \quad [u^j]S(u)=s_j=\frac{1+r_1}{z}\Bigl(\frac{r_1}{z}\Big)^{j+1};
\end{equation*}
the last quantity is the generating function of (restricted) Motzkin paths ending on level~$j$.

%% file: UHHU.tex
\section{UH and HU are forbidden, A329701}
\begin{figure}[h]

	\begin{center}
		\begin{tikzpicture}[scale=1.5,main node/.style={circle,draw,font=\Large\bfseries}]

			\fill (0,0) circle (0.08cm);

			\foreach \x in {0,...,8}
			{
				\draw[thick,-latex,blue ] (\x,1)  ..  controls (\x-0.25,1+0.5) and  (\x+0.25,1+0.5) .. (\x,1+0) ;	
				\draw[thick,-latex,blue ] (\x,0) to (\x,1+0) ;	
			}
			
			\foreach \x in {0,...,7}
			{
				\draw[thick,latex-, ] (\x,0)  to (\x+1,1+0) ;	
			}
			
			\foreach \x in {0,...,7}
			{
				
				\draw[thick,-latex ] (\x+1,-1)[out=150,in=-70] to (\x,0) ;	
			}

			\foreach \x in {0,1,2,3,4,5,6,7}
			{
				\draw[thick, latex-] (\x,0) to  (\x+1,0);	
				
				\draw[thick,red, -latex] (\x,0) [out=-30,in=110]to  (\x+1,-1);	
				\node at  (\x+0.2,0.1){\tiny$\x$};
			}			
			\foreach \x in {1,2,3,4,5,6,7}
			{
				
				\draw[thick,red, -latex] (\x,-1) to  (\x+1,-1);	
			}	
			\foreach \x in {0,1,2,3,4,5,6,7,8}
			{
			}

			\node at  (8+0.1,0.1){\tiny$8$};
			
			\foreach \x in {0,1,2,3,4,5,6,7,8}
			{
				\draw (\x,0) circle (0.05cm);
				\fill (\x,0) circle (0.05cm);
				\draw (\x,1) circle (0.05cm);
				\fill (\x,1) circle (0.05cm);
			}
			\foreach \x in {1,2,3,4,5,6,7,8}
			{
				\draw (\x,-1) circle (0.05cm);
				\fill (\x,-1) circle (0.05cm);
				
			}
		\end{tikzpicture}
	\end{center}
	\caption{Graph (automaton) to recognize  Motzkin paths with forbidden subwords UH and HU.}
\label{purpelUHHU}
\end{figure}
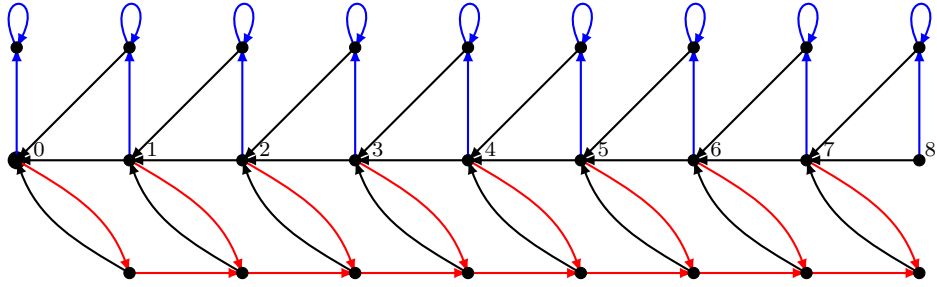
The following recursions can be read off the automaton (Figure \ref{purpelUHHU}), by considering the last step separately:
\begin{align*}
	f_i&=zf_{i}+zg_{i},\\
	g_{i}&=[i=0]+zf_{i+1}+zg_{i+1}+zh_{i+1},\\
	h_{i+1}&=zg_{i}+zh_{i}, \ i\ge0,\ h_0=0,
\end{align*}
and using double generating functions
\begin{align*}
	F(u)&=zF(u)+zG(u)=\frac{zG(u)}{1-z},\\
	G(u)&=1+\frac zu\big(F(u)-\frac z{1-z}G(0)\big)+\frac zu\big(G(u)-G(0)\big)+\frac zuH(u),\\
	H(u)&=zuG(u)+zuH(u).
\end{align*}
For $S=F+G+H$ we get
\begin{equation*}
	S=-{\frac { \left( -u+zu+zG(0) \right)  \left( {z}^{2}u-1 \right) }{			\left( u{z}^{3}-zu+{z}^{2}{u}^{2}-z+u-z{u}^{2} \right)  \left( -1+z
			\right) }}
\end{equation*}
We have the factorization $	( u{z}^{3}-zu+{z}^{2}{u}^{2}-z+u-z{u}^{2} )  =z(z-1)(u-r_1)(u-r_2)$
with 
\begin{equation*}
	r_1={\frac {1-z+{z}^{3}-W}{2z(1-z) }}, \ 
	r_2={\frac {1-z+{z}^{3}+W}{2z(1-z) }},\ W=\sqrt{(1+z-2z^2-z^3)(1-3z+2z^2-z^3)}.
\end{equation*}
After dividing out the factor $u-r_1$ from numerator and denominator and simplification,
\begin{gather*}
F(u)=\frac{(1-zu)r_1}{1-u(1-z)r_1},\ G(u)=\frac{(1-zu)(1-z)r_1}{z(1-u(1-z)r_1)},\ 
H(u)=\frac{(1-z)ur_1}{1-u(1-z)r_1},\\ S(u)=\frac{(1- z^2u)r_1}{z(1-u(1-z)r_1)}.
\end{gather*}
From this we find 
\begin{equation*}
S(0)=1+z+2 z^2+2 z^3+4 z^4+5 z^5+11 z^6+17 z^7+38 z^8+67 z^9+\cdots,
\end{equation*}
which is A329701. The generating function for meanders is
\begin{equation*}
	S(1)=1+2 z+3 z^2+5 z^3+9 z^4+17 z^5+34 z^6+69 z^7+144 z^8+302 z^9+\cdots,
\end{equation*}
which is not in OEIS.

For $j\ge1$,
\begin{equation*}
	[u^j]S(u)=[u^j]\frac{(1- z^2u)r_1}{z(1-u(1-z)r_1)}=
	\frac{(1-z)^jr_1^{j+1}}{z}-{z(1-z)^{j-1}r_1^j}
\end{equation*}
which is the generating function of restricted Motzkin paths ending at level $j$.

Now we move to paths of bounded height. We have 
\begin{equation*}
	S=\frac { \textsf{numerator} }{z(z-1){u}^{2}			 +u({z}^{3}-z+1)  -z}
\end{equation*}
and a similar shape for $F,G,H$. This leads to a recursion:
Set $s_j=[u^j]S$, then
\begin{align*}
	z(z-1)s_{j}+({z}^{3}-z+1)s_{j+1}-zs_{j+2}&=0,\quad j\ge1,\\
	z(z-1)s_{0}+({z}^{3}-z+1)s_{1}-zs_{2}&=-z^2,\\
	(z-1)s_0+zs_1&=-1,
\end{align*}
which is in matrix form
\begin{equation*}
	\begin{pmatrix}
		z-1 & z & 0&\cdots & &0  \\
		z(z-1) & z^3-z+1 &-z& \cdots && 0 \\
		&		z(z-1) & z^3-z+1 &-z& & 0 \\
		\vdots  & \ddots  & \ddots & \ddots  && \vdots\\
		0 & 0 & \cdots&&z(z-1) &z^3-z+1
	\end{pmatrix}
	\begin{pmatrix}
		s_0\\
		s_1\\
		s_2\\
		\vdots\\
		s_K
	\end{pmatrix}=
	\begin{pmatrix}
		-1\\
		-z^2\\
		0\\
		\vdots\\
		0
	\end{pmatrix}
\end{equation*}

Let $\mathcal{D}_K$ be the determinant with $K$ rows and columns, after deleting the first row and column. Expanding, we derive the recursion
($\mathcal{D}_0=1,\ \mathcal{D}_1=1-z+z^3$)
\begin{equation*}
	\mathcal{D}_K=(z^3-z+1)\mathcal{D}_{K-1}+z^2(z-1)\mathcal{D}_{K-2}
	=\frac{z^{K+1}(1-z)^{K+1}}{W}(r_2^{K+1}-r_1^{K+1}).
\end{equation*}
Rewriting this, we get, with $t_i=\frac{\mathcal{D}_{i-1}}{\mathcal{D}_{i}}$, $\tau_i =z^2(z-1)t_i$, $\tau_0=0$
\begin{equation*}
t_K=\frac{1}{1-z+z^3+z^2(z-1)t_{K-1})}, \quad \tau_K=\frac{z^2(z-1)}{1-z+z^3+\tau_{K-1})}.
\end{equation*}
Solving the system by Cramer's rule, we get  the solution for $s_0$ (return to the $x$-axis) as
\begin{equation*}
	s_0=\sigma_{K+1}=\frac{-\mathcal{D}_K+z^3\mathcal{D}_{K-1}}{(z-1)\mathcal{D}_K-z^2(z-1)\mathcal{D}_{K-1}}=
	\frac{-1+z^3t_{K-1}}{(z-1)-z^2(z-1)t_{K-1}}.
\end{equation*}
A simple computation (using Maple's convert-parfrac construction) shows
\begin{align*}
	\sigma_{K+1}&=\cfrac{z}{1-z}+\cfrac{1-z}{1-z-\tau_K},
\end{align*}
which is of the continued fraction type.

The recursion for the functions $f_0$ is in matrix form
{\small
\begin{equation*}
	\begin{pmatrix}
		1-z-z^2+z^3 & -z & 0&\cdots & &0  \\
		z(z-1) & z^3-z+1 &-z& \cdots && 0 \\
		&		z(z-1) & z^3-z+1 &-z& & 0 \\
		\vdots  & \ddots  & \ddots & \ddots  && \vdots\\
		0 & 0 & \cdots&&z(z-1) &z^3-z+1
	\end{pmatrix}
	\begin{pmatrix}
		f_0\\
		f_1\\
		f_2\\
		\vdots\\
		f_K
	\end{pmatrix}=
	\begin{pmatrix}
		z\\
		-z^2\\
		0\\
		\vdots\\
		0
	\end{pmatrix}
\end{equation*}
}
whence
\begin{align*}
f_0=\phi_{K+1}&=\frac{z\mD_{K}-z^3\mD_{K-1}}{(1-z-z^2+z^3)\mD_{K}+z^2(z-1)\mD_{K-1}}\\
&= \frac{z}{1-z}+\frac{z^3}{(1-z)(1-z^2-z^2t_K)},\quad\phi_0=\frac{z}{1-z}.
\end{align*}
Likewise
{\footnotesize
	\begin{equation*}
		\begin{pmatrix}
			1-z-z^2+z^3 & -z & 0&\cdots & &0  \\
			z(z-1) & z^3-z+1 &-z& \cdots && 0 \\
			&		z(z-1) & z^3-z+1 &-z& & 0 \\
			\vdots  & \ddots  & \ddots & \ddots  && \vdots\\
			0 & 0 & \cdots&&z(z-1) &z^3-z+1
		\end{pmatrix}
		\begin{pmatrix}
			g_0\\
			g_1\\
			g_2\\
			\vdots\\
			g_K
		\end{pmatrix}=
		\begin{pmatrix}
			1-z\\
			z(z-1)\\
			0\\
			\vdots\\
			0
		\end{pmatrix}
	\end{equation*}
}
and
\begin{align*}
g_0=\gamma_{K+1}&=\frac{(1-z)\mD_{K}+z^2(z-1)\mD_{K-1}}{(1-z-z^2+z^3)\mD_{K}+z^2(z-1)\mD_{K-1}}\\
&= 1+\frac{z^2}{1-z^2-z^2t_K},\quad\gamma_0=1.
\end{align*}
This leads to $\phi_K=\gamma_K\frac{z}{1-z}$, which is also combinatorially clear.  $h_0=\eta_K=0$ for any $K$, so no more computations are
required.

%% file: UUHH.tex
\section{UU and HH are forbidden, A329666}
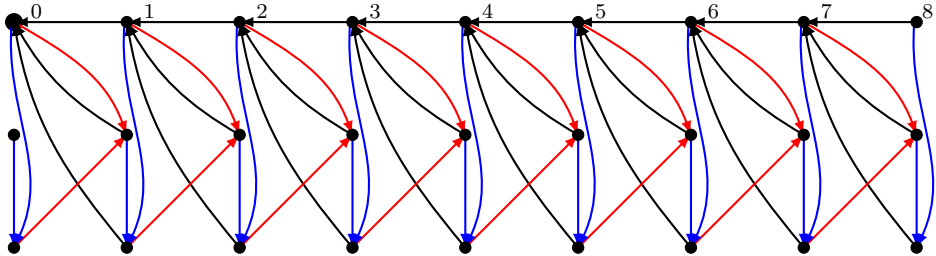
\begin{figure}[h]

	\begin{center}
		\begin{tikzpicture}[scale=1.5,main node/.style={circle,draw,font=\Large\bfseries}]

			\fill (0,0) circle (0.08cm);

			\foreach \x in {0,...,7}
			{
			}
			
			\foreach \x in {0,...,7}
			{
				
				\draw[thick,-latex ] (\x+1,-1)[out=150,in=-60] to (\x,0) ;	
			}

			\foreach \x in {0,1,2,3,4,5,6,7}
			{
				\draw[thick, latex-] (\x,0) to  (\x+1,0);	
				
				\draw[ thick,red, -latex] (\x,0) [out=-30,in=110]to  (\x+1,-1);	
				\node at  (\x+0.2,0.1){\tiny$\x$};
			}			
			
			\foreach \x in {0,1,2,3,4,5,6,7,8}
			{
				\draw[thick, -latex, blue] (\x,-1) to  (\x,-2);	
				
			}			
			
			\foreach \x in {1,2,3,4,5,6,7,8}
			{
				
				\draw[thick, -latex] [out=130,in=-80](\x,-2) to  (\x-1,0);	
				\draw[thick, -latex,red]  (\x-1,-2) to  (\x,-1);	
			}			
			
			\foreach \x in {0,1,2,3,4,5,6,7,8}
			{
				
				\draw[thick, -latex,blue] [out=-95,in=65](\x-0.02,0) to  (\x,-2);	
				
			}			
			
			\node at  (8+0.1,0.1){\tiny$8$};
			
			\foreach \x in {0,1,2,3,4,5,6,7,8}
			{
				\draw (\x,0) circle (0.05cm);
				\fill (\x,0) circle (0.05cm);
				\draw (\x,-1) circle (0.05cm);
				\fill (\x,-1) circle (0.05cm);
				\draw (\x,-2) circle (0.05cm);
				\fill (\x,-2) circle (0.05cm);
			}
		\end{tikzpicture}
	\end{center}
	\caption{Graph (automaton) to recognize  Motzkin paths with forbidden subwords UU and HH.}
\label{purpelUUHH}
\end{figure}
From Figure \ref{purpelUUHH} we derive
\begin{align*}
	f_j&=[j=0]+zf_{j+1}+zg_{j+1}+zh_{j+1},\\
	g_{j+1}&=zf_j+zh_j,\ g_0=0,\\
	h_j&=zf_j+zg_j
\end{align*}
and
\begin{align*}
	F(u)&=1+\frac zu[F(u)-F(0)]+\frac zuG(u)+\frac zu[H(u)-H(0)],\\
	G(u)&=zuF(u)+zuH(u),\\
	H(u)&=zF(u)+zG(u).
\end{align*}
Then
\begin{equation*}
	S(u)=F(u)+G(u)+H(u)=\frac{(1+z)(1+zu)(-u+z(1+z)F(0))}{z^2u^2+(z^3+z^2-1)u+z(1+z)}
\end{equation*}
Set $W=\sqrt { (1-2z-z^2+z^3) (1+2z+3z^2+ {z}^{3})}$ and
\begin{align*}
	r_1=\frac{1-z^2-z^3-W}{2z^2},\quad r_1=\frac{1-z^2-z^3+W}{2z^2},
\end{align*}
then $z^2u^2+(z^3+z^2-1)u+z(1+z)=z^2(u-r_1)(u-r_2)$.
Dividing out $u-r_1$, we compute $F(0)=\dfrac{1+zr_1}{1-z^2r_1}=\dfrac{r_1}{z(1+z)}$ and thus
$$S(0)=F(0)+H(0)=\frac{(1+z)(1+zr_1)}{1-z^2r_1}=\frac{r_1}{z}.$$
Using this, we get
\begin{gather*}
	F=\frac{(1+z)(1-z^2u)r_1}{z(1+z-uzr_1)},\quad G=\frac{(1+z)ur_1}{1+z-uzr_1},\quad H=\frac{(1+zu)r_1}{1+z-uzr_1}
\\	S(u)=\frac{(1+z)^2(1+zu)r_1}{z(1+z-uzr_1)}=\frac{(1+z)(1+zu)r_1}{z(1-u\frac{zr_1}{1+z})}.
\end{gather*}
\begin{equation*}
	[u^j]S(u)= \Bigl(\frac{z}{1+z}\Bigr)^{j-1}r_1^{j+1}+z\Bigl(\frac{z}{1+z}\Bigr)^{j-2}r_1^{j} ,\ j\ge1.
\end{equation*}
The denominator leads to the recursion
\begin{align*}
	z^2s_{j-2}+(z^3+z^2-1)s_{j-1}+z(1+z)s_j&=0,\quad j\ge3,\\
	z^2s_{0}+(z^3+z^2-1)s_{1}+z(1+z)s_2&=-z(1+z),\\
	-s_0+z(1+z)s_1=-1-z,
\end{align*}
which is in matrix form
\begin{equation*}
	\begin{pmatrix}
		-1 & z(1+z) & 0&\cdots & &0  \\
		z^2 & z^3+z^2-1 &z(1+z)& \cdots && 0 \\
		&		z^2 & z^3+z^2-1 &z(1+z)& & 0 \\
		\vdots  & \ddots  & \ddots & \ddots  && \vdots\\
		0 & 0 & \cdots&&z^2 &z^3+z^2-1
	\end{pmatrix}
	\begin{pmatrix}
		s_0\\
		s_1\\
		s_2\\
		\vdots\\
		s_K
	\end{pmatrix}=
	\begin{pmatrix}
		-1-z\\
		-z(1+z)\\
		0\\
		\vdots\\
		0
	\end{pmatrix}
\end{equation*}
Let $\mathcal{D}_K$ be the determinant with $K$ rows and columns, after deleting the first row and column. Expanding, we derive the recursion
\begin{equation*}
	\mathcal{D}_K=(z^3+z^2-1)\mathcal{D}_{K-1}-z^3(1+z)\mathcal{D}_{K-2},\ \mathcal{D}_0=1, \mathcal{D}_1=z^3+z^2-1.
\end{equation*}
It follows that
\begin{equation*}
	\mathcal{D}_K=\frac{(-z^2)^{K+1}}W(r_1^{i+1}-r_2^{i+1}).
\end{equation*}
Further
\begin{align*}
	s_0&=\frac{-(1+z)\mathcal{D}_K+z^2(1+z)^2\mathcal{D}_{K-1}}{-\mathcal{D}_K-z^3(1+z)\mathcal{D}_{K-1}}
	=(1+z)\frac{\mathcal{D}_K-z^2(1+z)\mathcal{D}_{K-1}}{\mathcal{D}_K+z^3(1+z)\mathcal{D}_{H-1}}
\end{align*}
with the abbreviations $t_i=\frac{\mathcal{D}_{i-1}}{\mathcal{D}_{i}}$ and $\tau_i=-(1+z)t_i$
\begin{align*}
	s_0=\sigma_{K+1}&=(1+z)\frac{1-z^2(1+z)t_{K}}{1+z^3(1+z)t_{K}}=-\frac{1+z}{z}+\frac{(1+z)^2}{z}\frac{1}{1-z^3\tau_{K}}
\end{align*}
as well as
\begin{equation*}
	t_K=\frac1{z^3+z^2-1-z^3(1+z)t_{K-1}} \quad\text{and}\quad \tau_K=\frac{1+z}{1-z^2-z^3-z^3\tau_{K-1}},\ \tau_0=0.
\end{equation*}
Now we move to the restricted Motzkin paths according to the three layers, computing them separately.
{\small
	\begin{equation*}
	\begin{pmatrix}
		-1+z^2 +2z^3+z^4& z(1+z) & 0&\cdots & &0  \\
		z^2 & z^3+z^2-1 &z(1+z)& \cdots && 0 \\
		&		z^2 & z^3+z^2-1 &z(1+z)& & 0 \\
		\vdots  & \ddots  & \ddots & \ddots  && \vdots\\
		0 & 0 & \cdots&&z^2 &z^3+z^2-1
	\end{pmatrix}
	\begin{pmatrix}
		f_0\\
		f_1\\
		f_2\\
		\vdots\\
		f_K
	\end{pmatrix}=
	\begin{pmatrix}
		-1\\
		z^2\\
		0\\
		\vdots\\
		0
	\end{pmatrix}
\end{equation*}
}
Therefore
\begin{align*}
f_0=\phi_{K+1}&=\frac{-\mD_{K}-z^3(1+z)\mD_{K-1}}{(-1+z^2 +2z^3+z^4)\mD_{K}-z^3(1+z)\mD_{K-1}}\\
&=1+\frac{z^2(1+z)^2}{1-z^2 -2z^3+z^3(1+z)t_K} ,\ \phi_0=1.
\end{align*}
Further, $g_0=\gamma_K=0$ for all $K$, so we can move to the instance of the third layer.
\begin{equation*}
	\begin{pmatrix}
		-1& z(1+z) & 0&\cdots & &0  \\
		z^2 & z^3+z^2-1 &z(1+z)& \cdots && 0 \\
		&		z^2 & z^3+z^2-1 &z(1+z)& & 0 \\
		\vdots  & \ddots  & \ddots & \ddots  && \vdots\\
		0 & 0 & \cdots&&z^2 &z^3+z^2-1
	\end{pmatrix}
	\begin{pmatrix}
		h_0\\
		h_1\\
		h_2\\
		\vdots\\
		h_K
	\end{pmatrix}=
	\begin{pmatrix}
		-z\\
		-z^2\\
		0\\
		\vdots\\
		0
	\end{pmatrix}
\end{equation*}
Therefore
\begin{align*}
h_0=\eta_{K+1}&=\frac{-z\mD_{K}+z^3(1+z)\mD_{K-1}}{-\mD_{K}-z^3\mD_{K-1}}\\
&=-1-z+\frac{1+2z}{1+z^3t_K},\ \eta_0=z.
\end{align*}

The generating function of Motzkin excursions  
\begin{equation*}
S(0)=1+z+z^2+3 z^3+4 z^4+7 z^5+15 z^6+26 z^7+50 z^8+102 z^9+196 z^{10}+\cdots
\end{equation*}
is A329666 in OEIS; the generating function of Motzkin meanders
\begin{equation*}
	S(1)=1+2z+3z^2+6z^3+11z^4+21z^5+42z^6+83z^7+167z^8+341z^9+697z^{10}+\cdots
\end{equation*}
is A329667  in OEIS.

%% file: UDHH.tex
\section{UD and HH are forbidden, A329676}
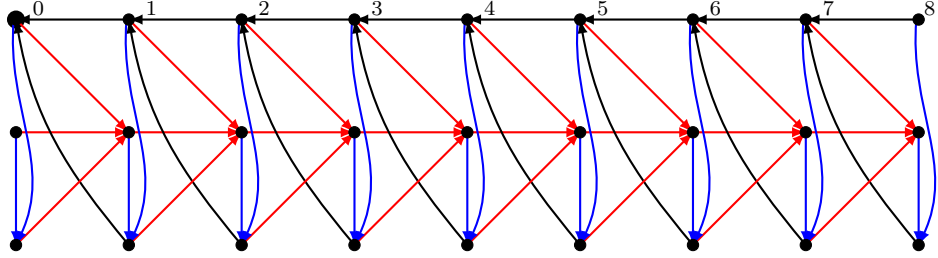
\begin{figure}[h]

	\begin{center}
		\begin{tikzpicture}[scale=1.5,main node/.style={circle,draw,font=\Large\bfseries}]

			\fill (0,0) circle (0.08cm);

			\foreach \x in {0,...,7}
			{
			}
			
			\foreach \x in {0,...,7}
			{
				
			}

			\foreach \x in {0,1,2,3,4,5,6,7}
			{
				\draw[thick, latex-] (\x,0) to  (\x+1,0);   
				\draw[ thick,red, -latex] (\x,-1)to  (\x+1,-1);    
				\draw[ thick,red, -latex] (\x,0)to  (\x+1,-1);         
				\node at  (\x+0.2,0.1){\tiny$\x$};
			}                        
			
			\foreach \x in {0,1,2,3,4,5,6,7,8}
			{
				\draw[thick, -latex, blue] (\x,-1) to  (\x,-2);         
				
			}                        
			
			\foreach \x in {1,2,3,4,5,6,7,8}
			{
				
				\draw[thick, -latex] [out=130,in=-80](\x,-2) to  (\x-1,0);    
				\draw[thick, -latex,red]  (\x-1,-2) to  (\x,-1);    
			}                        
			
			\foreach \x in {0,1,2,3,4,5,6,7,8}
			{
				
				\draw[thick, -latex,blue] [out=-95,in=65](\x-0.02,0) to  (\x,-2);    
				
			}                        
			
			\node at  (8+0.1,0.1){\tiny$8$};
			
			\foreach \x in {0,1,2,3,4,5,6,7,8}
			{
				\draw (\x,0) circle (0.05cm);
				\fill (\x,0) circle (0.05cm);
				\draw (\x,-1) circle (0.05cm);
				\fill (\x,-1) circle (0.05cm);
				\draw (\x,-2) circle (0.05cm);
				\fill (\x,-2) circle (0.05cm);
			}
		\end{tikzpicture}
	\end{center}
\caption{Graph (automaton) to recognize  Motzkin paths with forbidden subwords UD and HH.}
\label{purpelUDHH}
\end{figure}
We derive from Figure \ref*{purpelUDHH} the recursions
\begin{align*}
	f_i&=[i=0]+zf_{i+1}+zh_{i+1},\\
	g_{i+1}&=zf_{i}+zg_{i}+zh_{i},\\
	h_i&=zf_i+zg_i,
\end{align*}
and
\begin{align*}
	F(u)&=1+\frac zu(F(u)-F(0))+\frac zu(H(u)-H(0)),\\
	G(u)&=zuF(u)+zuG(u)+zuH(u),\\
	H(u)&=zF(u)+zG(u).
\end{align*}
and set $S=F+G+H$ and $	W:=\sqrt{(1+2z+3z^2+z^2)(1-2z-z^2+z^3)}$.
With
\begin{equation*}
	r_1=\frac{1+z^2+z^3-W}{2z(1+z)},\quad r_2=\frac{1+z^2+z^3+W}{2z(1+z)}
\end{equation*}
we have the factorization of the common numerator: $z(1+z)+(-1-z^3-z^2)u+zu^2(1+z)=z(1+z)(u-r_1)(u-r_2)$. Then
\begin{equation*}
	S(u)=\frac{(1+z)(zS(0)-u)}{z(1+z)(u-r_1)(u-r_2)}=\frac{zS(0)-u}{z(u-r_1)(u-r_2)}
\end{equation*}
Set $u=0$, divide out $u-r_1$, then $S(0)=\frac1{zr_2}$ and further
\begin{gather*}
	F(u)=\frac{(1-z(1+z)u)r_1}{z(1+z)(1-ur_1)},\ G(u)=\frac{ur_1}{1-ur_1},\ H(u)=\frac{r_1}{(1+z)(1-ur_1)}\\
	S(u)=\frac{r_1}{z(1-ur_1)}\Longrightarrow [u^j]S(u)=s_j= \frac{r_1^{j+1}}{z}.
\end{gather*}
Using the denominator, we find recursions
\begin{align*}
z(1+z)s_{j-2}-(1+z^2+z^3)s_{j-1}+z(1+z)s_{j}&=0,\ j\ge2,\\*
-s_0(1+z^2+z^3)+z(1+z)s_1&=-1-z.
\end{align*}
In matrix form:
{\small
	\begin{equation*}
		\begin{pmatrix}
			-z^3-z^2-1 & z(1+z) & 0&\cdots  &0  \\
			z(1+z) & -z^3-z^2-1 &z(1+z)& \cdots & 0 \\
			&		z(1+z) & -z^3-z^2-1 &z(1+z) & 0 \\
			\vdots  & \ddots  & \ddots & \ddots  & \vdots\\
			0 & 0 & \cdots&z(1+z) &-z^3-z^2-1
		\end{pmatrix}
		\begin{pmatrix}
			s_0\\
			s_1\\
			s_2\\
			\vdots\\
			s_K
		\end{pmatrix}=
		\begin{pmatrix}
			-1-z\\
			0\\
			0\\
			\vdots\\
			0
		\end{pmatrix}
	\end{equation*}
}

Perhaps a comment that we made before should be repeated here. While it is convenient to compute $s_0^{[K+1]}=f_0^{[K+1]}+g_0^{[K+1]}+h_0^{[K+1]}$
directly, it is slightly incorrect (and irrelevant for large $K$): One arrives in states from the first layer only with a $D$-step, 
in states from the second layer with a $U$-step, and in states from the third layer with an $H$-step. Consequently, the state $K+1$ in the first layer cannot
be reached, and the correct version would be $f_0^{[K]}+g_0^{[K+1]}+h_0^{[K+1]}$. Therefore we compute the respective quantities also independently.
We have  in matrix form:
{\footnotesize
\begin{equation*}
	\begin{pmatrix}
		z^4+z^3-1 & z(1+z) & 0&\cdots  &0  \\
		z(1+z) & -z^3-z^2-1 &z(1+z)& \cdots & 0 \\
		&		z(1+z) & -z^3-z^2-1 &z(1+z) & 0 \\
		\vdots  & \ddots  & \ddots & \ddots  & \vdots\\
		0 & 0 & \cdots&z(1+z) &-z^3-z^2-1
	\end{pmatrix}
	\begin{pmatrix}
		f_0\\
		f_1\\
		f_2\\
		\vdots\\
		f_K
	\end{pmatrix}=
	\begin{pmatrix}
		-1\\
		z(1+z)\\
		0\\
		\vdots\\
		0
	\end{pmatrix}
\end{equation*}
}

The quantity $g_0^{[K]}=0$ for all $K$, so we can move to the third layer.
Again, we have in matrix form
{\footnotesize
	\begin{equation*}
	\begin{pmatrix}
		-z^3-z^2-1 & z(1+z) & 0&\cdots  &0  \\
		z(1+z) & -z^3-z^2-1 &z(1+z)& \cdots & 0 \\
		&		z(1+z) & -z^3-z^2-1 &z(1+z) & 0 \\
		\vdots  & \ddots  & \ddots & \ddots  & \vdots\\
		0 & 0 & \cdots&z(1+z) &-z^3-z^2-1
	\end{pmatrix}
	\begin{pmatrix}
		h_0\\
		h_1\\
		h_2\\
		\vdots\\
		h_K
	\end{pmatrix}=
	\begin{pmatrix}
		-z\\
		0\\
		0\\
		\vdots\\
		0
	\end{pmatrix}
\end{equation*}
}

Let $\mathcal{D}_K$ be the determinant with $K$ rows and columns, after deleting the first row and column. Expanding, we derive the recursion
($\mathcal{D}_0=1$, $\mathcal{D}_1=-z^3-z^2-1$)
\begin{equation*}
	\mathcal{D}_K=(-z^3-z^2-1)\mathcal{D}_{K-1}-z^2(1+z)^2\mathcal{D}_{K-2}=
	\mathcal{D}_K=\frac{(-z(1+z))^{K+1}}W (r_1^{K+1}-r_2^{K+1}).
\end{equation*}
Using Cramer's rule as always,
\begin{align*}
	f_0&=\phi_{K+1}=\frac{-\mathcal{D}_K-z^2(1+z)^2\mathcal{D}_{K-1}}{(-1+z^3+z^4)\mathcal{D}_K-z^2(1+z)^2\mathcal{D}_{K-1}},
	\ \phi_0=1,\\*
	h_0&=\eta_{K+1}=\frac{-z\mathcal{D}_K}{(-z^3-z^2-1)\mathcal{D}_K-z^2(1+z)^2\mathcal{D}_{K-1}},\ \eta_0=z,\\*
	s_0&=\sigma_{K+1}=\frac{-(1+z)\mathcal{D}_K}{(-z^3-z^2-1)\mathcal{D}_K-z^2(1+z)^2\mathcal{D}_{K-1}},\ \sigma_0=1+z.
\end{align*}

Finally we move to expansions of the continuous fraction type. From the recursion for the determinants, we have, with
$t_K=\mD_{K-1}/\mD_{K}$ and $\tau_K=-z^2(1+z)^2t_K$,
\begin{equation*}
	{t_K}=\frac{-1}{1+z^2+z^3+z^2(1+z)^2t_{K-1}}, \quad - z^2(1+z)^2{t_K}=\frac{z^2(1+z)^2}{1+z^2+z^3+z^2(1+z)^2t_{K-1}},
\end{equation*}
and
\begin{equation*}
\tau_K=\frac{z^2(1+z)^2}{1+z^2+z^3-\tau_{K-1}},\quad \tau_0=0.
\end{equation*}
Furthermore,
\begin{align*}
	\sigma_{K+1}&=\frac{-(1+z)\mathcal{D}_K}{(-z^3-z^2-1)\mathcal{D}_K-z^2(1+z)^2\mathcal{D}_{K-1}}=	\frac{1+z}{1+z^2+z^3-\tau_{K}},\\
\phi_{K+1}&=\frac{-\mathcal{D}_K-z^2(1+z)^2\mathcal{D}_{K-1}}{(-1+z^3+z^4)\mathcal{D}_K-z^2(1+z)^2\mathcal{D}_{K-1}}\\&=
\frac{-z^2(1+z)^2t_{K}}{(-1+z^3+z^4)-z^2(1+z)^2t_{K}}=1-\frac{1-z^3-z^4}{1-z^3-z^4-\tau_K},\\
\sigma_{K+1}&=\frac{-(1+z)}{(-z^3-z^2-1)-z^2(1+z)^2t_K}=\frac{1+z}{1+z^2+z^3-\tau_K}.
\end{align*}
The generating function
\begin{equation*}
S(0)=1+z+z^3+2 z^4+2 z^5+5 z^6+10 z^7+16 z^8+34 z^9+68 z^{10}+\cdots
\end{equation*}
is in OEIS as sequence A329676; the sequence
\begin{equation*}
	S(1)=1+z+z^3+2 z^4+2 z^5+5 z^6+10 z^7+16 z^8+34 z^9+68 z^{10}+\cdots
\end{equation*}
is sequence A329675 in OEIS.

%% file: DUHH.tex
\section{DU and HH are forbidden, A329666}
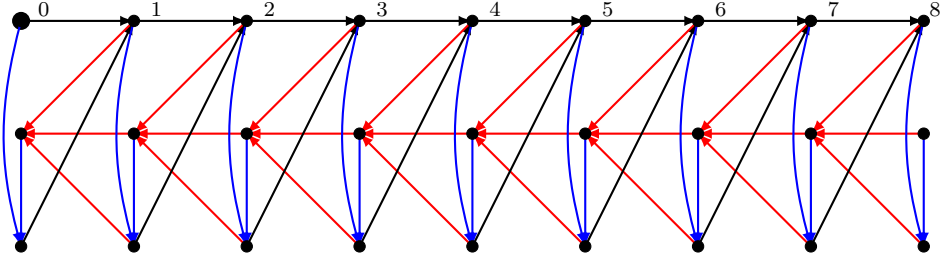
\begin{figure}[h]

	\begin{center}
		\begin{tikzpicture}[scale=1.5,main node/.style={circle,draw,font=\Large\bfseries}]

			\fill (0,0) circle (0.08cm);

			\foreach \x in {0,...,7}
			{
			}
			
			\foreach \x in {0,...,7}
			{
				
			}

			\foreach \x in {0,1,2,3,4,5,6,7}
			{
				\draw[thick, -latex] (\x,0) to  (\x+1,0);   
				\draw[ thick,red, -latex] (\x+1,-1)to  (\x,-1);    
				\draw[ thick,red, -latex] (\x+1,-2)to  (\x,-1);         
				\node at  (\x+0.2,0.1){\tiny$\x$};
			}                        
			
			\foreach \x in {0,1,2,3,4,5,6,7,8}
			{
				\draw[thick, -latex, blue] (\x,-1) to  (\x,-2);         
				
			}                        
			
			\foreach \x in {1,2,3,4,5,6,7,8}
			{
				
				\draw[thick, -latex] (\x-1,-2) to  (\x,0);    
				\draw[thick, -latex,red]  (\x,0) to  (\x-1,-1);    
			}                        
			
			\foreach \x in {0,1,2,3,4,5,6,7,8}
			{
				
				\draw[thick, -latex,blue] [out=-105,in=105](\x-0.00,0) to  (\x,-2);    
				
			}                        
			
			\node at  (8+0.1,0.1){\tiny$8$};
			
			\foreach \x in {0,1,2,3,4,5,6,7,8}
			{
				\draw (\x,0) circle (0.05cm);
				\fill (\x,0) circle (0.05cm);
				\draw (\x,-1) circle (0.05cm);
				\fill (\x,-1) circle (0.05cm);
				\draw (\x,-2) circle (0.05cm);
				\fill (\x,-2) circle (0.05cm);
			}
		\end{tikzpicture}
	\end{center}
\caption{Graph (automaton) to recognize  Motzkin paths with forbidden subwords DU and HH.}
\label{purpelDUHH}
\end{figure}
\begin{align*}
	f_{i+1}&=zf_i+zh_i,\quad f_0=1,\\
	g_i&=zf_{i+1}+zg_{i+1}+zh_{i+1},\\
	h_i&=zf_i+zg_i
\end{align*}
\begin{align*}
	F(u)&=1+zuF(u)+zuH(u),\\
	G(u)&=\frac{z}{u}\big(F(u)-F(0)\bigr)+\frac{z}{u}\big(G(u)-G(0)\bigr)+\frac{z}{u}\big(H(u)-H(0)\bigr),\\
	H(u)&=zF(u)+zG(z)
\end{align*}
Decomposing the common denominator as 
$$z(z+1)+(-z^3-1-z^2)u+z(z+1)u^2=z(z+1)(u-r_1)(u-r_2) 	  $$
with 
$$ W=\sqrt{(1+2z+3z^2+z^3)(1-2z-z^2+z^3)}$$
and
\begin{equation*}
	r_1=\frac{1+z^2+z^3-W}{2z(1+z)},\quad r_2=\frac{1+z^2+z^3+W}{2z(1+z)},
\end{equation*}
we obtain by dividing out $u-r_1$ and simplifying
\begin{gather*}
		F(u)=\frac1{1-ur_1},\ G(u)=\frac{ r_1-z(1+z)}{z^2(1-ur_1)},\ H(u)=\frac{ r_1-z}{z(1-ur_1)},\ S(u)=\frac{ (1+z)(r_1-z)}{z^2(1-ur_1)}.
\end{gather*}
This leads to (excursions)
\begin{equation*}
	S(0)=1+z+z^2+3z^3+4z^4+7z^5+15z^6+26z^7+50z^8+102z^9+196z^{10}+392z^{11}+\cdots,
\end{equation*}
which is sequence A329666. Further
\begin{equation*}
	S(1)=1+2 z+4 z^2+9 z^3+18 z^4+38 z^5+81 z^6+171 z^7+366 z^8+787 z^9+1693 z^{10}+\cdots,
\end{equation*}
and this generating function of meanders is  sequence A329668 in OEIS.
And we get the generating function of Motzkin excursions ending at level $j$:
\begin{equation*}
s_j=[u^j]\frac{ (1+z)(r_1-z)}{z^2(1-ur_1)}=\frac{ (1+z)(r_1-z)}{z^2}r_1^j.
\end{equation*}

Now we move to paths of bounded height (indices $>K$ cannot occur)
\begin{align*}
	z(1+z)s_1-s_0&=-1-z\\
	z(1+z)s_j-(1+z^2+z^3)s_{j-1}+z(1+z)s_{j-2}&=0,\ j\ge2,
\end{align*}
or in Matrix form
{\small
\begin{equation*}
	\begin{pmatrix}
		-1 & z(1+z) & 0&\cdots  &0  \\
		z(1+z) & -z^3-z^2-1 &z(1+z)& \cdots & 0 \\
		&		z(1+z) & -z^3-z^2-1 &z(1+z) & 0 \\
		\vdots  & \ddots  & \ddots & \ddots  & \vdots\\
		0 & 0 & \cdots&z(1+z) &-z^3-z^2-1
	\end{pmatrix}
	\begin{pmatrix}
		s_0\\
		s_1\\
		s_2\\
		\vdots\\
		s_K
	\end{pmatrix}=
	\begin{pmatrix}
		-1-z\\
		0\\
		0\\
		\vdots\\
		0
	\end{pmatrix}
\end{equation*}
}

Let, as always denote the determinant, after removing the first row resp.\ column by $\mD_K$, then
\begin{equation*}
	\mathcal{D}_K=(-z^3-z^2-1)\mathcal{D}_{K-1}-z^2(1+z)^2\mathcal{D}_{K-2},\ \mathcal{D}_0=1, \mathcal{D}_1=-z^3-z^2-1.
\end{equation*}
From this, with $t_K=\mD_{K-1}/\mD_{K}$ and $\tau_K=-z^2(1+z)^2t_{K}$
\begin{equation*}
	t_K=\frac{1}{-z^3-z^2-1-z^2(1+z)^2t_{K-1}},\quad \tau_K=\frac{z^2(1+z)^2}{1+z^2+z^3-\tau_{K-1}},\ \tau_0=0.
\end{equation*}
Of course, $\mD_K$ has also a Binet form:
\begin{equation*}
\mD_K=\frac{(-z(z+1))^{K+1}}{W}(r_1^{K+1}-r_2^{K+1})
\end{equation*}
\begin{equation*}
	\tau_i:=-z^2(1+z)^2\frac{\mathcal{D}_{i-1}}{\mathcal{D}_i},\quad
	\tau_i=\frac{z^2(1+z)^2}{1+z^2+z^3-\tau_{i-1}}, \ \tau_0=0
\end{equation*}
Let us also consider the  second and third layer in matrix form ($f_0^{[K]}=1$) for any $K$:
{\footnotesize
\begin{equation*}
	\begin{pmatrix}
		-1+z^3+z^4 & z(1+z) & 0&\cdots  &0  \\
		z(1+z) & -z^3-z^2-1 &z(1+z)& \cdots & 0 \\
		&		z(1+z) & -z^3-z^2-1 &z(1+z) & 0 \\
		\vdots  & \ddots  & \ddots & \ddots  & \vdots\\
		0 & 0 & \cdots&z(1+z) &-z^3-z^2-1
	\end{pmatrix}
	\begin{pmatrix}
		g_0\\
		g_1\\
		g_2\\
		\vdots\\
		g_K
	\end{pmatrix}=
	\begin{pmatrix}
		-z^2(z+1)^2\\
		0\\
		0\\
		\vdots\\
		0
	\end{pmatrix}
\end{equation*}
}
{\footnotesize
	\begin{equation*}
		\begin{pmatrix}
			-1 & z(1+z) & 0&\cdots  &0  \\
			z(1+z) & -z^3-z^2-1 &z(1+z)& \cdots & 0 \\
			&		z(1+z) & -z^3-z^2-1 &z(1+z) & 0 \\
			\vdots  & \ddots  & \ddots & \ddots  & \vdots\\
			0 & 0 & \cdots&z(1+z) &-z^3-z^2-1
		\end{pmatrix}
		\begin{pmatrix}
			h_0\\
			h_1\\
			h_2\\
			\vdots\\
			h_K
		\end{pmatrix}=
		\begin{pmatrix}
			-z\\
			0\\
			0\\
			\vdots\\
			0
		\end{pmatrix}
	\end{equation*}
}
By Cramers method:
\begin{align*}
s_0=\sigma_{K+1}&=\frac{-(1+z)\mD_{K}}{-\mD_{K}-z^2(1+z)^2\mD_{K-1}},\ \sigma_0=1+z\\
g_0=\gamma_{K+1}&=\frac{-z^2(1+z)^2\mD_{K}}{(-1+z^3+z^4)\mD_{K}-z^2(1+z)^2\mD_{K-1}},\ \gamma_0=0 	\\
h_0=\eta_{K+1}&=\frac{-z\mD_{K}}{-\mD_{K}-z^2(1+z)^2\mD_{K-1}},\ \eta_0=z
\end{align*}
Let us rewrite it in continued fraction style:
\begin{align*}
\sigma_{K+1}&=\frac{1+z}{1+z^2(1+z)^2t_{K}}= \frac{1+z}{1-\tau_{K}},\\
\gamma_{K+1}&=\frac{-z^2(1+z)^2}{(-1+z^3+z^4)-z^2(1+z)^2t_{K-1}}=\frac{z^2(1+z)^2}{1-z^3-z^4-\tau_{K}},\\
\eta_{K+1}&=\frac{-z}{-1-z^2(1+z)^2t_{K-1}}=\frac{z}{1-\tau_{K-1}}.
\end{align*}

%% file: DDHH.tex
\section{DD and HH are forbidden, A329666}


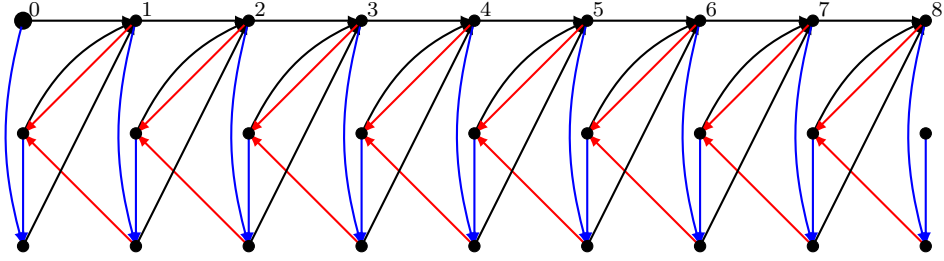
\begin{figure}[h]

	\begin{center}
		\begin{tikzpicture}[scale=1.5,main node/.style={circle,draw,font=\Large\bfseries}]

			\fill (0,0) circle (0.08cm);

			\foreach \x in {0,...,7}
			{
			}
			
			\foreach \x in {0,...,7}
			{
				
			}

			\foreach \x in {0,1,2,3,4,5,6,7}
			{
				\draw[thick, -latex] (\x,0) to  (\x+1,0);   
				\draw[ thick, -latex] (\x,-1)to [out=65,in=205] (\x+1,0);    
				\draw[ thick,red, -latex] (\x+1,-2)to  (\x,-1);         
				\node at  (\x+0.1,0.1){\tiny$\x$};
			}                        
			
			\foreach \x in {0,1,2,3,4,5,6,7,8}
			{
				\draw[thick, -latex, blue] (\x,-1) to  (\x,-2);         
				
			}                        
			
			\foreach \x in {1,2,3,4,5,6,7,8}
			{
				
				\draw[thick, -latex] (\x-1,-2) to  (\x,0);    
				\draw[thick, -latex,red]  (\x,0) to  (\x-1,-1);    
			}                        
			
			\foreach \x in {0,1,2,3,4,5,6,7,8}
			{
				
				\draw[thick, -latex,blue] [out=-105,in=105](\x-0.00,0) to  (\x,-2);    
				
			}                        
			
			\node at  (8+0.1,0.1){\tiny$8$};
			
			\foreach \x in {0,1,2,3,4,5,6,7,8}
			{
				\draw (\x,0) circle (0.05cm);
				\fill (\x,0) circle (0.05cm);
				\draw (\x,-1) circle (0.05cm);
				\fill (\x,-1) circle (0.05cm);
				\draw (\x,-2) circle (0.05cm);
				\fill (\x,-2) circle (0.05cm);
			}
		\end{tikzpicture}
	\end{center}
\caption{Graph (automaton) to recognize  Motzkin paths with forbidden subwords DD and HH.}
\label{purpelDDHH}
\end{figure}
From Figure~\ref{purpelDDHH} we derive the recursions
\begin{align*}
	f_{i+1}&=zf_i+zg_i+zh_i,\quad f_0=1,\\
	g_i&=zf_{i+1} +zh_{i+1},\\
	h_i&=zf_i+zg_i,
\end{align*}
and
\begin{align*}
	F(u)&=1+zuF(u)+zuG(u)+zuH(u),\\
	G(u)&=\frac{z}{u}[F(u)-F(0)] 	+\frac{z}{u}[H(u)-H(0)],\\
	H(u)&=zF(u)+zG(z).
\end{align*}
With $W=\sqrt{(1+2z+3z^2+z^3)(1-2z-z^2+z^3)}$ and
\begin{equation*}
	r_1=\frac{1-z^2-z^3-W}{2z(1+z)},\quad r_2=\frac{1-z^2-z^3+W}{2z(1+z)},
\end{equation*}
the denominator (after solving the system) factors:
\begin{equation*}
z^2+(z^2-1+z^3)u+z(1+z)u^2=z(1+z)(u-r_1)(u-r_2)
\end{equation*}
After diving out $u-r_1$ and simplifying,
\begin{gather*}
	F(u)=\frac{z}{z-u(1+z)r_1},\ G(u)=\frac{r_1-z^2}{z(z-u(1+z)r_1)},\	H(u)=\frac{r_1}{z-u(1+z)r_1},\\ 
	S(u)=\frac{(1+z)r_1}{z(z-u(1+z)r_1)}.
\end{gather*}
We can read off to get the generating function of Motzkin excursions ending at level $j$:
\begin{equation*}
	[u^j]	S(u)=[u^j]\frac{(1+z)r_1}{z^2(1-\frac{(1+z)r_1u}{z})}=\frac{(1+z)r_1}{z^2}\Bigl(\frac{(1+z)r_1}{z}\Big)^j=\frac{(1+z)^{j+1}r_1^{j+1}}{z^{j+2}}.
\end{equation*}
The generating function
\begin{equation*}
S(0)=1+z+z^2+3z^3+4z^4+7z^5+15z^6+26z^7+50z^8+102z^9+196z^{10}+\cdots
\end{equation*}
is A329666 in OEIS; 
\begin{equation*}
	S(1)=1+2z+4z^2+10z^3+23z^4+54z^5+129z^6+307z^7+733z^8+1757z^9+4213z^{10}+\cdots
\end{equation*}
is A329669 in OEIS.

Now we move to restricted paths (indices $>K$ are not allowed):
\begin{align*}
	z(1+z)s_{j-2}+(-1+z^2+z^3)s_{j-1}+z^2s_j&=0,\ j\ge2,\\* (-1+z^2+z^3)s_0 +z^2s_1&=-(1+z),
\end{align*}
or in Matrix form
\begin{equation*}
	\begin{pmatrix}
		z^3+z^2-1 & z^2 & 0&\cdots  &0  \\
		z(1+z) & z^3+z^2-1 &z^2& \cdots & 0 \\
		&		z(1+z) & z^3+z^2-1 &z^2 & 0 \\
		\vdots  & \ddots  & \ddots & \ddots  & \vdots\\
		0 & 0 & \cdots&z(1+z) &z^3+z^2-1
	\end{pmatrix}
	\begin{pmatrix}
		s_0\\
		s_1\\
		s_2\\
		\vdots\\
		s_K
	\end{pmatrix}=
	\begin{pmatrix}
		-1-z\\
		0\\
		0\\
		\vdots\\
		0
	\end{pmatrix}
\end{equation*}
Now $f_0=\phi_{K}=1$ for all $K$, so we can move to the next layer:
\begin{align*}
	z(1+z)g_{j-2}+(-1+z^2+z^3)g_{j-1}+z^2g_j&=0,\ j\ge2,\\* (-1+z^2+2z^3+z^4)g_0+z^2g_1&=-(1+z)^2z^2.
\end{align*}
\begin{equation*}
	{\footnotesize
	\begin{pmatrix}
		-1+z^2+2z^3+z^4 & z^2 & 0&\cdots  &0  \\
		z(1+z) & z^3+z^2-1 &z^2& \cdots & 0 \\
		&		z(1+z) & z^3+z^2-1 &z^2 & 0 \\
		\vdots  & \ddots  & \ddots & \ddots  & \vdots\\
		0 & 0 & \cdots&z(1+z) &z^3+z^2-1
	\end{pmatrix}
	\begin{pmatrix}
		g_0\\
		g_1\\
		g_2\\
		\vdots\\
		g_K
	\end{pmatrix}=
	\begin{pmatrix}
		-z^2(1+z)^2\\
		0\\
		0\\
		\vdots\\
		0
	\end{pmatrix}
}
\end{equation*}
Finally
\begin{align*}
	z(1+z)h_{j-2}+(-1+z^2+z^3)h_{j-1}+z^2h_j&=0,\ j\ge2,\\* (-1+z^2+z^3)h_0+z^2h_1&=-z.
\end{align*}
\begin{equation*}
	\begin{pmatrix}
		z^3+z^2-1 & z^2 & 0&\cdots  &0  \\
		z(1+z) & z^3+z^2-1 &z^2& \cdots & 0 \\
		&		z(1+z) & z^3+z^2-1 &z^2 & 0 \\
		\vdots  & \ddots  & \ddots & \ddots  & \vdots\\
		0 & 0 & \cdots&z(1+z) &z^3+z^2-1
	\end{pmatrix}
	\begin{pmatrix}
		h_0\\
		h_1\\
		h_2\\
		\vdots\\
		h_K
	\end{pmatrix}=
	\begin{pmatrix}
		-z\\
		0\\
		0\\
		\vdots\\
		0
	\end{pmatrix}
\end{equation*}

Finally we move to continued fractions, with the usual determinant, after deleting the first row resp.\ column.
\begin{equation*}
	\mathcal{D}_K=(z^3+z^2-1)\mathcal{D}_{K-1}-z^3(1+z)\mathcal{D}_{K-2}
\end{equation*}
and, with $t_K=\mD_{K-1}/\mD_{K}$ and $\tau_K=-z^3(1+z)t_K$
\begin{equation*}
t_K=\frac{1}{-1+z^2+z^3-z^3(1+z)t_{K-1}},\quad \tau_K=\frac{z^3(1+z)}{1-z^2-z^3-\tau_{K-1}},\ \tau_0=0.
\end{equation*}
Using this, we finally get
\begin{align*}
	s_0=\sigma_{K+1}&=\frac{-(1+z)\mathcal{D}_{K}}{(-1+z^2+z^3)\mathcal{D}_{K}-z^3(1+z)\mathcal{D}_{K-1}}\\
	&=\frac{1+z}{1-z^2-z^3-\tau_K},\quad \sigma_0=1+z,\\
	g_0=\gamma_{K+1}&=	\frac{-z^2(1+z)^2\mathcal{D}_{K}}{(-1+z^2+2z^3+z^4)\mathcal{D}_{K}-z^3(1+z)\mathcal{D}_{K-1}}\\
		&= \frac{z^2(1+z)^2}{1-z^2-2z^3-z^4-\tau_{K}},\gamma_0=0,\\
	h_0=\eta_{K+1}&=	\frac{-z\mathcal{D}_{K}}{(-1+z^2+z^3)\mathcal{D}_{K}-z^3(1+z)\mathcal{D}_{K-1}}\\
&= \frac{z}{1-z^2-z^3-\tau_{K}},\quad \eta_0=z.
		\end{align*}

%% file: UHHD.tex
\section{UH and HD are forbidden, A329702 }

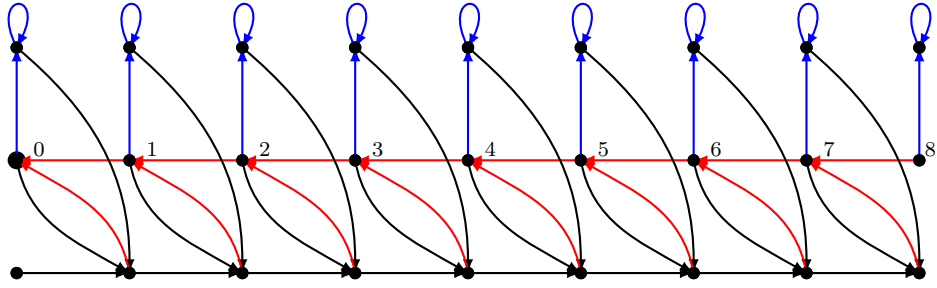
\begin{figure}[h]

	\begin{center}
		\begin{tikzpicture}[scale=1.5,main node/.style={circle,draw,font=\Large\bfseries}]

			\foreach \x in {0,1,2,3,4,5,6,7,8}
			{
				\draw (\x,0) circle (0.05cm);
				\fill (\x,0) circle (0.05cm);
				\draw (\x,-1) circle (0.05cm);
				\fill (\x,-1) circle (0.05cm);
				\draw (\x,1) circle (0.05cm);
				\fill (\x,1) circle (0.05cm);
			}

			\fill (0,0) circle (0.08cm);

			\foreach \x in {0,...,8}
			{
				\draw[thick,-latex,blue ] (\x,1)  ..  controls (\x-0.25,1+0.5) and  (\x+0.25,1+0.5) .. (\x,1+0) ;	
				\draw[thick,-latex,blue ] (\x,0)to (\x,1+0) ;	
			}
			
			\foreach \x in {0,...,7}
			{
				\draw[thick,latex-,red ] (\x,0)  to (\x+1,0) ;	
			}
			
			\foreach \x in {0,...,7}
			{
				
				\draw[thick,latex- ] (\x+1,-1)[out=150,in=-80] to (\x,0) ;	
				\draw[thick,-latex ] (\x,1) [out=-40,in=90]to (\x+1,-1) ;	
			}

			\foreach \x in {0,1,2,3,4,5,6,7}
			{
				\draw[thick, - latex] (\x,-1) to  (\x+1,-1);	
				\draw[thick,red, latex-] (\x,0) [out=-30,in=100]to  (\x+1,-1);	
				\node at  (\x+0.2,0.1){\tiny$\x$};
			}			
			
			\foreach \x in {0,1,2,3,4,5,6,7,8}
			{
			}

			\node at  (8+0.1,0.1){\tiny$8$};
			
			\foreach \x in {0,1,2,3,4,5,6,7,8}
			{
				\draw (\x,0) circle (0.05cm);
				\fill (\x,0) circle (0.05cm);
				\draw (\x,-1) circle (0.05cm);
				\fill (\x,-1) circle (0.05cm);
				\draw (\x,1) circle (0.05cm);
				\fill (\x,1) circle (0.05cm);
			}
		\end{tikzpicture}
	\end{center}
\caption{Graph (automaton) to recognize  Motzkin paths with forbidden subwords UH and HD.}
\label{purpelUHHD}

\end{figure}
The Figure~\ref{purpelUHHD} leads to the recursions
\begin{align*}
	f_i&=zf_i+zg_i,\\
	g_i&= [i=0]+zg_{i+1}+zh_{i+1},\ \\
	h_{i+1}&=zf_i+zg_i+zh_i, \ h_0=0
\end{align*}
and functional equations
\begin{align*}
	F(u)&=zF(u)+zG(u) \Rightarrow F(u)=\frac{zG(u)}{1-z}\\
	G(u)&=1+\frac zu(G(u)-G(0))+\frac zu H(u),\\
	H(u)&=zuF(u)+zuG(u)+zuH(u).
\end{align*}
The denominator of the solution(s) is $z(z-1)+(-z+1-z^3)u+z(z-1)u^2=z(z-1)(u-r_1)(u-r_2)$ with
\begin{equation*}
	r_1=\frac{1-z-z^3-W}{2z(1-z)},\ r_2=\frac{1-z-z^3+W}{2z(1-z)},\ W=\sqrt{(1+z-2z^2-z^3)(1-3z+2z^2-z^3)}
\end{equation*}
Dividing out $u-r_1$ and simplifying,
\begin{gather*}
	F(u)=\frac{(1-zu)r_1}{(1-z)(1-ur_1)},\ 
	G(u)= \frac{(1-zu)r_1}{z(1-ur_1)},\ 
	H(u)= \frac{ur_1}{(1-z)(1-ur_1)},\\ 
	S(u)=\frac{r_1}{z(1-z)(1-ur_1)}. 
\end{gather*}
The series
\begin{equation*}
S(0)=1+z+2z^2+3z^3+6z^4+10z^5+20z^6+36z^7+73z^8+139z^9+286z^{10}+567z^{11}+\cdots
\end{equation*}
is A329702 in OEIS
\begin{equation*}
S(1)=1+2z+4z^2+8z^3+17z^4+36z^5+78z^6+169z^7+371z^8+816z^9+1809z^{10}+\cdots
\end{equation*}
is not in OEIS. Also,
\begin{equation*}
s_j=[u^j]S(u)=[u^j]\frac{r_1}{z(1-z)(1-ur_1)}=\frac{r_1^{j+1}}{z(1-z)},
\end{equation*}
which is the generating function of restricted Motzkin paths ending at level $j$.

From the denominator, we derive recursions
\begin{align*}
	z(z-1)s_j +(1-z-z^3)s_{j-1}+z(z-1)s_{j-2}&=0,\ j\ge2\\
	(1-z-z^3)s_0+z(z-1)s_1&=1,
\end{align*}
which is in matrix form
\begin{equation*}
	\begin{pmatrix}
		1-z-z^3 & z(z-1) & 0&\cdots  &0  \\
		z(z-1) & 1-z-z^3 &z(z-1)& \cdots & 0 \\
		&		z(z-1) & 1-z-z^3 &z(z-1) & 0 \\
		\vdots  & \ddots  & \ddots & \ddots  & \vdots\\
		0 & 0 & \cdots&z(z-1) &1-z-z^3
	\end{pmatrix}
	\begin{pmatrix}
		s_0\\
		s_1\\
		s_2\\
		\vdots\\
		s_K
	\end{pmatrix}=
	\begin{pmatrix}
		1\\
		0\\
		0\\
		\vdots\\
		0
	\end{pmatrix}
\end{equation*}
The recursions for the first resp.\ second layer in matrix form are
{\small
\begin{equation*}
	\begin{pmatrix}
		1-z-z^2 & z(z-1) & 0&\cdots  &0  \\
		z(z-1) & 1-z-z^3 &z(z-1)& \cdots & 0 \\
		&		z(z-1) & 1-z-z^3 &z(z-1) & 0 \\
		\vdots  & \ddots  & \ddots & \ddots  & \vdots\\
		0 & 0 & \cdots&z(z-1) &1-z-z^3
	\end{pmatrix}
	\begin{pmatrix}
		f_0\\
		f_1\\
		f_2\\
		\vdots\\
		f_K
	\end{pmatrix}=
	\begin{pmatrix}
		z\\
		-z^2\\
		0\\
		\vdots\\
		0
	\end{pmatrix}
\end{equation*}
}
{\small
\begin{equation*}
	\begin{pmatrix}
		1-z-z^2 & z(z-1) & 0&\cdots  &0  \\
		z(z-1) & 1-z-z^3 &z(z-1)& \cdots & 0 \\
		&		z(z-1) & 1-z-z^3 &z(z-1) & 0 \\
		\vdots  & \ddots  & \ddots & \ddots  & \vdots\\
		0 & 0 & \cdots&z(z-1) &1-z-z^3
	\end{pmatrix}
	\begin{pmatrix}
		g_0\\
		g_1\\
		g_2\\
		\vdots\\
		g_K
	\end{pmatrix}=
	\begin{pmatrix}
		1-z\\
		z(z-1)\\
		0\\
		\vdots\\
		0
	\end{pmatrix}
\end{equation*}}

For the third layer we get $h_0^{[K]}=0$ for all $K$ without computations.
Let $\mathcal{D}_K$ be the usual determinant. Then
\begin{align*}
	\mathcal{D}_K&=(1-z-z^3)\mathcal{D}_{K-1}-z^2(z-1)^2\mathcal{D}_{K-2},\ \mathcal{D}_0=1,\ \mathcal{D}_1=1-z-z^3
	\\&=\frac{z^{K+1}(1-z)^{K+1}}W(r_2^{K+1}-r_1^{K+1}).
\end{align*}
Further, with $t_K=\mD_{K-1}/\mD_{K}$, $\tau_K=z^2(z-1)^2t_{K}$
\begin{equation*}
	t_K=\frac1{1-z-z^3-z^2(z-1)^2t_{K-1}},\quad \tau_K=\frac{z^2(z-1)^2}{1-z-z^3-\tau_{K-1}},\ \tau=0.
\end{equation*}
Then we can express the solutions in continued fraction form:
\begin{align*}
	s_0=\sigma_{K+1}&=\frac{(1-z-z^3)\mD_K}{(1-z-z^3)\mD_K-z^2(z-1)^2\mD_{K-1}}=
	\frac{1-z-z^3}{1-z-z^3-\tau_{K}},\\
	f_0=\phi_{K+1}&=\frac{z\mD_K+z^3(z-1)\mD_{K-1}}{(1-z-z^2)\mD_K-z^2(z-1)^2\mD_{K-1}}\\
	&=\frac{z+z^3(z-1)t_{K}}{1-z-z^2-z^2(z-1)^2t_{K}}=\frac{z}{1-z}+\frac{z^3}{(1-z)(1-z-z^2-\tau_K)},\\
	g_0=\gamma_{K+1}&=\frac{(1-z)\mD_K-z^2(z-1)^2\mD_{K-1}}{(1-z-z^2)\mD_K-z^2(z-1)^2\mD_{K-1}}\\
	&=\frac{1-z-z^2(z-1)^2t_{K}}{1-z-z^2-z^2(z-1)^2t_{K}}=1+\frac{z^2}{1-z-z^2-\tau_K}.
\end{align*}

%% file: DHHU.tex
\section{DH and HU are forbidden, A329701  }

\begin{figure}[h]

	\begin{center}
		\begin{tikzpicture}[scale=1.5,main node/.style={circle,draw,font=\Large\bfseries}]


			\foreach \x in {0,...,8}
			{
				\draw[thick,-latex,blue ] (\x,1)  ..  controls (\x-0.25,1+0.5) and  (\x+0.25,1+0.5) .. (\x,1+0) ;	
				\draw[thick,-latex,blue ] (\x,0)  to (\x,1+0) ;	
				
			}
			
			\foreach \x in {0,...,7}
			{
				\draw[thick,latex-,red ] (\x,-1)  to[in=-150,out=50] (\x+1,0) ;	
			}
			
			\foreach \x in {0,...,7}
			{
				
				\draw[thick,latex- ] (\x+1,0) to (\x,0) ;	
								\draw[thick,latex- ] (\x+1,0) to[in=20,out=-110] (\x,-1) ;	
				\draw[thick,-latex , red] (\x+1,1) to [in=70,out=-130](\x,-1) ;	
			}

			\foreach \x in {0,1,2,3,4,5,6,7}
			{
				\draw[thick,red, latex-] (\x,-1) to  (\x+1,-1);	
				\node at  (\x+0.2,0.1){\tiny$\x$};
			}			
			
			\foreach \x in {0,1,2,3,4,5,6,7,8}
			{
			}			
			

			
			\node at  (8+0.1,0.1){\tiny$8$};

			\draw (0,0) circle (0.07cm);
			\fill (0,0) circle (0.07cm);
			\foreach \x in {1,2,3,4,5,6,7,8}
			{
				\draw (\x,0) circle (0.05cm);
				\fill (\x,0) circle (0.05cm);
				\draw (\x,-1) circle (0.05cm);
				\fill (\x,-1) circle (0.05cm);
				\draw (\x,1) circle (0.05cm);
				\fill (\x,1) circle (0.05cm);
			}
		\end{tikzpicture}
	\end{center}
	\caption{Graph (automaton) to recognize  Motzkin paths with forbidden subwords DH and HU.}
	\label{purpelUUHH}
\end{figure}
 The corresponding recursions are
\begin{align*}
	f_j&=zf_j+zg_j,\\
	g_{j+1}&=zg_{j}+zh_{j},\ g_0=1,\\
	h_{j}&=zf_{j+1}+zg_{j+1}+zh_{j+1},
\end{align*}
and further
\begin{align*}
	F(u)&=zF(u)+zG(u),\\
G(u)&=1+zuG(u)+zuH(u),\\
	H(u)&=\frac zu\bigl(F(u)-F(0)\bigr)+\frac zu\bigl(G(u)-1\bigr)+\frac zu\bigl(H(u)-H(0)\bigr), 
\end{align*}
Solving the system there is a denominator which factors:
\begin{equation*}
z(z-1)+(-z^3-z+1)u+z(z-1)u^2=z(z-1)(u-r_1)(u-r_2)
\end{equation*}
with
\begin{gather*}
	r_1=\frac{1-z+z^3-W}{2z(1-z)},\quad r_2=\frac{1-z+z^3+W}{2z(1-z)},\\ W=\sqrt{(1-3z+2z^2-z^3) (1+z -2z^2 -z^3)}.
\end{gather*}
Dividing out the factor $u-r_1$ and simplifying, we find
\begin{gather*}
	F=\frac{z}{(1-z)(1-ur_1)},\ 	G=\frac{1}{1-ur_1},\ 	H=\frac{r1-z}{z(1-ur_1)},\\ 	S=\frac{(1-z)r_1+z^2}{z(1-z)(1-ur_1)}. 
	\end{gather*}
The generating function of excursions
\begin{equation*}
S(0)=1+z+2 z^2+2 z^3+4 z^4+5 z^5+11 z^6+17 z^7+38 z^8+67 z^9+148 z^{10}+\cdots
\end{equation*}
the generating function of meanders is
\begin{equation*}
	S(1)=1+2z+4z^2+7z^3+14z^4+28z^5+60z^6+128z^7+281z^8+615z^9+1365z^{10}+\cdots
\end{equation*}
is not in OEIS. We can also compute the generating function of Motzkin paths ending at level $j$:
\begin{equation*}
[u^j]S(u)=[u^j]\frac{(1-z)r_1+z^2}{z(1-z)(1-ur_1)}=\frac{(1-z)r_1+z^2}{z(1-z)}r_1^j=\frac{r_1^{j+1}}{z}+\frac{zr_1^j}{1-z}.
\end{equation*}
Now we move to height restricted paths:   
\begin{align*}
	-z(1-z)s_{j-2}+(1-z-z^3)s_{j-1}-z(1-z)s_{j}&=0,\ j\ge2\\
	(1-z)s_{0}+z(z-1)s_{1}&=1,
\end{align*}
and in matrix form
\begin{equation*}
	\begin{pmatrix}
		1-z & z(z-1) & 0&\cdots  &0  \\
		z(z-1) & 1-z-z^3 &z(z-1)& \cdots & 0 \\
		&		z(z-1) & 1-z-z^3 &z(z-1) & 0 \\
		\vdots  & \ddots  & \ddots & \ddots  & \vdots\\
		0 & 0 & \cdots&z(z-1) &1-z-z^3
	\end{pmatrix}
	\begin{pmatrix}
		s_0\\
		s_1\\
		s_2\\
		\vdots\\
		s_K
	\end{pmatrix}=
	\begin{pmatrix}
		1\\
		0\\
		0\\
		\vdots\\
		0
	\end{pmatrix}
\end{equation*}
The next two quantities need no computation, since $f_0^{[K]}=\frac z{1-z}$ and $g_0^{[K]}=1$ for all $K$. Finally
\begin{equation*}
	\begin{pmatrix}
		1-z-z^2 & z(z-1) & 0&\cdots  &0  \\
		z(z-1) & 1-z-z^3 &z(z-1)& \cdots & 0 \\
		&		z(z-1) & 1-z-z^3 &z(z-1) & 0 \\
		\vdots  & \ddots  & \ddots & \ddots  & \vdots\\
		0 & 0 & \cdots&z(z-1) &1-z-z^3
	\end{pmatrix}
	\begin{pmatrix}
		h_0\\
		h_1\\
		h_2\\
		\vdots\\
		h_K
	\end{pmatrix}=
	\begin{pmatrix}
		z^2\\
		0\\
		0\\
		\vdots\\
		0
	\end{pmatrix}
\end{equation*}
Let $\mD_K$ the usual determinant, then $\mD_{K}=(1-z-z^3)\mD_{K-1}-z^2(z-1)^2\mD_{K-2}$ and
\begin{equation*}
\mD_j=\frac{z^{j+1}(1-z)^{j+1}}{W}(r_2^{j+1}-r_1^{j+1}),\quad \mD_0=1,\ \mD_1=1-z-z^3.
\end{equation*}
Furthermore, with $t_K=\mD_{K-1}/\mD_{K}$ and $\tau_K=z^2(z-1)^2t_K$,
\begin{equation*}
t_K=\frac{1}{1-z-z^3-z^2(z-1)^2t_{K-1}} \Longrightarrow \tau_K=\frac{z^2(z-1)^2}{1-z-z^3-\tau_{K-1}},\ \tau_0=0.
\end{equation*}
Finally
\begin{align*}
s_0^{[K+1]}=\sigma_0&=\frac{\mD_{K}}{(1-z)\mD_{K}-z^2(z-1)^2\mD_{K-1}}=
\frac{1}{1-z-\tau_{K}},\ s_0^{[0]}=\frac1{1-z},\\
h_0^{[K+1]}=\eta_0&=\frac{z^2\mD_{K}}{(1-z-z^2)\mD_{K}-z^2(z-1)^2\mD_{K-1}}=
\frac{z^2}{1-z-z^2-\tau_{K}},\ h_0^{[0]}=0.
\end{align*}

%% file: UUDD.tex
\section{UU and DD are forbidden, A004149}

\begin{figure}[h]

	\begin{center}
		\begin{tikzpicture}[scale=1.5,main node/.style={circle,draw,font=\Large\bfseries}]

			\fill (0,-1) circle (0.08cm);

			\foreach \x in {0,1,2,3,4,5,6,7}
			{
				\draw[thick, -latex] (\x,0-1) [in=-165, out=75]to  (\x+1,0);   
				\draw[ thick,red, -latex] (\x+1,-1)[out=-110, in=20]to  (\x,-2);    

				\draw[ thick, latex-] (\x+1,0) [in=75,out=-125] to  (\x,-2);    
				\draw[ thick,red, -latex] (\x+1,0) [in=45,out=-105] to  (\x,-2);    
				
			}                        
			\foreach \x in {0,1,2,3,4,5,6,7,8}
			{
				\draw[thick,latex-,blue ] (\x,-1)  ..  controls (\x-0.3,-1+0.3) and  (\x-0.3,-1-0.3) .. (\x,-1) ;	
				\draw[ thick,blue, -latex] (\x,-2)to  (\x,-1);    
				
			}                        
			
			\foreach \x in {1,2,3,4,5,6,7,8}
			{
				
				\draw[ thick,blue, -latex] (\x,0)to  (\x,-1);

			}                        
			
			
			\foreach \x in {1,2,3,4,5,6,7,8}
			{
				
				\draw (\x,-2) circle (0.05cm);
				\fill (\x,-2) circle (0.05cm);
			}
			
			\foreach \x in {0,1,2,3,4,5,6,7,8}
			{
				
				\draw (\x,-2) circle (0.05cm);
				\fill (\x,-2) circle (0.05cm);
			}

			\foreach \x in {1,2,3,4,5,6,7,8}
			{
				\node at  (\x+0.1,0.1){\tiny$\x$};
				\draw (\x,0) circle (0.05cm);
				\fill (\x,0) circle (0.05cm);
				\draw (\x,-1) circle (0.05cm);
				\fill (\x,-1) circle (0.05cm);
			}
		\end{tikzpicture}
	\end{center}
\caption{Graph (automaton) to recognize  Motzkin paths with forbidden subwords UU and DD.}
\label{purpelUUHH}
\end{figure}
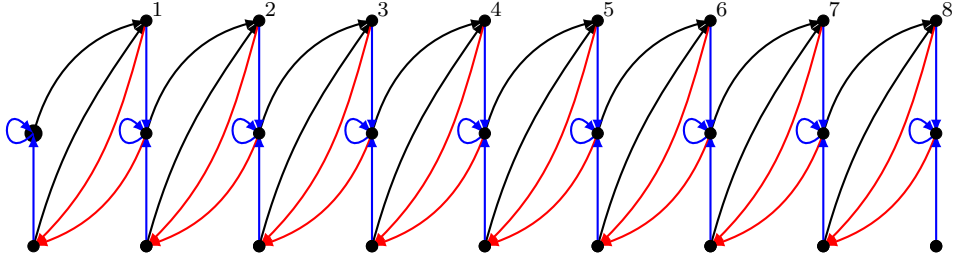
The recursions are 
\begin{align*}
	f_{i+1}&=zg_i+zh_i, \ f_0=0,\\
	g_i&=[i=0]+zf_i+zg_i+zh_i,\ i\ge0,\\
	h_{i}&=zf_{i+1}+zg_{i+1},
\end{align*}
leading to the functional equations
\begin{align*}
	F(u)&=zuG(u)+zuH(u),\\
	G(u)&=1+zF(u)+zG(u)+zH(u),\\
	H(u)&=\frac zuF(u)+\frac zu(G(u)-G(0)).
\end{align*}
Solving the system, a denominator pops up which can be factored:
\begin{equation*}
z^2+z^3u-u+zu+z^2u+z^2u^2=z^2(u-r_1)(u-r_2)
\end{equation*}
with
\begin{equation*}
 r_1=\frac{1-z-z^2-z^3-W}{2z^2},\ r_2=\frac{1-z-z^2-z^3+W}{2z^2},\quad 	W=\sqrt{(1-z^4)(1-2z-z^2)}.
\end{equation*}
Dividing out the factor $u-r_1$ and simplifying, we find
\begin{gather*}
	F(u)=\frac{ur_1}{z(1-ur_1)},\ G(u)=\frac{r_1+z}{z(1-ur_1)},\ H(u)=\frac{(1-z)r_1-z^2}{z^2(1-ur_1)},\ 
	S(u)=\frac{(1+zu)r_1}{z^2(1-ur_1)}.
\end{gather*}
The sequence of excursions 
\begin{equation*}
S(0)=1+z+2z^2+4z^3+8z^4+16z^5+33z^6+69z^7+146z^8+312z^9+673z^{10}+\cdots
\end{equation*}
 is  A308435 in OEIS; the sequence of meanders
 \begin{equation*}
 	S(1)=1+2z+4z^2+9z^3+20z^4+45z^5+102z^6+233z^7+535z^8+1234z^9+2857z^{10}+\cdots
 \end{equation*}
is A308435 in OEIS.
We can also compute the generating function of Motzkin paths ending at level $j\ge1$:
\begin{equation*}
s_j=[u^j]S(u)=[u^j]\frac{(1+zu)r_1}{z^2(1-ur_1)}=\frac{r_1^{j+1}}{z^2}+\frac{r_1^j}{z}.
\end{equation*}
 Now we move to the restricted versions, where indices $>K$ cannot occur.
 
From the denominator we derive the recursion
\begin{align*}
	z^2s_{j-2}+(z^3+z^2+z-1)s_{j-1}+z^2s_j&=0,\ j\ge3,\\
	z^2s_{0}+(z^3+z^2+z-1)s_{1}+z^2s_2&=-z,\\
	(z^2+z-1)s_0+z^2s_1&=-1,
\end{align*}
or in matrix form
{\small
\begin{equation*}
	\begin{pmatrix}
		z^2+z-1& z^2& 0&\cdots  &0  \\
		z^2 & z^3+z^2+z-1 &z^2& \cdots & 0 \\
		&		z^2 & z^3+z^2+z-1 &-z & 0 \\
		\vdots  & \ddots  & \ddots & \ddots  & \vdots\\
		0 & 0 & \cdots&z^2 &z^3+z^2+z-1
	\end{pmatrix}
	\begin{pmatrix}
		s_0\\
		s_1\\
		s_2\\
		\vdots\\
		s_K
	\end{pmatrix}=
	\begin{pmatrix}
		-1\\
		-z\\
		0\\
		\vdots\\
		0
	\end{pmatrix}
\end{equation*}}
The sequence for the first layer does not need computations, since we always have $f_0^{[K]}$, any $K$.
{\footnotesize
	\begin{equation*}
		\begin{pmatrix}
			z^2+z-1& z^2& 0&\cdots  &0  \\
			z^2 & z^3+z^2+z-1 &z^2& \cdots & 0 \\
			&		z^2 & z^3+z^2+z-1 &-z & 0 \\
			\vdots  & \ddots  & \ddots & \ddots  & \vdots\\
			0 & 0 & \cdots&z^2 &z^3+z^2+z-1
		\end{pmatrix}
		\begin{pmatrix}
			g_0\\
			g_1\\
			g_2\\
			\vdots\\
			g_K
		\end{pmatrix}=
		\begin{pmatrix}
			-1+z^2\\
			0\\
			0\\
			\vdots\\
			0
		\end{pmatrix}
\end{equation*}}
and
{\small
	\begin{equation*}
		\begin{pmatrix}
			2z-1& z^2(1-z)& 0&\cdots  &0  \\
			z^2 & z^3+z^2+z-1 &z^2& \cdots & 0 \\
			&		z^2 & z^3+z^2+z-1 &-z & 0 \\
			\vdots  & \ddots  & \ddots & \ddots  & \vdots\\
			0 & 0 & \cdots&z^2 &z^3+z^2+z-1
		\end{pmatrix}
		\begin{pmatrix}
			h_0\\
			h_1\\
			h_2\\
			\vdots\\
			h_K
		\end{pmatrix}=
		\begin{pmatrix}
			-z^2\\
			0\\
			0\\
			\vdots\\
			0
		\end{pmatrix}
\end{equation*}}
 Now we express the solutions in terms of continued fractions. With the usual determinant, 
\begin{equation*}
	\mathcal{D}_{j}=(z^3+z^2+z-1)\mathcal{D}_{j-1} -z^4\mathcal{D}_{j-2}=
	\frac{(-z^2)^{j+1}}{W}(r_1^{j+1}-r_2^{j+1})
\end{equation*}
Further, with $t_K=\mD_{K-1}/\mD_{K}$ and $\tau_K=-z^4t_K$,
\begin{equation*}
t_K=\frac{1}{z^3+z^2+z-1 -z^4t_{K-1}},\quad \tau_K=\frac{z^4}{1-z-z^2-z^3 -\tau_{K-1}}.
\end{equation*}
Now we can express the various solutions in continued fraction form:
\begin{align*}
	s_0=\sigma_{K+1}&=\frac{-\mathcal{D}_{K}+z^3\mathcal{D}_{K-1}}{(z^2+z-1)\mathcal{D}_{K}-z^4\mathcal{D}_{K-1}}\\
	&=\frac{-1+z^3t_{K}}{z^2+z-1-z^4t_{K-1}}=-\frac1z\frac{1-z^2}{z(1-z-z^2-\tau_K)},\ \sigma_0=\frac1{1-z},\\
	g_0=\gamma_{K+1}&=\frac{(-1+z^2)\mathcal{D}_{K}}{(z^2+z-1)\mathcal{D}_{K}-z^4\mathcal{D}_{K-1}}\\&
	=\frac{-1+z^2}{z^2+z-1-z^4t_{K}}=\frac{1-z^2}{1-z-z^2-\tau_{K}},\ \gamma_{0}=\frac1{1-z}\\
	h_0=\eta_{K+1}&=\frac{-z^2\mathcal{D}_{K}}{(2z-1)\mathcal{D}_{K}-z^4(1-z)\mathcal{D}_{K-1}}\\
	&=\frac{-z^2}{(2z-1)-z^4(1-z)t_{K}}=\frac{z^2}{1-2z-(1-z)\tau_{K}},\ \eta_0=0.
	\end{align*}

%% file: UUDU.tex
\section{UU and DU are forbidden, A023431}
\begin{figure}[h]

	\begin{center}
		\begin{tikzpicture}[scale=1.5,main node/.style={circle,draw,font=\Large\bfseries}]

			\fill (0,-1) circle (0.08cm);

			\foreach \x in {0,1,2,3,4,5,6,7}
			{
				\draw[thick, -latex] (\x,0-1) [in=-145, out=55]to  (\x+1,0);   
				\draw[ thick,red, -latex] (\x+1,-1)[out=-125, in=35]to  (\x,-2);    
				\draw[thick,latex-,blue ] (\x,-1)  ..  controls (\x-0.3,-1+0.3) and  (\x-0.3,-1-0.3) .. (\x,-1) ;	
				\draw[ thick, -latex,red] (\x+1,0) [in=75,out=-125] to  (\x,-2);    
				\draw[ thick,red, -latex] (\x+1,-2) to  (\x,-2);    
				
			}  
			
				\foreach \x in {8}
			{
				
				\draw[thick,latex-,blue ] (\x,-1)  ..  controls (\x-0.3,-1+0.3) and  (\x-0.3,-1-0.3) .. (\x,-1) ;

			}

			\foreach \x in {1,2,3,4,5,6,7}
			{				\node at  (\x+0.1,0.1){\tiny$\x$};
			}                
			\foreach \x in {0,1,2,3,4,5,6,7,8}
			{
				
				\draw[ thick,blue, -latex] (\x,-2)to  (\x,-1);    
				
			}                        
			\foreach \x in {1,2,3,4,5,6,7,8}
			{
				
				\draw[ thick,blue, -latex] (\x,0)to  (\x,-1);

			} 
			\node at  (8+0.2,0.1){\tiny$8$};
			
			\foreach \x in {0,1,2,3,4,5,6,7,8}
			{
				\draw (\x,-1) circle (0.05cm);
				\fill (\x,-1) circle (0.05cm);
				\draw (\x,-2) circle (0.05cm);
				\fill (\x,-2) circle (0.05cm);
			}
			
			\foreach \x in {1,2,3,4,5,6,7,8}
			{
				\draw (\x,0) circle (0.05cm);
				\fill (\x,0) circle (0.05cm);
			}
		\end{tikzpicture}
	\end{center}
\caption{Graph (automaton) to recognize  Motzkin paths with forbidden subwords UU and DU.}
\label{purpelUUDU}
\end{figure}
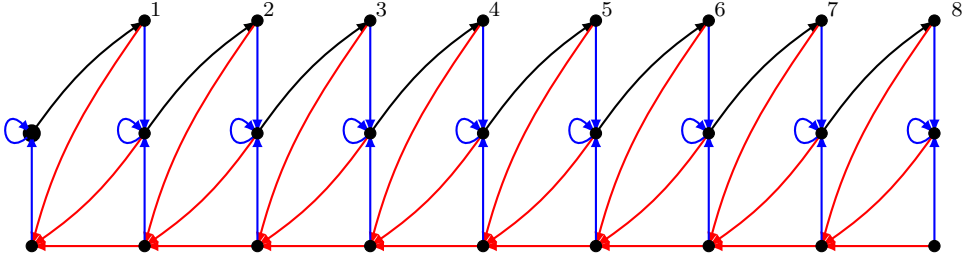
The recursions are
\begin{align*}
	f_{i+1}&=zg_i, \ f_0=0,\\
	g_i&=[i=0]+zf_i+zg_i+zh_i,\ i\ge0,\\
	h_{i}&=zf_{i+1}+zg_{i+1}+zh_{i+1},
\end{align*}
and the according functional equations
\begin{align*}
	F(u)&=zuG(u),\\
	G(u)&=1+zF(u)+zG(u)+zH(u),\\
	H(u)&=\frac zu(F(u)-F(0))+\frac zu(G(u)-G(0))+\frac zu(H(u)-H(0));
\end{align*}
the denominator factors
\begin{equation*}
	z^2u^2+(z-1)u+z=z^2(u-r_1)(u-r_2),
\end{equation*}
with
\begin{equation*}
\ r_1=\frac{1-z-W}{2z^2}\ r_2=\frac{1-z+W}{2z^2},\quad	W=\sqrt{1-2z+z^2-4z^3}.
\end{equation*}
Then, dividing out the factor $u-r_1$ and simplifying, we find
\begin{equation*}
	F=\frac{ur_1}{1-uzr_1},\ G=\frac{r_1}{z(1-uzr_1)},\ H=\frac{r_1(1-z)-z}{z^2(1-uzr_1)},\ S= \frac{(1+uz^2)r_1-z}{z^2(1-uzr_1)}.
\end{equation*}
The series of excursions (A023431 in OEIS) is
\begin{equation*}
S(0)=1+z+2z^2+4z^3+7z^4+13z^5+26z^6+52z^7+104z^8+212z^9+438z^{10}+\cdots
\end{equation*}
and of meanders (not in the OEIs)
\begin{equation*}
	S(1)=1+2z+4z^2+8z^3+16z^4+33z^5+69z^6+145z^7+307z^8+655z^9+1405z^{10}+\cdots.
\end{equation*}

Now we move to Motzkin paths of bounded height, provided by the denominator:
\begin{align*}
	z^2s_{j-2}+(z-1)s_{j-1}+zs_j&=0,\ j\ge3,\\
	z^2s_{0}+(z-1)s_{1}+zs_2&=-z,\\
	(z-1)s_0 +zs_1&=-1,
\end{align*}
or in matrix form
\begin{equation*}
	\begin{pmatrix}
		z-1& z& 0&\cdots  &0  \\
		z^2 & z-1 &z& \cdots & 0 \\
		&		z^2 & z-1 &z & 0 \\
		\vdots  & \ddots  & \ddots & \ddots  & \vdots\\
		0 & 0 & \cdots&z^2 &z-1
	\end{pmatrix}
	\begin{pmatrix}
		s_0\\
		s_1\\
		s_2\\
		\vdots\\
		s_K
	\end{pmatrix}=
	\begin{pmatrix}
		-1\\
		-z\\
		0\\
		\vdots\\
		0
	\end{pmatrix}
\end{equation*}
Since $f_0^{[K]}=0$ for any $K$, we move to the next one:
\begin{align*}
	z^2s_{j-2}+(z-1)s_{j-1}+zs_j&=0,\ j\ge3\\
	z^2s_{0}+(z-1)s_{1}+zs_2&=-z,\\
	(z-1)s_0 +zs_1&=-1,
\end{align*}
or in matrix form
\begin{equation*}
	\begin{pmatrix}
		z-1& z& 0&\cdots  &0  \\
		z^2 & z-1 &z& \cdots & 0 \\
		&		z^2 & z-1 &z & 0 \\
		\vdots  & \ddots  & \ddots & \ddots  & \vdots\\
		0 & 0 & \cdots&z^2 &z-1
	\end{pmatrix}
	\begin{pmatrix}
		g_0\\
		g_1\\
		g_2\\
		\vdots\\
		g_K
	\end{pmatrix}=
	\begin{pmatrix}
		-1\\
		0\\
		0\\
		\vdots\\
		0
	\end{pmatrix}
\end{equation*}
{\small	
	\begin{equation*}
	\begin{pmatrix}
		1-2z+z^2-z^3& z(z-1)& 0&\cdots  &0  \\
		z^2 & z-1 &z& \cdots & 0 \\
		&		z^2 & z-1 &z & 0 \\
		\vdots  & \ddots  & \ddots & \ddots  & \vdots\\
		0 & 0 & \cdots&z^2 &z-1
	\end{pmatrix}
	\begin{pmatrix}
		h_0\\
		h_1\\
		h_2\\
		\vdots\\
		h_K
	\end{pmatrix}=
	\begin{pmatrix}
		z^2\\
		0\\
		0\\
		\vdots\\
		0
	\end{pmatrix}
\end{equation*}
}
With the usual determinant, we derive the recursion
\begin{equation*}
	\mD_{K}=(z-1)\mD_{K-1}-z^3\mD_{K-2}=\frac{(-z^2)^{K+1}}{W}(r_1^{K+1}-r_2^{K+1})
\end{equation*}
and also, with $t_K=\mD_{K-1}/\mD_{K}$ and $\tau_K=-z^3t_K$
\begin{equation*}
	t_K=\frac{1}{-1+z-z^3t_{K-1}}, \quad \tau_K=\frac{z^3}{1-z-\tau_{K-1}}
\end{equation*}
Now we can express the solutions:
\begin{align*}
	s_0=\sigma_{K+1}&=\frac{-\mD_{K}+z^2\mD_{K-1}}{(z-1)\mD_{K}-z^3\mD_{K-1}}=
	\frac{-1+z^2t_K}{z-1-z^3t_K}\\*&=-\frac1z+\frac1{z(1-z-\tau_K)},\quad \sigma_0=\frac1{1-z}\\
	g_0=\gamma_{K+1}&=\frac{-\mD_{K}}{(z-1)\mD_{K}-z^3\mD_{K-1}}=\frac{1}{1-z-\tau_{K-1}},\quad \gamma_0=\frac1{1-z},\\
	h_0=\eta_{K+1}&=\frac{z^2\mD_{K}}{(1-2z+z^2-z^3)\mD_{K}-z^3(z-1)\mD_{K-1}}\\&=\frac{z^2}{1-2z+z^2-z^3-z^3(z-1)t_{K}}
=\frac{z^2}{1-2z+z^2-z^3-(1-z)\tau_{K}},\quad \eta_0=0.
\end{align*}

%% file: UUHD.tex
\section{UU and HD are forbidden, A217282}
\begin{figure}[h]

	\begin{center}
		\begin{tikzpicture}[scale=1.5,main node/.style={circle,draw,font=\Large\bfseries}]

			\fill (0,-1) circle (0.08cm);

			\foreach \x in {0,1,2,3,4,5,6,7}
			{
				\draw[thick, -latex] (\x,-1) [in=-160, out=75]to  (\x+1,0);   
				\draw[ thick,red, -latex] (\x+1,-0)[out=-135, in=45]to  (\x,-1);    
								\draw[ thick,red, -latex] (\x+1,-1)to  (\x,-1);    
								\draw[ thick,-latex] (\x,-2)to  (\x+1,0);    


						}                        
								\foreach \x in {0,1,2,3,4,5,6,7,8}
								{				\draw[thick,latex-,blue ] (\x,-2)  ..  controls (\x-0.3,-2-0.3) and  (\x+0.3,-2-0.3) .. (\x,-2) ;
										}

			\foreach \x in {1,2,3,4,5,6,7,8}
			{
	\draw[ thick,blue, -latex] (\x+0.025,0)[in=110,out=-110]to  (\x-0.025,-2);        
	\draw[ thick,blue, -latex] (\x,-1)to  (\x,-2);        
				\node at  (\x+0.1,0.1){\tiny$\x$};
			}
			\foreach \x in {0,1,2,3,4,5,6,7,8}
			{

				\draw[ thick,blue, -latex] (\x,-1)to  (\x,-2);        
				\node at  (\x+0.1,0.1){\tiny$\x$};
			}
			
			
			\foreach \x in {0,1,2,3,4,5,6,7,8}
			{
				\draw (\x,-1) circle (0.05cm);
				\fill (\x,-1) circle (0.05cm);
				\draw (\x,-2) circle (0.05cm);
				\fill (\x,-2) circle (0.05cm);
			}
			\foreach \x in {1,2,3,4,5,6,7,8}{				\draw (\x,0) circle (0.05cm);
				\fill (\x,0) circle (0.05cm);
			}
		\end{tikzpicture}
	\end{center}
\caption{Graph (automaton) to recognize  Motzkin paths with forbidden subwords UU and HD.}
\label{purpelUUHD}
\end{figure}
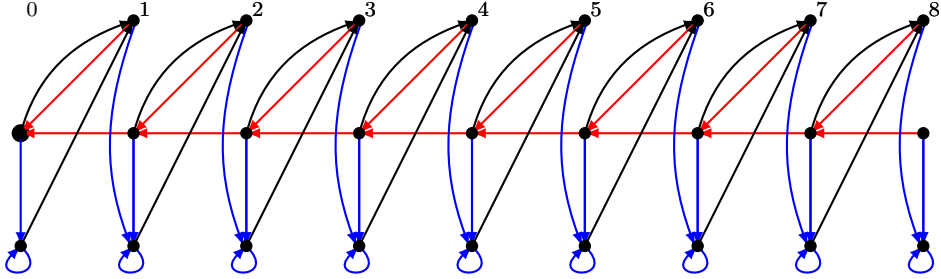

The recursions that can be read off from Figure~\ref{purpelUUHD} are
\begin{align*}
	f_{i+1}&=zg_i+zh_i, \ f_0=0,\\
	g_i&=[i=0]+zf_{i+1}+zg_{i+1},\ i\ge0,\\
	h_{i}&=zf_{i}+zg_{i}+zh_{i},
\end{align*}
and the corresponding functional equations
\begin{align*}
	F(u)&=zuG(u)+zuH(u),\\
	G(u)&=1+\frac zuF(u)+\frac zu (G(u)-G(0)),\\
	H(u)&=zF(u)+zG(u)+zH(u).
\end{align*}
The factorization of the ubiquitous denominator:
\begin{equation*}
z^2u^2  -(z+1)(-1+z)^2u -z(-1+z)=z^2(u-r_1)(u-r_2),
\end{equation*}
with 
\begin{equation*}
r_1=\frac{1-z-z^2+z^3-W}{2z^2},\ r_2=\frac{1-z-z^2+z^3+W}{2z^2},  \quad W=\sqrt{1-2z-{z}^{2}+3{z}^{4}-2{z}^{5}+{z}^{6}}
\end{equation*}
After dividing out $u-r_1$ and simplifying, we find
\begin{gather*}
	F(u)=\frac{ur_1}{ 1-z-uzr_1},\
	G(u)=\frac{(1-z-z^2u)r_1}{z(1-z-uzr_1)},\
	H(u)=\frac{(1+zu)r_1}{1-z-uzr_1},\ 
	S(u) =\frac{(1+zu)r_1}{z(1-z-uzr_1)}.
\end{gather*}
The sequence of excursions
\begin{equation*}
S(0)=1+z+2z^2+3z^3+5z^4+9z^5+16z^6+30z^7+57z^8+110z^9+216z^{10}+\cdots,
\end{equation*}
is A217282 in OEIS, and
\begin{equation*}
	S(1)=1+2z+4z^2+8z^3+16z^4+33z^5+68z^6+142z^7+298z^8+629z^9+1334z^{10}+\cdots,
\end{equation*}
the series of Motzkin meanders is not in OEIS. The enumeration of Motzkin paths ending at level $j\ge1$ is
\begin{align*}
[u^j]S(u)&=[u^j]\frac{(1+zu)r_1}{z(1-z)\big(1-\frac{uzr_1}{1-z}\big)}=\frac1{z(1-z)}\Big(\frac{zr_1}{1-z}\Big)^j
+\frac1{1-z}\Big(\frac{zr_1}{1-z}\Big)^{j-1}r_1\\
&=\frac1{z^2}\Big(\frac{zr_1}{1-z}\Big)^{j+1}
+\frac1{z}\Big(\frac{zr_1}{1-z}\Big)^{j}.
\end{align*}

Now we move to the restricted paths.
\begin{align*}
z^2s_{j-2}-(-1+z)^2(z+1)s_{j-1}+z(1-z)s_j=0,\ j\ge3\\
z^2s_{0}-(-1+z)^2(z+1)s_{1}+z(1-z)s_0=-z,\\
(z-1)s_0+z(1-z)=-1
\end{align*}
or in matrix form
{\footnotesize
	\begin{equation*}
	\begin{pmatrix}
		z-1& z(1-z)& 0&\cdots  &0  \\
		z^2 & -(-1+z)^2(z+1) &z(1-z)& \cdots & 0 \\
		&		z^2 & -(-1+z)^2(z+1) &z(1-z) & 0 \\
		\vdots  & \ddots  & \ddots & \ddots  & \vdots\\
		0 & 0 & \cdots&z^2 &-(-1+z)^2(z+1)
	\end{pmatrix}
	\begin{pmatrix}
		s_0\\
		s_1\\
		s_2\\
		\vdots\\
		s_K
	\end{pmatrix}=
	\begin{pmatrix}
		-1\\
		-z\\
		0\\
		\vdots\\
		0
	\end{pmatrix}
\end{equation*}
}
{\scriptsize
	\begin{equation*}
		\begin{pmatrix}
			-1+z+z^2& z(1-z)& 0&\cdots  &0  \\
			z^2 & -(-1+z)^2(z+1) &z(1-z)& \cdots & 0 \\
			&		z^2 & -(-1+z)^2(z+1) &z(1-z) & 0 \\
			\vdots  & \ddots  & \ddots & \ddots  & \vdots\\
			0 & 0 & \cdots&z^2 &-(-1+z)^2(z+1)
		\end{pmatrix}
		\begin{pmatrix}
			g_0\\
			g_1\\
			g_2\\
			\vdots\\
			g_K
		\end{pmatrix}=
		\begin{pmatrix}
			-1+z\\
			z^2\\
			0\\
			\vdots\\
			0
		\end{pmatrix}
	\end{equation*}
}
{\footnotesize
	\begin{equation*}
		\begin{pmatrix}
			-1+z& z(1-z)& 0&\cdots  &0  \\
			z^2 & -(-1+z)^2(z+1) &z(1-z)& \cdots & 0 \\
			&		z^2 & -(-1+z)^2(z+1) &z(1-z) & 0 \\
			\vdots  & \ddots  & \ddots & \ddots  & \vdots\\
			0 & 0 & \cdots&z^2 &-(-1+z)^2(z+1)
		\end{pmatrix}
		\begin{pmatrix}
			h_0\\
			h_1\\
			h_2\\
			\vdots\\
			h_K
		\end{pmatrix}=
		\begin{pmatrix}
			-z\\
			-z^2\\
			0\\
			\vdots\\
			0
		\end{pmatrix}
	\end{equation*}
}

For the usual determinant we derive a recursion and an explicit Binet form:
\begin{equation*}
	\mathcal{D}_{K}=-(1-z)^2(z+1)\mathcal{D}_{K-1}-z^3(1-z)\mathcal{D}_{K-2}=\frac{(-z^2)^{K+1}}{W}(r_1^{K+1}-r_2^{K+1}).
\end{equation*}
Alternatively, with $t_K=\mD_{K-1}/\mD_K$ and $\tau_K=-z^3t_K$
\begin{equation*}
t_K=\frac{1}{-z^3+z+z^2-1-z^3(1-z)t_{K-1}} \Longrightarrow \tau_K=\frac{z^3}{1-z-z^2+z^3-(1-z)\tau_{K-1}},\ \tau_0=0.
\end{equation*}
Finally, $\phi_0=0$ for all $K$,   and
\begin{align*}
	s_0=\sigma_{K+1}&=\frac{-\mD_K+z^2(1-z)\mD_{K-1}}{(z-1)\mD_K-z^3(1-z)\mD_{K-1}}
	=\frac{-1+z^2(1-z)t_{K}}{(z-1) -z^3(1-z)t_{K}}\\
		&=-\frac1z+\frac{1}{z(1-z)(1-\tau_K)},\quad \sigma_0=\frac1{1-z},\\
	g_0=\gamma_{K+1}&=\frac{(-1+z)\mD_K-z^3(1-z)\mD_{K-1}}{(-1+z+z^2)\mD_K-z^3(1-z)\mD_{K-1}}
	=\frac{(-1+z)-z^3(1-z)t_K}{(-1+z+z^2)-z^3(1-z)t_K}\\
	&=1+\frac{z^2}{1-z-z^2-(1-z)\tau_K},\quad\gamma_0=1,\\
	h_0=\eta_{K+1}&=\frac{-z\mD_K+z^3(1-z)\mD_{K-1}}{(z-1)\mD_K-z^3(1-z)\mD_{K-1}}
	=\frac{-z+z^3(1-z)t_{K}}{(z-1)-z^3(1-z)t_{K}}\\
	&=-1+\frac{1}{(1-z)(1-\tau_K)},\quad \eta_0=\frac{z}{1-z}.
\end{align*}

%% file: UUUD.tex
\section{UU and UD are forbidden, A023431}
\begin{figure}[h]

	\begin{center}
		\begin{tikzpicture}[scale=1.5,main node/.style={circle,draw,font=\Large\bfseries}]

			\fill (0,-1) circle (0.08cm);

			\foreach \x in {0,1,2,3,4,5,6,7}
			{
				\draw[thick, -latex] (\x,0-1) [in=-145, out=55]to  (\x+1,0);   
				\draw[ thick,red, -latex] (\x+1,-1)[out=-125, in=35]to  (\x,-2);    

				\draw[ thick, latex-] (\x+1,0) [in=75,out=-125] to  (\x,-2);    
				\draw[ thick,red, -latex] (\x+1,-2) to  (\x,-2);    

						}                        
								\foreach \x in {0,1,2,3,4,5,6,7,8}
								{				\draw[thick,latex-,blue ] (\x,-1)  ..  controls (\x-0.3,-1+0.3) and  (\x-0.3,-1-0.3) .. (\x,-1) ;
												\draw[ thick,blue, -latex] (\x,-2)to  (\x,-1);    }

			\foreach \x in {1,2,3,4,5,6,7,8}
			{
				\draw[ thick,blue, -latex] (\x,0)to  (\x,-1);        
				\node at  (\x+0.1,0.1){\tiny$\x$};
			}

			
			\foreach \x in {0,1,2,3,4,5,6,7,8}
			{
				\draw (\x,-1) circle (0.05cm);
				\fill (\x,-1) circle (0.05cm);
				\draw (\x,-2) circle (0.05cm);
				\fill (\x,-2) circle (0.05cm);
			}
			\foreach \x in {1,2,3,4,5,6,7,8}{				\draw (\x,0) circle (0.05cm);
				\fill (\x,0) circle (0.05cm);
			}
		\end{tikzpicture}
	\end{center}
\caption{Graph (automaton) to recognize  Motzkin paths with forbidden subwords UU and UD.}
\label{purpelUUUD}
\end{figure}
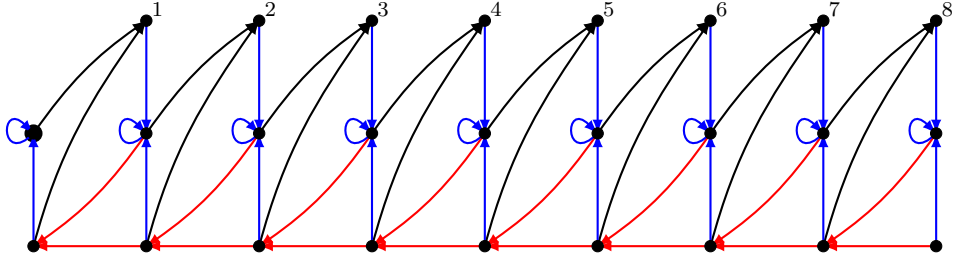
The recursions that can be read off from Figure~\ref{purpelUUUD} are
\begin{align*}
	f_{i+1}&=zg_i+zh_i, \ f_0=0,\\
	g_i&=[i=0]+zf_i+zg_i+zh_i,\ i\ge0,\\
	h_{i}&=zg_{i+1}+zh_{i+1},
\end{align*}
and the corresponding functional equations
\begin{align*}
	F(u)&=zuG(u)+zuH(u),\\
	G(u)&=1+zF(u)+zG(u)+zH(u),\\
	H(u)&=\frac zu(G(u)-G(0))+\frac zu(H(u)-H(0)).
\end{align*}
Solving this, there is a denominator $z+(-1+z)u+z^2u^2=z^2(u-r_1)(u-r_2)$
with
\begin{equation*}
	r_1=\frac{1-z-W}{2z^2},\ r_2=\frac{1-z+W}{2z^2},\quad W=\sqrt{1-2z+z^2-4z^3}.
\end{equation*}
Dividing the factor $u-r_1$ out of the solutions and simplifying, we get
\begin{gather*}
	F(u)=\frac{ur_1}{1-uzr_1},\ G(u)=\frac{1+r_1}{1-uzr_1},\ H(u)=\frac{(1-z)r_1-z}{z(1-uzr_1)},\ S(u)=\frac{(1+zu)r_1}{z(1-uzr_1)}.
\end{gather*}
The sequence of excursions 
\begin{equation*}
S(0)=1+z+z^2+2z^3+4z^4+7z^5+13z^6+26z^7+52z^8+104z^9+212z^{10}+\cdots
\end{equation*}
is A023431 in OEIS, the sequence of meanders
\begin{equation*}
	S(1)=1+2z+3z^2+6z^3+12z^4+24z^5+49z^6+102z^7+214z^8+452z^9+962z^{10}+\cdots
\end{equation*}
is not in OEIS.

From the denominator we derive the recursion
\begin{align*}
	z^2s_{j-2}+(z-1)s_{j-1}+zs_j=0,\ j\ge3\\
	z^2s_{0}+(z-1)s_{1}+zs_2&=-z,\\
	(-1+z-z^2)s_0+zs_1&=-1,
\end{align*}
or in matrix form
\begin{equation*}
	\begin{pmatrix}
		-1+z-z^2& z& 0&\cdots  &0  \\
		z^2 & z-1 &z& \cdots & 0 \\
		&		z^2 & z-1 &z & 0 \\
		\vdots  & \ddots  & \ddots & \ddots  & \vdots\\
		0 & 0 & \cdots&z^2 &z-1
	\end{pmatrix}
	\begin{pmatrix}
		s_0\\
		s_1\\
		s_2\\
		\vdots\\
		s_K
	\end{pmatrix}=
	\begin{pmatrix}
		-1\\
		-z\\
		0\\
		\vdots\\
		0
	\end{pmatrix}
\end{equation*}
$\phi_K=0$ for any $K$, so we further get
\begin{equation*}
	\begin{pmatrix}
		-1+z-z^2& z& 0&\cdots  &0  \\
		z^2 & z-1 &z& \cdots & 0 \\
		&		z^2 & z-1 &z & 0 \\
		\vdots  & \ddots  & \ddots & \ddots  & \vdots\\
		0 & 0 & \cdots&z^2 &z-1
	\end{pmatrix}
	\begin{pmatrix}
		g_0\\
		g_1\\
		g_2\\
		\vdots\\
		g_K
	\end{pmatrix}=
	\begin{pmatrix}
		-1-z^2\\
		0\\
		0\\
		\vdots\\
		0
	\end{pmatrix}
\end{equation*}
\begin{equation*}
	\begin{pmatrix}
		-1+2z-z^2+z^3& z(1-z)& 0&\cdots  &0  \\
		z^2 & z-1 &z& \cdots & 0 \\
		&		z^2 & z-1 &z & 0 \\
		\vdots  & \ddots  & \ddots & \ddots  & \vdots\\
		0 & 0 & \cdots&z^2 &z-1
	\end{pmatrix}
	\begin{pmatrix}
		h_0\\
		h_1\\
		h_2\\
		\vdots\\
		h_K
	\end{pmatrix}=
	\begin{pmatrix}
		-z^3\\
		0\\
		0\\
		\vdots\\
		0
	\end{pmatrix}
\end{equation*}
The traditional determinant satisfies
\begin{equation*}
	\mathcal{D}_{K}=(z-1)\mathcal{D}_{K-1}-z^3\mathcal{D}_{K-2}=\frac{(-z^2)^{K+1}}{W}(r_2^{K+1}-r_1^{K+1}).
\end{equation*}
Further, with $t_K=\mD_{K-1}/\mD_{K}$ and $\tau_K=-z^3t_K$
\begin{equation*}
	t_{K}=\frac1{(z-1)-z^3t_{K-1}}, \Longrightarrow  \tau_K=\frac{z^3}{1-z-\tau_{K-1}},\quad \tau_0=0.
\end{equation*}
Now we can express the solutions in continued fraction form:
\begin{align*}
	s_0=\sigma_{K+1}&=\frac{-\mathcal{D}_{K}+z^2\mathcal{D}_{K-1}}{(-1+z-z^2)\mathcal{D}_{K}-z^3\mathcal{D}_{K-1}}
	=\frac{-1+z^2t_{K}}{(-1+z-z^2)-z^3t_{K}}\\
	&=-\frac1z+\frac{1+z^2}{z(1-z+z^2-\tau_K)},\quad \sigma_0=\frac1{1-z},\\
	g_0=\gamma_{K+1}&=\frac{-(1+z^2)\mathcal{D}_{K}}{(-1+z-z^2)\mathcal{D}_{K}-z^3\mathcal{D}_{K-1}}
	=\frac{1+z^2}{1-z+z^2-\tau_{K}},\quad \gamma_0=\frac1{1-z},\\
	h_0=\eta_{K+1}&=\frac{-z^3\mathcal{D}_{K}}{(-1+2z-z^2+z^3)\mathcal{D}_{K}-z^3(1-z)\mathcal{D}_{K-1}}
	\\*&=\frac{z^3}{1-2z+z^2-z^3-(1-z)\tau_{K}},\quad \eta_0=0.
\end{align*}

%% file: UDDD.tex
\section{UD and DD are forbidden, A023431}
\begin{figure}[h]

	\begin{center}
		\begin{tikzpicture}[scale=1.5,main node/.style={circle,draw,font=\Large\bfseries}]

			\fill (0,-1) circle (0.08cm);

			\foreach \x in {0,1,2,3,4,5,6,7}
			{
				\draw[thick, -latex] (\x,0-1) [in=-145, out=55]to  (\x+1,0);   
				\draw[ thick,red, -latex] (\x+1,-1) to  (\x,-2);    
				\draw[thick,latex-,blue ] (\x,-1)  ..  controls (\x-0.3,-1+0.3) and  (\x-0.3,-1-0.3) .. (\x,-1) ;	
				\draw[ thick, latex-] (\x+1,0)   to  (\x,-2);    
				\node at  (\x+0.1,0.1){\tiny$\x$};
			}                        
			\foreach \x in {1,2,3,4,5,6,7}
{
	\draw[ thick, latex-] (\x+1,0) to  (\x,0);    
}

			\foreach \x in {0,1,2,3,4,5,6,7,8}
			{

				\draw[ thick,blue, -latex] (\x,-2)to  (\x,-1);    
				
			}                        

			\foreach \x in {1,2,3,4,5,6,7,8}
{
	
	\draw[ thick,blue, -latex] (\x,0)to  (\x,-1);

}                        
			
			\node at  (8+0.1,0.1){\tiny$8$};
			
			\foreach \x in {0,1,2,3,4,5,6,7,8}
			{
				\draw (\x,-1) circle (0.05cm);
				\fill (\x,-1) circle (0.05cm);
				\draw (\x,-2) circle (0.05cm);
				\fill (\x,-2) circle (0.05cm);
			}
		\foreach \x in {1,2,3,4,5,6,7,8}
		{
			\draw (\x,0) circle (0.05cm);
			\fill (\x,0) circle (0.05cm);
		}
		\end{tikzpicture}
	\end{center}
	\caption{Graph (automaton) to recognize Motzkin paths with forbidden subwords UU and DD.}
	\label{purpelUDDD}
\end{figure}
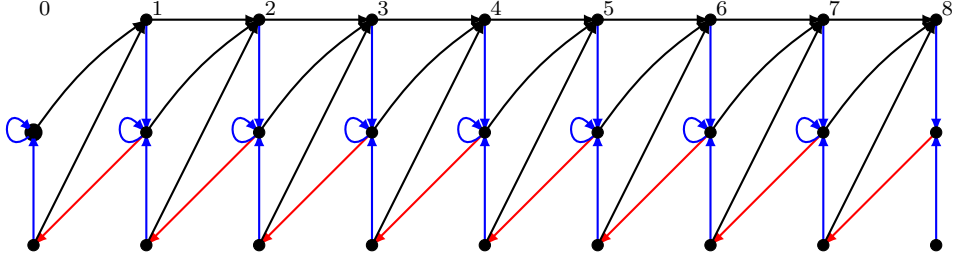
The recursion related to Figure \ref{purpelUDDD} are
\begin{align*}
	f_{i+1}&=zf_i+zg_i+zh_i, \ f_0=0,\\
	g_i&=[i=0]+zf_i+zg_i+zh_i,\ i\ge0\\
	h_{i}&=zg_{i+1}
\end{align*}
and the system of functional equations is
\begin{align*}
	F(u)&=zuF(u)+zuG(u)+zuH(u),\\
	G(u)&=1+zF(u)+zG(u)+zH(u),\\
	H(u)&=\frac zu(G(u)-G(0)).
\end{align*}
The ubiquitous denominator is $zu^2+(z-1)u+z^2=z(u-r_1)(u-r_2)$ with
\begin{equation*}
r_1=\frac{1-z-W}{2z},\ r_2=\frac{1-z+W}{2z},\quad W=\sqrt{1-2z+z^2-4z^3}.
\end{equation*}
Dividing out the factor $u-r_1$ and simplifying, we end up with
\begin{equation*}
	F(u)=\frac{ur_1}{z-ur_1},\ G(u)=1+\frac{r_1}{z-ur_1},\ H(u)=\frac{(1-z)r_1-z^2}{z(z-ur_1)},\ S(u)=\frac{r_1}{z(z-ur_1)}.
\end{equation*}
We find
\begin{equation*}
S(0)=1+z+z^2+2z^3+4z^4+7z^5+13z^6+26z^7+52z^8+104z^9+212z^{10}+\cdots
\end{equation*}
which is A023431 in OEIS
and
\begin{equation*}
	S(1)=1+2z+4z^2+9z^3+21z^4+49z^5+115z^6+272z^7+646z^8+1538z^9+3670z^{10}+\cdots
\end{equation*}
which is not in OEIS. Furthermore
\begin{equation*}
[u^j]S(u)=[u^j]\frac{r_1}{z(z-ur_1)}=[u^j]\frac{r_1}{z^2(1-\frac{ur_1}{z})}=\frac1z\Big(\frac{r_1}{z}\Big)^{j+1}
\end{equation*}

Since, before diving out the factor, we have
\begin{equation*}
	S(u) = \frac{-u-z-uz^2G(0)+zG(0)}{-u+zu^2+z^2+zu},
\end{equation*}
we get a recursion from it:
\begin{align*}
	zs_{j-2}+(z-1)s_{j-1}+z^2s_j&=0,\ j\ge2\\
	(z-1)s_0+z^2s_1&=-1,
\end{align*}
which is in matrix form
\begin{equation*}
	\begin{pmatrix}
		z-1& z^2& 0&\cdots  &0  \\
		z & z-1 &z^2& \cdots & 0 \\
		&		z & z-1 &z^2 & 0 \\
		\vdots  & \ddots  & \ddots & \ddots  & \vdots\\
		0 & 0 & \cdots&z &z-1
	\end{pmatrix}
	\begin{pmatrix}
		s_0\\
		s_1\\
		s_2\\
		\vdots\\
		s_K
	\end{pmatrix}=
	\begin{pmatrix}
		-1\\
		0\\
		0\\
		\vdots\\
		0
	\end{pmatrix}
\end{equation*}
The usual determinant is
\begin{equation*}
	\mD_{K}=(z-1)\mD_{K-1}-z^3\mD_{K-2}=\frac{(-z)^{K+1}}{W}(r_1^{K+1}-r_2^{K+1}).
\end{equation*}
With $t_K=\mD_{K-1}/\mD_{K}$ and $\tau_K=-z^3t_K$
\begin{equation*}
	t_K=\frac{1}{z-1-z^3t_{K-1}}\ \Longrightarrow \ \tau_K=\frac{z^3}{1-z-\tau_{K-1}},\ \tau_0=0
\end{equation*}
\begin{equation*}
	\begin{pmatrix}
		z-1& z^2& 0&\cdots  &0  \\
		z & z-1 &z^2& \cdots & 0 \\
		&		z & z-1 &z^2 & 0 \\
		\vdots  & \ddots  & \ddots & \ddots  & \vdots\\
		0 & 0 & \cdots&z &z-1
	\end{pmatrix}
	\begin{pmatrix}
		g_0\\
		g_1\\
		g_2\\
		\vdots\\
		g_K
	\end{pmatrix}=
	\begin{pmatrix}
		-1\\
		z\\
		0\\
		\vdots\\
		0
	\end{pmatrix}
\end{equation*}
\begin{equation*}
	\begin{pmatrix}
		-1+2z-z^2+z^3& z^2(1-z)& 0&\cdots  &0  \\
		z & z-1 &z^2& \cdots & 0 \\
		&		z & z-1 &z^2 & 0 \\
		\vdots  & \ddots  & \ddots & \ddots  & \vdots\\
		0 & 0 & \cdots&z &z-1
	\end{pmatrix}
	\begin{pmatrix}
		h_0\\
		h_1\\
		h_2\\
		\vdots\\
		h_K
	\end{pmatrix}=
	\begin{pmatrix}
		-z^3\\
		0\\
		0\\
		\vdots\\
		0
	\end{pmatrix}
\end{equation*}
From this, we can express the solutions (return to the $x$-axis) in continued fraction form
\begin{align*}
	s_0=\sigma_{K+1}&=\frac{-\mD_{K}}{(z-1)\mD_{K}-z^3\mD_{K-1}}=\frac{1}{1-z+z^3t_{H}}=\frac{1}{1-z-\tau_{H}},\ \sigma_0=\frac1{1-z},\\
	g_0=\gamma_{K+1}&=\frac{-\mD_{K}-z^3\mD_{K-1}}{(z-1)\mD_{K}-z^3\mD_{K-1}}
	=\frac{-1-z^3t_{K}}{z-1-z^3t_{K}}=1+\frac{z}{1-z-\tau_K},\  \gamma_0=\frac1{1-z},\\
	h_0=\eta_{K+1}&=\frac{-z^3\mD_{K}}{(-1+2z-z^2+z^3)\mD_{K}-z^3(1-z)\mD_{K-1}}=
	\frac{z^3}{1-2z+z^2-z^3-(1-z)\tau_{K}},\ \eta_0=0.
\end{align*}

%% file: UDHU.tex
\section{UD and HU are forbidden, not in OEIS}

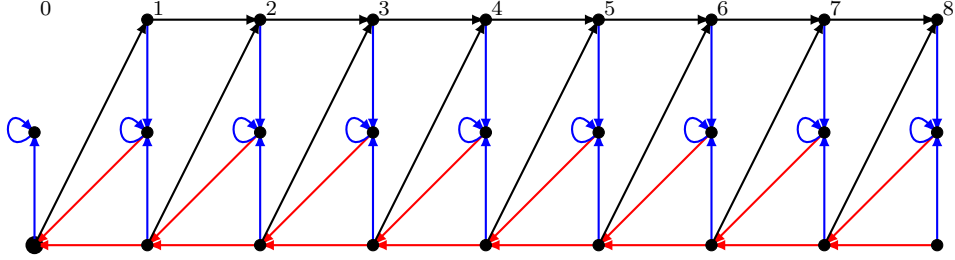
\begin{figure}[h]

	\begin{center}
		\begin{tikzpicture}[scale=1.5,main node/.style={circle,draw,font=\Large\bfseries}]

			\fill (0,-2) circle (0.08cm);

			\foreach \x in {0,1,2,3,4,5,6,7}
			{
				\draw[ thick,red, -latex] (\x+1,-1) to  (\x,-2);    
				\draw[thick,latex-,blue ] (\x,-1)  ..  controls (\x-0.3,-1+0.3) and  (\x-0.3,-1-0.3) .. (\x,-1) ;	
				\draw[ thick, latex-] (\x+1,0)   to  (\x,-2);    
					\draw[ thick,red, latex-] (\x,-2)to  (\x+1,-2);         
				\node at  (\x+0.1,0.1){\tiny$\x$};
			}

		\foreach \x in {8}
		{
			\draw[thick,latex-,blue ] (\x,-1)  ..  controls (\x-0.3,-1+0.3) and  (\x-0.3,-1-0.3) .. (\x,-1) ;	
		}

			\foreach \x in {1,2,3,4,5,6,7}
{
	\draw[ thick, latex-] (\x+1,0) to  (\x,0);    
}

			\foreach \x in {0,1,2,3,4,5,6,7,8}
			{

				\draw[ thick,blue, -latex] (\x,-2)to  (\x,-1);    
				
			}                        

			\foreach \x in {1,2,3,4,5,6,7,8}
{
	
	\draw[ thick,blue, -latex] (\x,0)to  (\x,-1);

}                        
			
			\node at  (8+0.1,0.1){\tiny$8$};
			
			\foreach \x in {0,1,2,3,4,5,6,7,8}
			{
				\draw (\x,-1) circle (0.05cm);
				\fill (\x,-1) circle (0.05cm);
				\draw (\x,-2) circle (0.05cm);
				\fill (\x,-2) circle (0.05cm);
			}
		\foreach \x in {1,2,3,4,5,6,7,8}
		{
			\draw (\x,0) circle (0.05cm);
			\fill (\x,0) circle (0.05cm);
		}
		\end{tikzpicture}
	\end{center}
\caption{Graph (automaton) to recognize  Motzkin paths with forbidden subwords UD and HU.}
	\label{purpelUDHU}
\end{figure}
The recursion related to Figure \ref{purpelUDHU} are
\begin{align*}
	f_{i+1}&=zf_i+zh_i, \ f_0=0,\\
	g_i&=zf_i+zg_i+zh_i,\ i\ge0\\
	h_{i}&=[i=0]+zg_{i+1}+zh_{i+1}
\end{align*}
and the system of functional equations is
\begin{align*}
	F(u)&=zuF(u)+zuH(u),\\
	G(u)&=zF(u)+zG(u)+zH(u),\\
	H(u)&=1+\frac zu(G(u)-G(0))+\frac zu(H(u)-H(0)).
\end{align*}
The ubiquitous denominator is $z(z-1)u^2+(1-z)(1+z^2)u-z=z(z-1)(u-r_1)(u-r_2)$ with
\begin{gather*}
r_1=\frac{(1-z)(1+z^2)-W}{2z(1-z)},\ r_2=\frac{(1-z)(1+z^2)+W}{2z(1-z)},\\ W=\sqrt{(1-z)(1-z-2z^2-2z^3+z^4 -z^5)}.
\end{gather*}
Dividing out the factor $u-r_1$ and simplifying, we end up with
\begin{gather*}
	F(u)=\frac{1}{1-u(1-z)r_1}\,\ G(u)=\frac{r_1}{1-u(1-z)r_1},\ H(u)=\frac{(1-z)r_1-z}{z(1-u(1-z)r_1)},\\ S(u)=\frac{r_1}{z(1-u(1-z)r_1)}.
\end{gather*}
We find
\begin{equation*}
S(0)=1+z+z^2+2z^3+3z^4+5z^5+9z^6+16z^7+30z^8+57z^9+110z^{10}+\cdots
\end{equation*}
which is not in OEIS
and
\begin{equation*}
	S(1)=1+2z+3z^2+5z^3+9z^4+17z^5+33z^6+65z^7+130z^8+263z^9+537z^{10}+\cdots
\end{equation*}
which is also not in OEIS. Furthermore
\begin{equation*}
[u^j]S(u)=[u^j]\frac{r_1}{z(1-u(1-z)r_1)}=\frac{r_1}{z}\big((1-z)r_1\big)^j=\frac{1}{z(1-z)}\big((1-z)r_1\big)^{j+1}.
\end{equation*}
Since, before diving out the factor, we have
\begin{equation*}
	S(u) =  \frac{u+z+uz^2G(0)-zG(0)}{z(z-1)u^2+(1-z)(1+z^2)u-z},
\end{equation*}
we get a recursion from it:
\begin{align*}
	z(z-1)s_{j-2}+(1-z)(1+z^2)s_{j-1}-zs_j&=0,\ j\ge2\\
	(z-1)(1+z^2)s_0+zs_1&=-1 
\end{align*}
or in matrix form
{\scriptsize
\begin{equation*}
	\begin{pmatrix}
		(z-1)(1+z^2)& z& 0&\cdots  &0  \\
		z(z-1) & (1-z)(1+z^2) &-z& \cdots & 0 \\
		&		z(z-1) & (1-z)(1+z^2) &-z & 0 \\
		\vdots  & \ddots  & \ddots & \ddots  & \vdots\\
		0 & 0 & \cdots&z(z-1) &(1-z)(1+z^2)
	\end{pmatrix}
	\begin{pmatrix}
		s_0\\
		s_1\\
		s_2\\
		\vdots\\
		s_K
	\end{pmatrix}=
	\begin{pmatrix}
		-1\\
		0\\
		0\\
		\vdots\\
		0
	\end{pmatrix}
\end{equation*}
}
The quantity $\phi_K=0$ for all $K$, so be concentrate on the remaining two.
{\scriptsize
	\begin{equation*}
		\begin{pmatrix}
			(z-1)(1+z^2)& z& 0&\cdots  &0  \\
			z(z-1) & (1-z)(1+z^2) &-z& \cdots & 0 \\
			&		z(z-1) & (1-z)(1+z^2) &-z & 0 \\
			\vdots  & \ddots  & \ddots & \ddots  & \vdots\\
			0 & 0 & \cdots&z(z-1) &(1-z)(1+z^2)
		\end{pmatrix}
		\begin{pmatrix}
			g_0\\
			g_1\\
			g_2\\
			\vdots\\
			g_K
		\end{pmatrix}=
		\begin{pmatrix}
			-z\\
			0\\
			0\\
			\vdots\\
			0
		\end{pmatrix}
	\end{equation*}
}
and
{\scriptsize
	\begin{equation*}
		\begin{pmatrix}
			z^3+z-1& z& 0&\cdots  &0  \\
			z(z-1) & (1-z)(1+z^2) &-z& \cdots & 0 \\
			&		z(z-1) & (1-z)(1+z^2) &-z & 0 \\
			\vdots  & \ddots  & \ddots & \ddots  & \vdots\\
			0 & 0 & \cdots&z(z-1) &(1-z)(1+z^2)
		\end{pmatrix}
		\begin{pmatrix}
			h_0\\
			h_1\\
			h_2\\
			\vdots\\
			h_K
		\end{pmatrix}=
		\begin{pmatrix}
			-z^3\\
			0\\
			0\\
			\vdots\\
			0
		\end{pmatrix}
	\end{equation*}
}
The usual determinant is
\begin{equation*}
	\mD_{K}=(z-1)\mD_{K-1}-z^3\mD_{K-2}=\frac{(-z)^{K+1}}{W}(r_2^{K+1}-r_1^{K+1}).
\end{equation*}
With $t_K=\mD_{K-1}/\mD_{K}$ and $\tau_K=-z^3t_K$
\begin{equation*}
	t_K=\frac{1}{z-1-z^3t_{K-1}}\ \Longrightarrow \ \tau_K=\frac{z^3}{1-z-\tau_{K-1}},\ \tau_0=0.
\end{equation*}
\begin{equation*}
	\begin{pmatrix}
		z-1& z^2& 0&\cdots  &0  \\
		z & z-1 &z^2& \cdots & 0 \\
		&		z & z-1 &z^2 & 0 \\
		\vdots  & \ddots  & \ddots & \ddots  & \vdots\\
		0 & 0 & \cdots&z &z-1
	\end{pmatrix}
	\begin{pmatrix}
		g_0\\
		g_1\\
		g_2\\
		\vdots\\
		g_K
	\end{pmatrix}=
	\begin{pmatrix}
		-1\\
		z\\
		0\\
		\vdots\\
		0
	\end{pmatrix}
\end{equation*}
\begin{equation*}
	\begin{pmatrix}
		-1+2z-z^2+z^3& z^2(1-z)& 0&\cdots  &0  \\
		z & z-1 &z^2& \cdots & 0 \\
		&		z & z-1 &z^2 & 0 \\
		\vdots  & \ddots  & \ddots & \ddots  & \vdots\\
		0 & 0 & \cdots&z &z-1
	\end{pmatrix}
	\begin{pmatrix}
		h_0\\
		h_1\\
		h_2\\
		\vdots\\
		h_K
	\end{pmatrix}=
	\begin{pmatrix}
		-z^3\\
		0\\
		0\\
		\vdots\\
		0
	\end{pmatrix}
\end{equation*}
From this, we can express the solutions (return to the $x$-axis) in continued fraction form
\begin{align*}
	s_0=\sigma_{K+1}&=\frac{-\mD_K}{(z-1)(1+z^2)\mD_{K}-z^2(z-1)\mD_{K-1}}=
	\frac{-1}{(z-1)(1+z^2)-z^2(z-1)t_{K}}\\&=
	\frac{z}{(1-z)(z+z^3-\tau_K)},\ \sigma_0=\frac1{1-z},\\
	g_0=\gamma_{K+1}&=\frac{-z\mD_K}{(z-1)(1+z^2)\mD_{K}-z^2(z-1)\mD_{K-1}}=z\sigma_{K+1},\\
		h_0=\eta_{K+1}&=\frac{-z^3\mD_{K}}{(-1+2z-z^2+z^3)\mD_{K}-z^3(1-z)\mD_{K-1}}\\&=
	\frac{-z^3}{(-1+2z-z^2+z^3)-z^3(1-z)t_{K}}
	=\frac{z^3}{1-2z+z^2-z^3-(1-z)\tau_{K}},\ \eta_0=0.
\end{align*}

%% file: DUHD.tex
\section{DU and HD are forbidden, A004149}
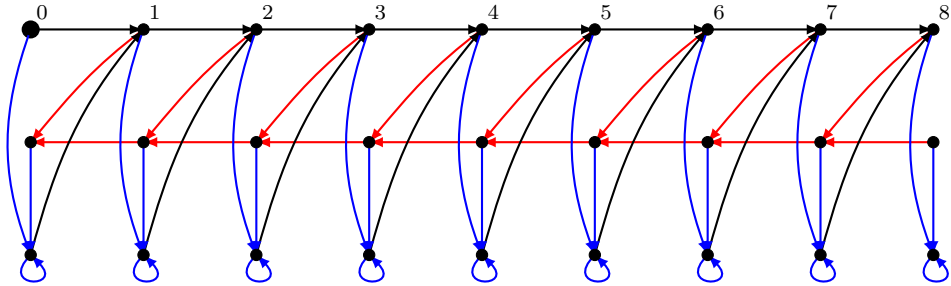
\begin{figure}[h]

	\begin{center}
		\begin{tikzpicture}[scale=1.5,main node/.style={circle,draw,font=\Large\bfseries}]

			\fill (0,0) circle (0.08cm);

			\foreach \x in {0,1,2,3,4,5,6,7}
			{
				\draw[thick, -latex] (\x,0) to  (\x+1,0);   
				      				\draw[thick, latex-,red] (\x,-1) [out=50,in=217]to  (\x+1,0);   
				      				\draw[thick, latex-,red] (\x,-1) to  (\x+1,-1);   
					      				\draw[thick, -latex][out=75,in=-130] (\x,-2) to  (\x+1,0);   
				\node at  (\x+0.1,0.15){\tiny$\x$};
			}

			\foreach \x in {0,1,2,3,4,5,6,7,8}
{
		\draw[thick,-latex,blue ] (\x,-2)  ..  controls (\x-0.3,-0.3-2) and  (\x+0.3,-0.3-2) .. (\x,-2) ;	
	\draw[ thick,blue, latex-] (\x,-2)[in=-110,out=110]to  (\x,0);    
		\draw[ thick,blue, -latex] (\x,-1) to  (\x,-2);

}                        
			
			\node at  (8+0.1,0.15){\tiny$8$};
			
			\foreach \x in {0,1,2,3,4,5,6,7,8}
			{
				\draw (\x,-1) circle (0.05cm);
				\fill (\x,-1) circle (0.05cm);
				\draw (\x,-2) circle (0.05cm);
				\fill (\x,-2) circle (0.05cm);
			}
		\foreach \x in {0,1,2,3,4,5,6,7,8}
		{
			\draw (\x,0) circle (0.05cm);
			\fill (\x,0) circle (0.05cm);
		}
		\end{tikzpicture}
	\end{center}
	\caption{Graph (automaton) to recognize  Motzkin paths with forbidden subwords DU and HD.}
	\label{purpelDUHD}
\end{figure}
The recursion related to Figure \ref{purpelDUHD} are
\begin{align*}
	f_{i+1}&=zf_i+zh_i, \ f_0=1,\\
	g_i&=zf_{i+1}+zg_{i+1},\ i\ge0,\\
	h_{i}&=zf_i+zg_{i}+zh_{i},
\end{align*}
and the system of functional equations is
\begin{align*}
	F(u)&=1+zuF(u)+zuH(u),\\
	G(u)&=\frac zu(G(u)-G(0))+\frac zu(H(u)-H(0)),\\
	H(u)&=zF(u)+zG(u)+z H(u).
\end{align*}
The ubiquitous denominator is $zu^2+(-1+z-z^2-z^3)u+z=z(u-r_1)(u-r_2)$ with
\begin{equation*}
r_1=\frac{1-z+z^2+z^3-W}{2z},\ r_2=\frac{1-z+z^2+z^3-W}{2z},\quad W=\sqrt{(1-z^4)(1-2z-z^2)}.
\end{equation*}
Dividing out the factor $u-r_1$ and simplifying, we end up with
\begin{equation*}
	F(u)=\frac{1}{1-ur_1},\ 	G(u)=\frac{(1-z)r_1-z}{z^2(1-ur_1)},\ 	H(u)=\frac{r_1-z}{z(1-ur_1)},\ 
S(u)=\frac{r_1-z}{z^2(1-ur_1)}.
\end{equation*}
We find
\begin{equation*}
S(0)=1+z+z^2+2z^3+4z^4+8z^5+16z^6+33z^7+69z^8+146z^9+312z^{10}+\cdots
\end{equation*}
which is A004149 in OEIS
and
\begin{equation*}
	S(1)=1+2z+4z^2+9z^3+20z^4+45z^5+102z^6+233z^7+535z^8+1234z^9+2857z^{10}+\cdots
\end{equation*}
1,2,4,9,20,45,102,233
which is A308435  in OEIS. Furthermore
\begin{equation*}
[u^j]S(u)=[u^j]\frac{r_1-z}{z^2(1-ur_1)}=\frac{r_1^{j+1}}{z^2}-\frac{r_1^{j}}{z}.
\end{equation*}
Since, before diving out the factor, we have
\begin{equation*}
	S(u) = \frac{\textsf{numerator}}{zu^2+(-1+z-z^2-z^3)u+z},
\end{equation*}
we get a recursion from it:
\begin{align*}
	zs_{j-2}+(-1+z-z^2-z^3)s_{j-1}+zs_j&=0,\ j\ge2,\\
(-1+z-z^3)s_0+zs_1&=z^2-1,
\end{align*}
which is in matrix form
{\scriptsize
\begin{equation*}
	\begin{pmatrix}
		-1+z-z^3& z& 0&\cdots  &0  \\
		z & -1+z-z^2-z^3&z& \cdots & 0 \\
		&		z & -1+z-z^2-z^3&z & 0 \\
		\vdots  & \ddots  & \ddots & \ddots  & \vdots\\
		0 & 0 & \cdots&z &-1+z-z^2-z^3
	\end{pmatrix}
	\begin{pmatrix}
		s_0\\
		s_1\\
		s_2\\
		\vdots\\
		s_K
	\end{pmatrix}=
	\begin{pmatrix}
		z^2-1\\
		0\\
		0\\
		\vdots\\
		0
	\end{pmatrix}
\end{equation*}}
The usual determinant is
\begin{equation*}
	\mD_{K}=(-1+z-z^2-z^3)\mD_{K-1}-z^2\mD_{K-2}=\frac{(-z)^{K+1}}{W}(r_1^{K+1}-r_2^{K+1}).
\end{equation*}
With $t_K=\mD_{K-1}/\mD_{K}$ and $\tau_K=-z^2t_K$,
\begin{equation*}
	t_K=\frac{1}{-1+z-z^2-z^3-z^2t_{K-1}}\ \Longrightarrow \ \tau_K=\frac{z^2}{1-z+z^2+z^3-\tau_{K-1}},\ \tau_0=0.
\end{equation*}
Since $\phi_K=1$ for all $K$, we concentrate on the remaining two:
{\scriptsize
	\begin{equation*}
		\begin{pmatrix}
			-1+2z-z^2+z^4& z(1-z)& 0&\cdots  &0  \\
			z & -1+z-z^2-z^3&z& \cdots & 0 \\
			&		z & -1+z-z^2-z^3&z & 0 \\
			\vdots  & \ddots  & \ddots & \ddots  & \vdots\\
			0 & 0 & \cdots&z &-1+z-z^2-z^3
		\end{pmatrix}
		\begin{pmatrix}
			g_0\\
			g_1\\
			g_2\\
			\vdots\\
			g_K
		\end{pmatrix}=
		\begin{pmatrix}
			-z^3\\
			0\\
			0\\
			\vdots\\
			0
		\end{pmatrix}
\end{equation*}}
{\scriptsize
	\begin{equation*}
		\begin{pmatrix}
			-1+z-z^3& z& 0&\cdots  &0  \\
			z & -1+z-z^2-z^3&z& \cdots & 0 \\
			&		z & -1+z-z^2-z^3&z & 0 \\
			\vdots  & \ddots  & \ddots & \ddots  & \vdots\\
			0 & 0 & \cdots&z &-1+z-z^2-z^3
		\end{pmatrix}
		\begin{pmatrix}
			h_0\\
			h_1\\
			h_2\\
			\vdots\\
			h_K
		\end{pmatrix}=
		\begin{pmatrix}
			z(z^2-1)\\
			0\\
			0\\
			\vdots\\
			0
		\end{pmatrix}
\end{equation*}}

From this, we can express the solutions (return to the $x$-axis) in continued fraction form
\begin{align*}
	s_0=\sigma_{K+1}&=\frac{(z^2-1)\mD_{K}}{	(-1+z-z^3)\mD_{K}-z^2\mD_{K-1}}
	=\frac{(z^2-1)}{	(-1+z-z^3)-z^2t_{K}}\\&=\frac{1-z^2}{	1-z+z^3-\tau_{K}},\	\ \sigma_0=\frac1{1-z},\\
	g_0=\gamma_{K+1}&=\frac{-z^3\mD_{K}}{(-1+2z-z^2+z^4)\mD_{K}-z^2(1-z)\mD_{K-1}}\\&
	=\frac{z^3}{1-2z+z^2-z^4-(1-z)\tau_{K}}	,\  \gamma_0=0,\\
	h_0=\eta_{K+1}&=\frac{z(z^2-1)\mD_{K}}{(-1+z-z^3)\mD_{K}-z^2\mD_{K-1}}\\&=
	\frac{z(z^2-1)}{(-1+z-z^3)-z^2t_{K}}=\frac{z(1-z^2)}{1-z+z^3-\tau_{K}},\ \eta_0=\frac{z}{1-z}.
\end{align*}

%% file: DUUH.tex
\section{DU and UH are forbidden, A217282}

\begin{figure}[h]

	\begin{center}
		\begin{tikzpicture}[scale=1.5,main node/.style={circle,draw,font=\Large\bfseries}]

			\fill (0,-2) circle (0.08cm);

						\foreach \x in {1,2,3,4,5,6,7}
						{				\draw[thick, -latex] (\x,0) to  (\x+1,0);   
										\node at  (\x+0.1,0.15){\tiny$\x$};}
			\foreach \x in {0,1,2,3,4,5,6,7}
			{

				      				\draw[thick, latex-,red] (\x,-1) [out=50,in=217]to  (\x+1,0);   
				      				\draw[thick, latex-,red] (\x,-1) to  (\x+1,-1);   
					      				\draw[thick, -latex,red] (\x+1,-2) to  (\x,-1);   

			}

			\foreach \x in {0,1,2,3,4,5,6,7,8}
{
		\draw[thick,-latex,blue ] (\x,-2)  ..  controls (\x-0.3,-0.3-2) and  (\x+0.3,-0.3-2) .. (\x,-2) ;	
   
		\draw[ thick,blue, -latex] (\x,-1) to  (\x,-2);

}          

			\foreach \x in {0,1,2,3,4,5,6,7}
              {	\draw[ thick,  -latex] (\x,-2)[in=-110,out=50]to  (\x+1,0); }
			
			\node at  (8+0.1,0.15){\tiny$8$};

			\foreach \x in {0,1,2,3,4,5,6,7,8}
			{
				\draw (\x,-1) circle (0.05cm);
				\fill (\x,-1) circle (0.05cm);
				\draw (\x,-2) circle (0.05cm);
				\fill (\x,-2) circle (0.05cm);
			}
		\foreach \x in {1,2,3,4,5,6,7,8}
		{
			\draw (\x,0) circle (0.05cm);
			\fill (\x,0) circle (0.05cm);
		}
		\end{tikzpicture}
	\end{center}
	\caption{Graph (automaton) to recognize  Motzkin paths with forbidden subwords DU and UH.}
	\label{purpelDUUH}
\end{figure}
The recursion related to Figure \ref{purpelDUUH} are
\begin{align*}
	f_{i+1}&=zf_i+zh_i, \ f_0=0,\\
	g_i&=zf_{i+1}+zg_{i+1}+zh_{i+1},\ i\ge0\\
	h_{i}&=1+zg_{i}+zh_{i}
\end{align*}
and the system of functional equations is
\begin{align*}
	F(u)&=zuF(u)+zuH(u),\\
	G(u)&=\frac zu(F(u)-F(0))+\frac zu(G(u)-G(0))+\frac zu(H(u)-H(0)),\\
	H(u)&=1+zG(u)+z H(u).
\end{align*}
The ubiquitous denominator is $zu^2+(-1+z-z^2-z^3)u+z=z(z-1)(u-r_1)(u-r_2)$ with
\begin{gather*}
r_1=\frac{1-z+z^2-z^3-W}{2z(1-z)},\ r_2=\frac{1-z+z^2-z^3-W}{2z(1-z)},\\ W=\sqrt{(1-z)(1-z-2z^2-2z^3+z^4-z^5)}.
\end{gather*}
Dividing out the factor $u-r_1$ and simplifying, we end up with
\begin{gather*}
	F(u)=\frac{ur_1}{1-u(1-z)r_1},\ G(u)=\frac{(1-z)r_1-z}{z^2(1-u(1-z)r_1)},\ H(u)=\frac{(1-zu)r_1}{z(1-u(1-z)r_1)},\\
	 S(u)=\frac{r_1-z}{z^2(1-u(1-z)r_1)}.
\end{gather*}
We find
\begin{equation*}
S(0)=1+z+2z^2+3z^3+5z^4+9z^5+16z^6+30z^7+57z^8+110z^9+216z^{10}+\cdots
\end{equation*}
which is A217282 in OEIS
and
\begin{equation*}
	S(1)=1+2z+4z^2+7z^3+13z^4+25z^5+49z^6+98z^7+198z^8+404z^9+831z^{10}+\cdots
\end{equation*}
which is not  in OEIS. Furthermore
\begin{equation*}
[u^j]S(u)=[u^j]\frac{r_1-z}{z^2(1-u(1-z)r_1)}=
\frac{(1-z)^jr_1^{j+1}}{z^2}-\frac{(1-z)^jr_1^{j}}{z}
,\quad j\ge1.
\end{equation*}
Since, before diving out the factor, we have
\begin{equation*}
	S(u) = \frac{\textsf{numerator}}{-z-(z-1)(1+z^2)u+z(z-1)u^2},
\end{equation*}
we get a recursion from it:
\begin{align*}
	z(z-1)s_{j-2}+(1-z+z^2-z^3)s_{j-1}-zs_j&=0,\ j\ge2,\\
(z-1)s_0+zs_1&=-1,
\end{align*}
which is in matrix form
{\footnotesize
	\begin{equation*}
		\begin{pmatrix}
			z-1& z& 0&\cdots  &0  \\
			z(z-1) & 1-z+z^2-z^3&-z& \cdots & 0 \\
			&		z(z-1) & 1-z+z^2-z^3&-z & 0 \\
			\vdots  & \ddots  & \ddots & \ddots  & \vdots\\
			0 & 0 & \cdots&z(z-1) &1-z+z^2-z^3
		\end{pmatrix}
		\begin{pmatrix}
			s_0\\
			s_1\\
			s_2\\
			\vdots\\
			s_K
		\end{pmatrix}=
		\begin{pmatrix}
			-1\\
			0\\
			0\\
			\vdots\\
			0
		\end{pmatrix}
\end{equation*}}
Likewise,
{\footnotesize
	\begin{equation*}
		\begin{pmatrix}
			z^3+z-1& z& 0&\cdots  &0  \\
			z(z-1) & 1-z+z^2-z^3&-z& \cdots & 0 \\
			&		z(z-1) & 1-z+z^2-z^3&-z & 0 \\
			\vdots  & \ddots  & \ddots & \ddots  & \vdots\\
			0 & 0 & \cdots&z(z-1) &1-z+z^2-z^3
		\end{pmatrix}
		\begin{pmatrix}
			g_0\\
			g_1\\
			g_2\\
			\vdots\\
			g_K
		\end{pmatrix}=
		\begin{pmatrix}
			-z^2\\
			0\\
			0\\
			\vdots\\
			0
		\end{pmatrix}
\end{equation*}}
and
{\footnotesize
	\begin{equation*}
		\begin{pmatrix}
			z^3+z-1& z& 0&\cdots  &0  \\
			z(z-1) & 1-z+z^2-z^3&-z& \cdots & 0 \\
			&		z(z-1) & 1-z+z^2-z^3&-z & 0 \\
			\vdots  & \ddots  & \ddots & \ddots  & \vdots\\
			0 & 0 & \cdots&z(z-1) &1-z+z^2-z^3
		\end{pmatrix}
		\begin{pmatrix}
			h_0\\
			h_1\\
			h_2\\
			\vdots\\
			h_K
		\end{pmatrix}=
		\begin{pmatrix}
			-1\\
			-z\\
			0\\
			\vdots\\
			0
		\end{pmatrix}
\end{equation*}}
The usual determinant is
\begin{equation*}
	\mD_{K}=(1-z+z^2-z^3)\mD_{K-1}+z^2(z-1)\mD_{K-2}=\frac{(z(1-z))^{K+1}}{W}(r_2^{K+1}-r_1^{K+1}).
\end{equation*}
With $t_K=\mD_{K-1}/\mD_{K}$ and $\tau_K=z^2(1-z)t_K$,
\begin{equation*}
	t_K=\frac{1}{1-z+z^2-z^3+z^2(z-1)t_{K-1}}\ \Longrightarrow \ \tau_K=\frac{z^2(1-z)}{1-z+z^2-z^3-\tau_{K-1}},\ \tau_0=0.
\end{equation*}
From this, we can express the solutions (return to the $x$-axis) in continued fraction form
\begin{align*}
	s_0=\sigma_{K+1}&=\frac{-\mD_{K}}{	(z-1)\mD_{K}-z^2(z-1)\mD_{K-1}}
	=\frac{1}{	1-z-\tau_{K}},\	\ \sigma_0=\frac1{1-z},\\
	g_0=\gamma_{K+1}&=\frac{-z^2\mD_{K}}{(z^3+z-1)\mD_{K}-z^2(z-1)\mD_{K-1}}
	=\frac{z^2}{1-z-z^3-\tau_{K}},\  \gamma_0=0,\\
	h_0=\eta_{K+1}&=\frac{-\mD_{K}+z^2\mD_{K-1}}{(z^3+z-1)\mD_{K}-z^2(z-1)\mD_{K-1}}=
	\frac{-1+z^2t_{K}}{(z^3+z-1)-z^2(z-1)t_{K}}\\&
	=\frac1{1-z}+\frac{z^3}{(1-z)(1-z-z^3-\tau_K)},\ \eta_0=\frac{z}{1-z}.
\end{align*}

%% file: HUDD.tex
\section{HU and DD are forbidden, A217282}
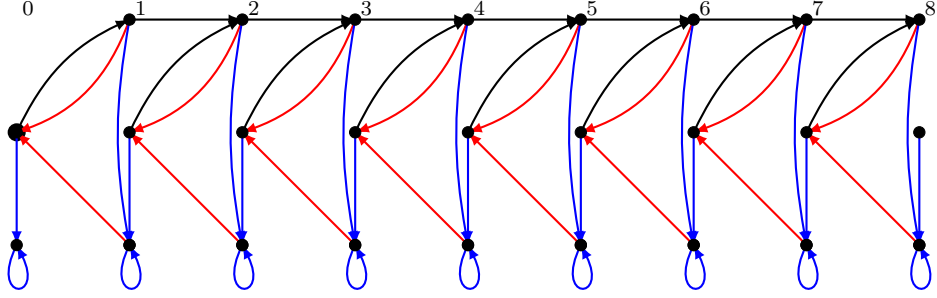
\begin{figure}[h]

	\begin{center}
		\begin{tikzpicture}[scale=1.5,main node/.style={circle,draw,font=\Large\bfseries}]

			\fill (0,0) circle (0.08cm);

			\foreach \x in {0,...,8}
			{
				\draw[thick,-latex,blue ] (\x,-1)  ..  controls (\x-0.25,-1-0.5) and  (\x+0.25,-1-0.5) .. (\x,-1+0) ;	
			}
			
			\foreach \x in {0,...,7}
			{
				\draw[thick,-latex, ] (\x,0) [in=-155,out=65] to (\x+1,1+0) ;	
				
				 				 \draw[thick,  -latex,red] (\x+1,1) [out=-110,in=20]to   (\x,0);
			}
			\foreach \x in {1,...,7}
			{ \draw[thick,  -latex] (\x,1)  to   (\x+1,1);}

			\foreach \x in {0,...,7}
			{
				
			}

			\foreach \x in {0,1,2,3,4,5,6,7}
			{
				
				\node at  (\x+0.1,1.1){\tiny$\x$};
			}			
			\foreach \x in {0,1,2,3,4,5,6,7}
			{
				
				\draw[thick,red, latex-] (\x,0) to  (\x+1,-1);	
			}	
		
		\foreach \x in {1,2,3,4,5,6,7,8}
		{\draw[thick, -latex,blue] (\x,1) [out=-100, in=100]to  (\x,-1);	}
			\foreach \x in {0,1,2,3,4,5,6,7,8}
			{
				
								\draw[thick, -latex,blue] (\x,0)  to  (\x,-1);	
			}

			\node at  (8+0.1,1.1){\tiny$8$};
			
			\foreach \x in {0,1,2,3,4,5,6,7,8}
			{
				\draw (\x,0) circle (0.05cm);
				\fill (\x,0) circle (0.05cm);
					\draw (\x,-1) circle (0.05cm);
					\fill (\x,-1) circle (0.05cm);
				
			}
			\foreach \x in {1,2,3,4,5,6,7,8}
			{
				\draw (\x,1) circle (0.05cm);
				\fill (\x,1) circle (0.05cm);
				\draw (\x,-1) circle (0.05cm);
				\fill (\x,-1) circle (0.05cm);
				
			}
		\end{tikzpicture}
	\end{center}
	\caption{Graph (automaton) to recognize  Motzkin paths with forbidden subwords HU and DD.}
\label{purpelHUDD}
\end{figure}
The following recursions can be read off the automaton (Figure \ref{purpelHUDD}), by considering the last step separately:
\begin{align*}
	f_{i+1}&=zf_{i}+zg_{i},\ f_0=0,\\
	g_{i}&=[i=0]+zf_{i+1}+zh_{i+1},\\
	h_{i}&=zf_{i}+zg_{i}+zh_{i}, 
\end{align*}
and using double generating functions
\begin{align*}
	F(u)&=zuF(u)+zuG(u),\\
	G(u)&=1+\frac zu\big(F(u)-F(0)\big)+\frac zu\big(H(u)-H(0)\big),\\
	H(u)&=zF(u)+zG(u)+zH(u).
\end{align*}
We have the factorization $-z^2+(z+1)(1-z)^2u+z(z-1)u^2=z(z-1)(u-r_1)(u-r_2)$
with 
\begin{gather*}
	r_1={\frac {1-z-z^2+z^{3}-W}{2z(1-z) }}, \quad
	r_2={\frac {1-z-z^2+z^{3}+W}{2z(1-z) }},\\ W=\sqrt{(1-z)(1-z-2z^2-2z^3+z^4-z^5)}.
\end{gather*}
Solving the system of functional equations and dividing out the factor $u-r_1$ we get after simplification
\begin{gather*}
	F(u)=\frac{u(1-z)r_1}{z-u(1-z)r_1}, \ G(u)=\frac{(1-z)(1-zu)r_1}{z(z-u(1-z)r_1)},\ H(u)=\frac{r_1}{z-u(1-z)r_1},\\ S(u)=\frac{r_1}{z(z-u(1-z)r_1)}.
\end{gather*}
From this we find 
\begin{equation*}
S(0)=1+z+2z^2+3z^3+5z^4+9z^5+16z^6+30z^7+57z^8+110z^9+216z^{ 10} +\cdots,
\end{equation*}
which is A217282. The generating function for meanders is
\begin{equation*}
	S(1)=1+2 z+4 z^2+8 z^3+16 z^4+33 z^5+69 z^6+146 z^7+312 z^8+671 z^9+1451 z^{10}+\cdots,
\end{equation*}
which is not in OEIS. Furthermore, for $j\ge1$,
\begin{equation*}
[u^j]S(u)=[u^j]\frac{r_1}{z^2(1-\frac{u(1-z)r_1}{z})}=\frac{(1-z)^jr_1^{j+1}}{z^{j+2}},
\end{equation*}
which is the generating function of restricted Motzkin paths ending at level $j$.

For $S=F+G+H$ we get
\begin{equation*}
	S=\frac {\textsf{numerator}}{-z^2+(z+1)(1-z)^2u+z(z-1)u^2},
\end{equation*}
and from this a recursion
 \begin{align*}
	z(z-1)s_{j-2}+(z+1)(1-z)^2s_{j-1}-z^2s_{j}&=0,\quad j\ge2,\\
-(z+1)(1-z)^2s[0]+z^2s[1]&=-1,
\end{align*}
which is in matrix form
{\scriptsize
	\begin{equation*}
	\begin{pmatrix}
		-(z+1)(1-z)^2 &z^2 & 0&\cdots & &0  \\
		z(z-1) & (z+1)(1-z)^2 &-z^2& \cdots && 0 \\
		&		z(z-1) & (z+1)(1-z)^2 &-z^2& & 0 \\
		\vdots  & \ddots  & \ddots & \ddots  && \vdots\\
		0 & 0 & \cdots&&z(z-1) &(z+1)(1-z)^2
	\end{pmatrix}
	\begin{pmatrix}
		s_0\\
		s_1\\
		s_2\\
		\vdots\\
		s_K
	\end{pmatrix}=
	\begin{pmatrix}
		-1\\
		0\\
		0\\
		\vdots\\
		0
	\end{pmatrix}
\end{equation*}
}
We have $f_0=\phi_K=0$ for all $K$. Furthermore,
{\scriptsize
	\begin{equation*}	
		\begin{pmatrix}
			z^2+z-1&z^2 & 0&\cdots & &0  \\
			z(z-1) & (z+1)(1-z)^2 &-z^2& \cdots && 0 \\
			&		z(z-1) & (z+1)(1-z)^2 &-z^2& & 0 \\
			\vdots  & \ddots  & \ddots & \ddots  && \vdots\\
			0 & 0 & \cdots&&z(z-1) &(z+1)(1-z)^2
		\end{pmatrix}
		\begin{pmatrix}
			g_0\\
			g_1\\
			g_2\\
			\vdots\\
			g_K
		\end{pmatrix}=
		\begin{pmatrix}
			z-1\\
			z(z-1)\\
			0\\
			\vdots\\
			0
		\end{pmatrix}
	\end{equation*}
}
and
{\scriptsize
	\begin{equation*}	
		\begin{pmatrix}
			-(z+1)(-1+z)^2			&z^2 & 0&\cdots & &0  \\
			z(z-1) & (z+1)(1-z)^2 &-z^2& \cdots && 0 \\
			&		z(z-1) & (z+1)(1-z)^2 &-z^2& & 0 \\
			\vdots  & \ddots  & \ddots & \ddots  && \vdots\\
			0 & 0 & \cdots&&z(z-1) &(z+1)(1-z)^2
		\end{pmatrix}
		\begin{pmatrix}
			h_0\\
			h_1\\
			h_2\\
			\vdots\\
			h_K
		\end{pmatrix}=
		\begin{pmatrix}
			-z\\
			0\\
			0\\
			\vdots\\
			0
		\end{pmatrix}
	\end{equation*}
}
Let $\mathcal{D}_K$ be the determinant with $K$ rows and columns, after deleting the first row and column. Expanding, we derive the recursion
($\mathcal{D}_0=1,\ \mathcal{D}_1= (z+1)(1-z)^2$)
\begin{equation*}
	\mathcal{D}_K=(z+1)(1-z)^2\mathcal{D}_{K-1}+z^3(z-1)\mathcal{D}_{K-2}
	=\frac{z^{K+1}(1-z)^{K+1}}{W}(r_2^{K+1}-r_1^{K+1}).
\end{equation*}
Rewriting this, we get, with $t_K=\frac{\mathcal{D}_{K-1}}{\mathcal{D}_{K}}$, $\tau_K =z^3(1-z)t_K$, $\tau_0=0$
\begin{equation*}
t_K=\frac{1}{1-z-z^2+z^3+z^3(z-1)t_{K-1}}, \quad \tau_K=\frac{z^3(1-z)}{1-z+z^3-\tau_{K-1}}.
\end{equation*}
Solving the system by Cramer's rule, we get  the solution for $s_0$ (return to the $x$-axis) as
\begin{align*}
	s_0=\sigma_{K+1}&=\frac{-\mD_K }{(-1+z+z^2-z^3)\mD_K -z^3(z-1)\mD_{K-1}}\\&=
	\frac{-1}{(-1+z+z^2-z^3) -z^3(z-1)t_{K}}=	\frac{1}{1-z-z^2+z^3 -\tau_{K}},\ \sigma_0=\frac1{1-z},\\
	g_0=\gamma_{K+1}&=\frac{(z-1)\mD_K-z^3(z-1)\mD_{K-1} }{(-1+z+z^2)\mD_K -z^3(z-1)\mD_{K-1}}=
	\frac{(z-1)-z^3(z-1)t_{K} }{(-1+z+z^2) -z^3(z-1)t_{K}}\\&=
	1+\frac{z^2}{1-z-z^2-\tau_K},\ \gamma_0 =1,\\
	h_0=\eta_{K+1}&=\frac{-z\mD_K}{(z^2-z^3-1+z)\mD_K -z^3(z-1)\mD_{K-1}}\\&=
	\frac{-z}{(z^2-z^3-1+z) -z^3(z-1)t_{K}}=	\frac{z}{1-z-z^2+z^3-\tau_{K}},\ \eta_0=\frac{z}{1-z}.
	\end{align*}

%% file: roitner2.bbl
\begin{thebibliography}{1}
	
	
 
	
	
	
	\bibitem{ABBG}
	Andrei Asinowski, Axel Bacher, Cyril Banderier, and Bernhard Gittenberger.
	\newblock Analytic combinatorics of lattice paths with
	forbidden patterns, the vectorial kernel method, and generating functions for pushdown automata.
	\newblock {\em Algorithmica}, 82(3): 386-428 2020.
	
	\bibitem{AsyWally}
	Andrei Asinowski, Cyril Banderier, and Valerie Roitner.
	\newblock Generating functions for lattice paths with several forbidden
	patterns.
	\newblock {\em S\'{e}m. Lothar. Combin.}, 84B:Art. 95, 12, 2020.
	
	\bibitem{BF}
	Cyril Banderier and Philippe Flajolet.
	\newblock Basic analytic combinatorics of directed lattice paths.
	\newblock {\em Theoret. Comput. Sci.}, 281(1-2):37--80, 2002.
	\newblock Selected papers in honour of Maurice Nivat.
	
	\bibitem{Prodinger-kernel}
	Helmut Prodinger.
	\newblock The kernel method: a collection of examples.
	\newblock {\em S\'{e}m. Lothar. Combin.}, 50: Art. B50f, 19, 2003/04.
	
	\bibitem{garden}
	Helmut Prodinger.
	\newblock A walk in my lattice path garden.
	\newblock {\em S\'{e}m. Lothar. Combin.}, 87B: Art. 1, 49, 2023.
	
	\bibitem{Vally}
	Valerie Roitner.
	\newblock {\em Studies on several parameters in lattice paths}.
	\newblock PhD thesis, Dissertation, Technische Universit\"at Wien, 2020.
	
	\bibitem{Drew}
	Andrew~V. Sills.
	\newblock Finite {R}ogers-{R}amanujan type identities.
	\newblock {\em Electron. J. Combin.}, 10:Research Paper 13, 122, 2003.
	
\end{thebibliography}
